\documentclass[a4paper,11pt,reqno]{amsart}
\newtheorem{lem}{Lemma}
\newtheorem{theo}{Theorem}
\newtheorem{rem}{Remark}

\numberwithin{equation}{section} 
 
\newcommand{\sgn}{\operatorname{sgn}} 

\newcommand{\RhombusVert}{\rlvec(0.8660254037844 0.5)
                          \rlvec(0.8660254037844 -0.5)
                          \rlvec(-0.8660254037844 -0.5)
                          \rlvec(-0.8660254037844  0.5)}
\newcommand{\RhombusPos}{\rlvec(0.8660254037844 0.5) 
                         \rlvec(0 -1)
                         \rlvec(-0.8660254037844 -0.5)
                         \rlvec(0 1)}   
\newcommand{\RhombusNeg}{\rlvec(0 1) 
                         \rlvec(0.8660254037844 -0.5)
                         \rlvec(0 -1)
                         \rlvec(-0.8660254037844 0.5)}  
\newcommand{\DreieckBreit}{\rlvec(0 1)
                           \rlvec(0.8660254037844 -0.5)
                           \rlvec(-0.8660254037844 -0.5)}
\newcommand{\DreieckSpitz}{\rlvec(0.8660254037844 0.5)
                           \rlvec(0 -1)
                           \rlvec(-0.8660254037844 0.5)}  
\RequirePackage{graphics}

\def\t@xderror #1{%
   \GenericError{%
      \space\space\space\@spaces\@spaces\@spaces
   }{%
      TeXdraw Error: #1%
   }{%
      See the TeXdraw manual for an explanation.%
   }{\@ehc}%
}

\input texdraw
\catcode`\@=11

\thinlines
\newskip\Einheit \Einheit=0.6cm
\newcount\xcoord \newcount\ycoord
\newdimen\xdim \newdimen\ydim \newdimen\PfadD@cke \newdimen\Pfadd@cke
\def\PfadDicke#1{\PfadD@cke#1 \divide\PfadD@cke by2 \Pfadd@cke\PfadD@cke \multiply\PfadD@cke by2}
\long\def\LOOP#1\REPEAT{\def\BODY{#1}\ITERATE}
\def\ITERATE{\BODY \let\next\ITERATE \else\let\next\relax\fi \next}
\let\REPEAT=\fi
\def\Punkt{\hbox{\raise-2pt\hbox to0pt{\hss\scriptsize$\bullet$\hss}}}
\def\DuennPunkt(#1,#2){\unskip
  \raise#2 \Einheit\hbox to0pt{\hskip#1 \Einheit
          \raise-2.5pt\hbox to0pt{\hss\normalsize$\bullet$\hss}\hss}}
\def\NormalPunkt(#1,#2){\unskip
  \raise#2 \Einheit\hbox to0pt{\hskip#1 \Einheit
          \raise-3pt\hbox to0pt{\hss\large$\bullet$\hss}\hss}}
\def\DickPunkt(#1,#2){\unskip
  \raise#2 \Einheit\hbox to0pt{\hskip#1 \Einheit
          \raise-4pt\hbox to0pt{\hss\Large$\bullet$\hss}\hss}}
\def\Kreis(#1,#2){\unskip
  \raise#2 \Einheit\hbox to0pt{\hskip#1 \Einheit
          \raise-4pt\hbox to0pt{\hss\Large$\circ$\hss}\hss}}
\def\Diagonale(#1,#2)#3{\unskip\leavevmode
  \xcoord#1\relax \ycoord#2\relax
      \raise\ycoord \Einheit\hbox to0pt{\hskip\xcoord \Einheit
         \unitlength\Einheit
         \line(1,1){#3}\hss}}
\def\AntiDiagonale(#1,#2)#3{\unskip\leavevmode
  \xcoord#1\relax \ycoord#2\relax \advance\xcoord by -0.05\relax
      \raise\ycoord \Einheit\hbox to0pt{\hskip\xcoord \Einheit
         \unitlength\Einheit
         \line(1,-1){#3}\hss}}
\def\Pfad(#1,#2),#3\endPfad{\unskip\leavevmode
  \xcoord#1 \ycoord#2 \thicklines\ZeichnePfad#3\endPfad\thinlines}
\def\ZeichnePfad#1{\ifx#1\endPfad\let\next\relax
  \else\let\next\ZeichnePfad
    \ifnum#1=1
      \raise\ycoord \Einheit\hbox to0pt{\hskip\xcoord \Einheit
         \vrule height\Pfadd@cke width1 \Einheit depth\Pfadd@cke\hss}%
      \advance\xcoord by 1
    \else\ifnum#1=2  \advance\ycoord by -1
       \raise\ycoord
             \Einheit\hbox to0pt{\hskip\xcoord \Einheit 
        \hbox{\hskip-\PfadD@cke  \vrule height1 \Einheit width\PfadD@cke depth0pt}\hss}%
    \else\ifnum#1=3
      \raise\ycoord \Einheit\hbox to0pt{\hskip\xcoord \Einheit
         \unitlength\Einheit
         \line(1,1){1}\hss}
      \advance\xcoord by 1
      \advance\ycoord by 1
    \else\ifnum#1=4
      \raise\ycoord \Einheit\hbox to0pt{\hskip\xcoord \Einheit
         \unitlength\Einheit
         \line(1,-1){1}\hss}
      \advance\xcoord by 1
      \advance\ycoord by -1
    \fi\fi\fi\fi
  \fi\next}
\def\hSSchritt{\leavevmode\raise-.4pt\hbox to0pt{\hss.\hss}\hskip.2\Einheit
  \raise-.4pt\hbox to0pt{\hss.\hss}\hskip.2\Einheit
  \raise-.4pt\hbox to0pt{\hss.\hss}\hskip.2\Einheit
  \raise-.4pt\hbox to0pt{\hss.\hss}\hskip.2\Einheit
  \raise-.4pt\hbox to0pt{\hss.\hss}\hskip.2\Einheit}
\def\vSSchritt{\vbox{\baselineskip.2\Einheit\lineskiplimit0pt
\hbox{.}\hbox{.}\hbox{.}\hbox{.}\hbox{.}}}
\def\DSSchritt{\leavevmode\raise-.4pt\hbox to0pt{%
  \hbox to0pt{\hss.\hss}\hskip.2\Einheit
  \raise.2\Einheit\hbox to0pt{\hss.\hss}\hskip.2\Einheit
  \raise.4\Einheit\hbox to0pt{\hss.\hss}\hskip.2\Einheit
  \raise.6\Einheit\hbox to0pt{\hss.\hss}\hskip.2\Einheit
  \raise.8\Einheit\hbox to0pt{\hss.\hss}\hss}}
\def\dSSchritt{\leavevmode\raise-.4pt\hbox to0pt{%
  \hbox to0pt{\hss.\hss}\hskip.2\Einheit
  \raise-.2\Einheit\hbox to0pt{\hss.\hss}\hskip.2\Einheit
  \raise-.4\Einheit\hbox to0pt{\hss.\hss}\hskip.2\Einheit
  \raise-.6\Einheit\hbox to0pt{\hss.\hss}\hskip.2\Einheit
  \raise-.8\Einheit\hbox to0pt{\hss.\hss}\hss}}
\def\SPfad(#1,#2),#3\endSPfad{\unskip\leavevmode
  \xcoord#1 \ycoord#2 \ZeichneSPfad#3\endSPfad}
\def\ZeichneSPfad#1{\ifx#1\endSPfad\let\next\relax
  \else\let\next\ZeichneSPfad
    \ifnum#1=1
      \raise\ycoord \Einheit\hbox to0pt{\hskip\xcoord \Einheit
         \hSSchritt\hss}%
      \advance\xcoord by 1
    \else\ifnum#1=2   \advance\ycoord by -1
      \raise\ycoord \Einheit\hbox to0pt{\hskip\xcoord \Einheit
        \hbox{\hskip-2pt \vSSchritt}\hss}%
    \else\ifnum#1=3
      \raise\ycoord \Einheit\hbox to0pt{\hskip\xcoord \Einheit
         \DSSchritt\hss}
      \advance\xcoord by 1
      \advance\ycoord by 1
    \else\ifnum#1=4
      \raise\ycoord \Einheit\hbox to0pt{\hskip\xcoord \Einheit
         \dSSchritt\hss}
      \advance\xcoord by 1
      \advance\ycoord by -1
    \fi\fi\fi\fi
  \fi\next}
\def\Koordinatenachsen(#1,#2){\unskip
 \hbox to0pt{\hskip-.5pt\vrule height#2 \Einheit width.5pt depth1 \Einheit}%
 \hbox to0pt{\hskip-1 \Einheit \xcoord#1 \advance\xcoord by1
    \vrule height0.25pt width\xcoord \Einheit depth0.25pt\hss}}
\def\Koordinatenachsen(#1,#2)(#3,#4){\unskip
 \hbox to0pt{\hskip-.5pt \ycoord-#4 \advance\ycoord by1
    \vrule height#2 \Einheit width.5pt depth\ycoord \Einheit}%
 \hbox to0pt{\hskip-1 \Einheit \hskip#3\Einheit 
    \xcoord#1 \advance\xcoord by1 \advance\xcoord by-#3 
    \vrule height0.25pt width\xcoord \Einheit depth0.25pt\hss}}
\def\Gitter(#1,#2){\unskip \xcoord0 \ycoord0 \leavevmode
  \LOOP\ifnum\ycoord<#2
    \loop\ifnum\xcoord<#1
      \raise\ycoord \Einheit\hbox to0pt{\hskip\xcoord \Einheit\Punkt\hss}%
      \advance\xcoord by1
    \repeat
    \xcoord0
    \advance\ycoord by1
  \REPEAT}
\def\Gitter(#1,#2)(#3,#4){\unskip \xcoord#3 \ycoord#4 \leavevmode
  \LOOP\ifnum\ycoord<#2
    \loop\ifnum\xcoord<#1
      \raise\ycoord \Einheit\hbox to0pt{\hskip\xcoord \Einheit\Punkt\hss}%
      \advance\xcoord by1
    \repeat
    \xcoord#3
    \advance\ycoord by1
  \REPEAT}
\def\Label#1#2(#3,#4){\unskip \xdim#3 \Einheit \ydim#4 \Einheit
  \def\lo{\advance\xdim by-.5 \Einheit \advance\ydim by.5 \Einheit}%
  \def\llo{\advance\xdim by-.25cm \advance\ydim by.5 \Einheit}%
  \def\loo{\advance\xdim by-.5 \Einheit \advance\ydim by.25cm}%
  \def\o{\advance\ydim by.25cm}%
  \def\ro{\advance\xdim by.5 \Einheit \advance\ydim by.5 \Einheit}%
  \def\rro{\advance\xdim by.25cm \advance\ydim by.5 \Einheit}%
  \def\roo{\advance\xdim by.5 \Einheit \advance\ydim by.25cm}%
  \def\l{\advance\xdim by-.30cm}%
  \def\r{\advance\xdim by.30cm}%
  \def\lu{\advance\xdim by-.5 \Einheit \advance\ydim by-.6 \Einheit}%
  \def\llu{\advance\xdim by-.25cm \advance\ydim by-.6 \Einheit}%
  \def\luu{\advance\xdim by-.5 \Einheit \advance\ydim by-.30cm}%
  \def\u{\advance\ydim by-.30cm}%
  \def\ru{\advance\xdim by.5 \Einheit \advance\ydim by-.6 \Einheit}%
  \def\rru{\advance\xdim by.25cm \advance\ydim by-.6 \Einheit}%
  \def\ruu{\advance\xdim by.5 \Einheit \advance\ydim by-.30cm}%
  #1\raise\ydim\hbox to0pt{\hskip\xdim
     \vbox to0pt{\vss\hbox to0pt{\hss$#2$\hss}\vss}\hss}%
}
\catcode`\@=12

\begin{document}

\author[Ilse Fischer]{\box\Adr}

\newbox\Adr
\setbox\Adr\vbox{
\centerline{ \Large Ilse Fischer$^\dag$  \thanks{$^\dag$Partially supported by the Austrian Science Foundation FWF, grant P13190-MAT}}
\vspace{0.3cm}
\centerline{Institut f\"ur Mathematik der Universit\"at Klagenfurt,}
\centerline{Universit\"atsstra{\ss}e 65-67, A-9020 Klagenfurt, Austria.}
\centerline{E-mail: {\tt Ilse.Fischer@uni-klu.ac.at}}
}

\

\title[Enumeration of rhombus tilings]{Enumeration of rhombus tilings of a hexagon which contain a 
fixed rhombus in the centre}
\date{}

\begin{abstract}
We compute the number of rhombus tilings of a hexagon with side lengths 
$a$,$b$,$c$,$a$,$b$,$c$ which contain the central rhombus and the 
number of rhombus tilings of a hexagon with side lengths $a$,$b$,$c$,$a$,$b$,$c$ which 
contain the `almost central' rhombus above the centre.
\end{abstract}

\maketitle
{\parindent0cm
{\small {\it  Mathematics Subject Classification}: 05A15; 05A16, 05A19, 05B45, 33C20, 52C20}

 }

\section{Introduction}
\label{int}
Let $a$, $b$ and $c$ be positive integers and consider a hexagon with side lengths 
$a$,$b$,$c$,$a$,$b$,$c$ whose angles are $120^\circ$ (see Figure~\ref{fig:1}). 

The subject of our interest is rhombus tilings of such a
hexagon using rhombi with all sides of length $1$ and angles $60^\circ$ and $120^\circ$. 
Figure~\ref{fig:2} shows an example of a rhombus tiling of a hexagon with $a=3$, $b=5$ 
and $c=4$.

A first natural question to be asked is how many rhombus tiling of a fixed hexagon exist.
A well known bijection between such rhombus tilings and plane partitions contained in an 
$a \times b \times c$ box \cite{rhombustilings-planepartitions} 
and MacMahon's enumeration of plane partitions
\cite[Sec.~429, $q\rightarrow 1$; proof in Sec.~494]{MacMahon}  give the following answer: 
{\em The 
number of all rhombus tilings of a hexagon with side lengths 
$a$,$b$,$c$,$a$,$b$,$c$ equals }
\begin{equation}
\label{mac}
\prod_{i=1}^{a} \prod_{j=1}^{b} \prod_{k=1}^{c} \frac{i+j+k-1}{i+j+k-2}.
\end{equation}

On a next level, one may ask for the number of rhombus tilings with special properties. 
In this paper, we address this question. We study rhombus tilings  of a hexagon
which contain the central rhombus and  rhombus tilings of a hexagon which 
contain the `almost central' rhombus above the centre. By the `central rhombus' we mean the 
rhombus whose centre is equal to the centre of the hexagon (see Figure~\ref{fig:3},
where the central rhombus is marked; furthermore, the example of a rhombus tiling in 
Figure~\ref{fig:2} contains the central rhombus). By the `almost central' rhombus above the centre
we mean the horizontal rhombus whose lowest vertex is the centre of the hexagon (see Figure~\ref{fig:20}).

%Sechseck in Dreiecke unterteilt;
\begin{figure}[hp]
\centertexdraw{
\drawdim cm
\DreieckBreit
\rmove(0 1)
\DreieckSpitz
\rmove(0.8660254037844 -0.5)
\DreieckBreit
\rmove(0 1)
\DreieckSpitz
\rmove(0.8660254037844 -0.5)
\DreieckBreit
\rmove(0 1)
\DreieckSpitz
\rmove(0.8660254037844 -0.5)
\DreieckBreit
\rmove(-0.8660254037844 -0.5)
\rmove(-0.8660254037844 -0.5)
\rmove(-0.8660254037844 -0.5)
\rmove(0 -1)
\DreieckBreit
\rmove(0 1)
\DreieckSpitz
\rmove(0.8660254037844 -0.5)
\DreieckBreit
\rmove(0 1)
\DreieckSpitz
\rmove(0.8660254037844 -0.5)
\DreieckBreit
\rmove(0 1)
\DreieckSpitz
\rmove(0.8660254037844 -0.5)
\DreieckBreit
\rmove(0 1)
\DreieckSpitz
\rmove(0.8660254037844 -0.5)
\DreieckBreit
\rmove(-0.8660254037844 -0.5)
\rmove(-0.8660254037844 -0.5)
\rmove(-0.8660254037844 -0.5)
\rmove(-0.8660254037844 -0.5)
\rmove(0 -1)
\DreieckBreit
\rmove(0 1)
\DreieckSpitz
\rmove(0.8660254037844 -0.5)
\DreieckBreit
\rmove(0 1)
\DreieckSpitz
\rmove(0.8660254037844 -0.5)
\DreieckBreit
\rmove(0 1)
\DreieckSpitz
\rmove(0.8660254037844 -0.5)
\DreieckBreit
\rmove(0 1)
\DreieckSpitz
\rmove(0.8660254037844 -0.5)
\DreieckBreit
\rmove(0 1)
\DreieckSpitz
\rmove(0.8660254037844 -0.5)
\DreieckBreit
\rmove(-0.8660254037844 -0.5)
\rmove(-0.8660254037844 -0.5)
\rmove(-0.8660254037844 -0.5)
\rmove(-0.8660254037844 -0.5)
\rmove(-0.8660254037844 -0.5)
\rmove(0 -1)
\DreieckBreit
\rmove(0 1)
\DreieckSpitz
\rmove(0.8660254037844 -0.5)
\DreieckBreit
\rmove(0 1)
\DreieckSpitz
\rmove(0.8660254037844 -0.5)
\DreieckBreit
\rmove(0 1)
\DreieckSpitz
\rmove(0.8660254037844 -0.5)
\DreieckBreit
\rmove(0 1)
\DreieckSpitz
\rmove(0.8660254037844 -0.5)
\DreieckBreit
\rmove(0 1)
\DreieckSpitz
\rmove(0.8660254037844 -0.5)
\DreieckBreit
\rmove(0 1)
\DreieckSpitz
\rmove(0.8660254037844 -0.5)
\DreieckBreit
\rmove(-0.8660254037844 -0.5)
\rmove(-0.8660254037844 -0.5)
\rmove(-0.8660254037844 -0.5)
\rmove(-0.8660254037844 -0.5)
\rmove(-0.8660254037844 -0.5)
\rmove(-0.8660254037844 -0.5)
\DreieckSpitz
\rmove(0.8660254037844 -0.5)
\DreieckBreit
\rmove(0 1)
\DreieckSpitz
\rmove(0.8660254037844 -0.5)
\DreieckBreit
\rmove(0 1)
\DreieckSpitz
\rmove(0.8660254037844 -0.5)
\DreieckBreit
\rmove(0 1)
\DreieckSpitz
\rmove(0.8660254037844 -0.5)
\DreieckBreit
\rmove(0 1)
\DreieckSpitz
\rmove(0.8660254037844 -0.5)
\DreieckBreit
\rmove(0 1)
\DreieckSpitz
\rmove(0.8660254037844 -0.5)
\DreieckBreit
\rmove(0 1)
\DreieckSpitz
\rmove(0.8660254037844 -0.5)
\DreieckBreit
\rmove(-0.8660254037844 -0.5)
\rmove(-0.8660254037844 -0.5)
\rmove(-0.8660254037844 -0.5)
\rmove(-0.8660254037844 -0.5)
\rmove(-0.8660254037844 -0.5)
\rmove(-0.8660254037844 -0.5)
\DreieckSpitz
\rmove(0.8660254037844 -0.5)
\DreieckBreit
\rmove(0 1)
\DreieckSpitz
\rmove(0.8660254037844 -0.5)
\DreieckBreit
\rmove(0 1)
\DreieckSpitz
\rmove(0.8660254037844 -0.5)
\DreieckBreit
\rmove(0 1)
\DreieckSpitz
\rmove(0.8660254037844 -0.5)
\DreieckBreit
\rmove(0 1)
\DreieckSpitz
\rmove(0.8660254037844 -0.5)
\DreieckBreit
\rmove(0 1)
\DreieckSpitz
\rmove(0.8660254037844 -0.5)
\DreieckBreit
\rmove(0 1)
\DreieckSpitz
\rmove(0.8660254037844 -0.5)
\rmove(-0.8660254037844 -0.5)
\rmove(-0.8660254037844 -0.5)
\rmove(-0.8660254037844 -0.5)
\rmove(-0.8660254037844 -0.5)
\rmove(-0.8660254037844 -0.5)
\rmove(-0.8660254037844 -0.5)
\DreieckSpitz
\rmove(0.8660254037844 -0.5)
\DreieckBreit
\rmove(0 1)
\DreieckSpitz
\rmove(0.8660254037844 -0.5)
\DreieckBreit
\rmove(0 1)
\DreieckSpitz
\rmove(0.8660254037844 -0.5)
\DreieckBreit
\rmove(0 1)
\DreieckSpitz
\rmove(0.8660254037844 -0.5)
\DreieckBreit
\rmove(0 1)
\DreieckSpitz
\rmove(0.8660254037844 -0.5)
\DreieckBreit
\rmove(0 1)
\DreieckSpitz
\rmove(0.8660254037844 -0.5)
\rmove(-0.8660254037844 -0.5)
\rmove(-0.8660254037844 -0.5)
\rmove(-0.8660254037844 -0.5)
\rmove(-0.8660254037844 -0.5)
\rmove(-0.8660254037844 -0.5)
\DreieckSpitz
\rmove(0.8660254037844 -0.5)
\DreieckBreit
\rmove(0 1)
\DreieckSpitz
\rmove(0.8660254037844 -0.5)
\DreieckBreit
\rmove(0 1)
\DreieckSpitz
\rmove(0.8660254037844 -0.5)
\DreieckBreit
\rmove(0 1)
\DreieckSpitz
\rmove(0.8660254037844 -0.5)
\DreieckBreit
\rmove(0 1)
\DreieckSpitz
\rmove(0.8660254037844 -0.5)
\rmove(-0.8660254037844 -0.5)
\rmove(-0.8660254037844 -0.5)
\rmove(-0.8660254037844 -0.5)
\rmove(-0.8660254037844 -0.5)
\DreieckSpitz
\rmove(0.8660254037844 -0.5)
\DreieckBreit
\rmove(0 1)
\DreieckSpitz
\rmove(0.8660254037844 -0.5)
\DreieckBreit
\rmove(0 1)
\DreieckSpitz
\rmove(0.8660254037844 -0.5)
\DreieckBreit
\rmove(0 1)
\DreieckSpitz
\rmove(0.8660254037844 -0.5)
\move(-0.5 -1)
\htext{$c$}
\move(1.299038105677 1.75)
\rmove(0  0.5)
\htext{$a$}
\rmove(0 -0.5)
\rmove(1.299038105677 0.75)
\rmove(1.299038105677 -0.75)
\rmove(0.8660254037844 -0.5)
\rmove(0 0.5)
\htext{$b$}
\rmove(0 -0.5)
\rmove(1.299038105677 -0.75)
\rmove(0.8660254037844 -0.5)
\rmove(0 -2)
\rmove(0.5 0)
\htext{$c$}
\rmove(-0.5 0)
\rmove(0 -2)
\rmove(-1.299038105677 -0.75)
\rmove(0 -0.5)
\htext{$a$}
\rmove(0.5 0)
\rmove(-1.299038105677 -0.75)
\rmove(-1.299038105677 0.75)
\rmove(-0.8660254037844 0.5)
\rmove(-0.8660254037844 0.5)
\rmove(0 -0.5)
\htext{$b$}
\htext (-2.0 -6.5){\small  A hexagon with side lengths $a,b,c,a,b,c$, where $a=3$, $b=5$, $c=4$.}
}
\caption{}
\label{fig:1}
\end{figure}

\begin{figure}[hp]
\centertexdraw{
\drawdim cm
\linewd 0.08
\rmove(-5 0)
\RhombusPos
\rmove(0.8660254037844 -0.5)
\RhombusNeg
\rmove(0 1)
\RhombusVert
\rmove(0.8660254037844 0.5)
\RhombusVert

\rmove(-0.8660254037844 -0.5)
\rmove(-0.8660254037844 -0.5)
\rmove(0 -1)

\RhombusPos
\rmove(0.8660254037844 -0.5)
\RhombusNeg
\rmove(0.8660254037844 0.5)
\rmove(0 1)
\RhombusPos
\rmove(0.8660254037844 0.5)
\RhombusPos

\rmove(-0.8660254037844 -0.5)
\rmove(-0.8660254037844 -0.5)
\rmove(-0.8660254037844 -0.5)
\rmove(0 -1)

\RhombusPos
\rmove(0.8660254037844 -0.5)
\RhombusNeg
\rmove(0.8660254037844 0.5)
\RhombusNeg
\rmove(0 1)
\RhombusVert
\rmove(0.8660254037844 0.5)
\RhombusVert
\rmove(0.8660254037844 0.5)
\RhombusNeg

\rmove(-0.8660254037844 -0.5)
\rmove(-0.8660254037844 -0.5)
\rmove(-0.8660254037844 -0.5)
\rmove(-0.8660254037844 -0.5)
\rmove(0 -1)

\RhombusVert
\rmove(0.8660254037844 0.5)
\rmove(0.8660254037844 -0.5)
\RhombusNeg
\rmove(0.8660254037844 0.5)
\rmove(0 1)
\RhombusPos
\rmove(0.8660254037844 -0.5)
\RhombusNeg
\rmove(0 1)
\RhombusVert
\rmove(0.8660254037844 0.5)
\RhombusNeg

\rmove(-0.8660254037844 -0.5)
\rmove(-0.8660254037844 -0.5)
\rmove(-0.8660254037844 -0.5)
\rmove(-0.8660254037844 -0.5)
\rmove(-0.8660254037844 -0.5)
\rmove(0 -1)

\RhombusNeg
\rmove(0.8660254037844 0.5)
\RhombusPos
\rmove(0.8660254037844 -0.5)
\RhombusNeg
\rmove(0.8660254037844 0.5)
\RhombusNeg
\rmove(0 1)
%Mitte
\RhombusVert
\rmove(0.8660254037844 0.5)
\rmove(0.8660254037844 0.5)
\RhombusPos
\rmove(0.8660254037844 0.5)
\RhombusNeg
\rmove(0.8660254037844 -0.5)
\RhombusNeg

\rmove(-0.8660254037844 -0.5)
\rmove(-0.8660254037844 -0.5)
\rmove(-0.8660254037844 -0.5)
\rmove(-0.8660254037844 -0.5)
\rmove(-0.8660254037844 -0.5)
\rmove(-0.8660254037844 -0.5)

\RhombusVert
\rmove(0.8660254037844 0.5)
\rmove(0.8660254037844 -0.5)
\RhombusNeg
\rmove(0.8660254037844 0.5)
\RhombusNeg
\rmove(0 1)
\RhombusVert
\rmove(0.8660254037844 0.5)
\RhombusVert
\rmove(0.8660254037844 0.5)
\RhombusNeg
\rmove(0.8660254037844 -0.5)
\RhombusNeg

\rmove(-0.8660254037844 -0.5)
\rmove(-0.8660254037844 -0.5)
\rmove(-0.8660254037844 -0.5)
\rmove(-0.8660254037844 -0.5)
\rmove(-0.8660254037844 -0.5)

\RhombusVert
\rmove(0.8660254037844 0.5)
\rmove(0.8660254037844 -0.5)
\RhombusNeg
\rmove(0.8660254037844 0.5)
\rmove(0 1)
\RhombusPos
\rmove(0.8660254037844 0.5)
\RhombusPos
\rmove(0.8660254037844 -0.5)
\RhombusNeg

\rmove(-0.8660254037844 -0.5)
\rmove(-0.8660254037844 -0.5)
\rmove(-0.8660254037844 -0.5)
\rmove(-0.8660254037844 -0.5)

\RhombusVert 
\rmove(0.8660254037844 0.5)
\rmove(0.8660254037844 0.5)
\RhombusPos
\rmove(0.8660254037844 0.5)
\RhombusVert

\rmove(-0.8660254037844 -0.5)
\rmove(-0.8660254037844 -0.5)
\rmove(0 -1)

\RhombusVert
\rmove(0.8660254037844 0.5)
\RhombusVert
\rmove(0.8660254037844 0.5)
\RhombusNeg
\rmove(0.8660254037844 0.5)
\RhombusPos

\move(0 -1)
\rmove(-5 0)
\linewd 0.01
\DreieckBreit
\rmove(0 1)
\DreieckSpitz
\rmove(0.8660254037844 -0.5)
\DreieckBreit
\rmove(0 1)
\DreieckSpitz
\rmove(0.8660254037844 -0.5)
\DreieckBreit
\rmove(0 1)
\DreieckSpitz
\rmove(0.8660254037844 -0.5)
\DreieckBreit

\rmove(-0.8660254037844 -0.5)
\rmove(-0.8660254037844 -0.5)
\rmove(-0.8660254037844 -0.5)
\rmove(0 -1)

\DreieckBreit
\rmove(0 1)
\DreieckSpitz
\rmove(0.8660254037844 -0.5)
\DreieckBreit
\rmove(0 1)
\DreieckSpitz
\rmove(0.8660254037844 -0.5)
\DreieckBreit
\rmove(0 1)
\DreieckSpitz
\rmove(0.8660254037844 -0.5)
\DreieckBreit
\rmove(0 1)
\DreieckSpitz
\rmove(0.8660254037844 -0.5)
\DreieckBreit

\rmove(-0.8660254037844 -0.5)
\rmove(-0.8660254037844 -0.5)
\rmove(-0.8660254037844 -0.5)
\rmove(-0.8660254037844 -0.5)
\rmove(0 -1)

\DreieckBreit
\rmove(0 1)
\DreieckSpitz
\rmove(0.8660254037844 -0.5)
\DreieckBreit
\rmove(0 1)
\DreieckSpitz
\rmove(0.8660254037844 -0.5)
\DreieckBreit
\rmove(0 1)
\DreieckSpitz
\rmove(0.8660254037844 -0.5)
\DreieckBreit
\rmove(0 1)
\DreieckSpitz
\rmove(0.8660254037844 -0.5)
\DreieckBreit
\rmove(0 1)
\DreieckSpitz
\rmove(0.8660254037844 -0.5)
\DreieckBreit

\rmove(-0.8660254037844 -0.5)
\rmove(-0.8660254037844 -0.5)
\rmove(-0.8660254037844 -0.5)
\rmove(-0.8660254037844 -0.5)
\rmove(-0.8660254037844 -0.5)
\rmove(0 -1)

\DreieckBreit
\rmove(0 1)
\DreieckSpitz
\rmove(0.8660254037844 -0.5)
\DreieckBreit
\rmove(0 1)
\DreieckSpitz
\rmove(0.8660254037844 -0.5)
\DreieckBreit
\rmove(0 1)
\DreieckSpitz
\rmove(0.8660254037844 -0.5)
\DreieckBreit
\rmove(0 1)
\DreieckSpitz
\rmove(0.8660254037844 -0.5)
\DreieckBreit
\rmove(0 1)
\DreieckSpitz
\rmove(0.8660254037844 -0.5)
\DreieckBreit
\rmove(0 1)
\DreieckSpitz
\rmove(0.8660254037844 -0.5)
\DreieckBreit

\rmove(-0.8660254037844 -0.5)
\rmove(-0.8660254037844 -0.5)
\rmove(-0.8660254037844 -0.5)
\rmove(-0.8660254037844 -0.5)
\rmove(-0.8660254037844 -0.5)
\rmove(-0.8660254037844 -0.5)

\DreieckSpitz
\rmove(0.8660254037844 -0.5)
\DreieckBreit
\rmove(0 1)
\DreieckSpitz
\rmove(0.8660254037844 -0.5)
\DreieckBreit
\rmove(0 1)
\DreieckSpitz
\rmove(0.8660254037844 -0.5)
\DreieckBreit
\rmove(0 1)
\DreieckSpitz
\rmove(0.8660254037844 -0.5)
\DreieckBreit
\rmove(0 1)
\DreieckSpitz
\rmove(0.8660254037844 -0.5)
\DreieckBreit
\rmove(0 1)
\DreieckSpitz
\rmove(0.8660254037844 -0.5)
\DreieckBreit
\rmove(0 1)
\DreieckSpitz
\rmove(0.8660254037844 -0.5)
\DreieckBreit

\rmove(-0.8660254037844 -0.5)
\rmove(-0.8660254037844 -0.5)
\rmove(-0.8660254037844 -0.5)
\rmove(-0.8660254037844 -0.5)
\rmove(-0.8660254037844 -0.5)
\rmove(-0.8660254037844 -0.5)

\DreieckSpitz
\rmove(0.8660254037844 -0.5)
\DreieckBreit
\rmove(0 1)
\DreieckSpitz
\rmove(0.8660254037844 -0.5)
\DreieckBreit
\rmove(0 1)
\DreieckSpitz
\rmove(0.8660254037844 -0.5)
\DreieckBreit
\rmove(0 1)
\DreieckSpitz
\rmove(0.8660254037844 -0.5)
\DreieckBreit
\rmove(0 1)
\DreieckSpitz
\rmove(0.8660254037844 -0.5)
\DreieckBreit
\rmove(0 1)
\DreieckSpitz
\rmove(0.8660254037844 -0.5)
\DreieckBreit
\rmove(0 1)
\DreieckSpitz
\rmove(0.8660254037844 -0.5)

\rmove(-0.8660254037844 -0.5)
\rmove(-0.8660254037844 -0.5)
\rmove(-0.8660254037844 -0.5)
\rmove(-0.8660254037844 -0.5)
\rmove(-0.8660254037844 -0.5)
\rmove(-0.8660254037844 -0.5)

\DreieckSpitz
\rmove(0.8660254037844 -0.5)
\DreieckBreit
\rmove(0 1)
\DreieckSpitz
\rmove(0.8660254037844 -0.5)
\DreieckBreit
\rmove(0 1)
\DreieckSpitz
\rmove(0.8660254037844 -0.5)
\DreieckBreit
\rmove(0 1)
\DreieckSpitz
\rmove(0.8660254037844 -0.5)
\DreieckBreit
\rmove(0 1)
\DreieckSpitz
\rmove(0.8660254037844 -0.5)
\DreieckBreit
\rmove(0 1)
\DreieckSpitz
\rmove(0.8660254037844 -0.5)

\rmove(-0.8660254037844 -0.5)
\rmove(-0.8660254037844 -0.5)
\rmove(-0.8660254037844 -0.5)
\rmove(-0.8660254037844 -0.5)
\rmove(-0.8660254037844 -0.5)

\DreieckSpitz
\rmove(0.8660254037844 -0.5)
\DreieckBreit
\rmove(0 1)
\DreieckSpitz
\rmove(0.8660254037844 -0.5)
\DreieckBreit
\rmove(0 1)
\DreieckSpitz
\rmove(0.8660254037844 -0.5)
\DreieckBreit
\rmove(0 1)
\DreieckSpitz
\rmove(0.8660254037844 -0.5)
\DreieckBreit
\rmove(0 1)
\DreieckSpitz
\rmove(0.8660254037844 -0.5)

\rmove(-0.8660254037844 -0.5)
\rmove(-0.8660254037844 -0.5)
\rmove(-0.8660254037844 -0.5)
\rmove(-0.8660254037844 -0.5)

\DreieckSpitz
\rmove(0.8660254037844 -0.5)
\DreieckBreit
\rmove(0 1)
\DreieckSpitz
\rmove(0.8660254037844 -0.5)
\DreieckBreit
\rmove(0 1)
\DreieckSpitz
\rmove(0.8660254037844 -0.5)
\DreieckBreit
\rmove(0 1)
\DreieckSpitz
\rmove(0.8660254037844 -0.5)
\htext (-6  -7.5){\small A rhombus tiling of a hexagon with side lengths $a,b,c,a,b,c$.}
}
\caption{}
\label{fig:2}
\end{figure}

The main results of this paper are the following two theorems.
\begin{theo} 
\label{th:hp}
Let $a$,$b$,$c$ be positive integers with $a \equiv b \hspace{2mm} (\text{mod} 
\hspace{1mm} 2)$ and $a \not\equiv c \hspace{2mm} (\text{mod} \hspace{1mm} 2)$. 
Then the number of rhombus tilings of a hexagon with side lengths
$a$,$b$,$c$,$a$,$b$,$c$ which contain the rhombus in the 
centre is
\begin{multline}
\label{number:odd}
\left( \prod_{i=1}^a \prod_{j=1}^b \prod_{k=1}^c \frac{i+j+k-1}{i+j+k-2}\right)
   \frac{(1)_c}{(b+1)_{c+a-1}}  
\binom{\frac{b+c-1}{2}}{\frac{b-1}{2}} \binom{\frac{a+b+c-2}{2}}{\frac{b-1}{2}} 2^{a-1}\\ 
\times \sum_{k=0}^{(a-1)/2} \left[ \left(\frac{c+1}{2}\right)_k 
      \left(\frac{1+b+c}{2}\right)_{k}  
      \left(\frac{c+2k+2}{2}\right)_{(a-2k-1)/2} \right. \\ 
       \hspace{2cm} \left. \times  
      \left(\frac{b+c+2k+3}{2}\right)_{(a-2k-1)/2}
      \frac{(\frac{1}{2})_{(a-2k-1)/2}}{(1)_{(a-2k-1)/2}} \right]  
\end{multline}
in case that $a$ is odd, and 
\begin{multline}
\label{number:even}
\left( \prod_{i=1}^a \prod_{j=1}^b \prod_{k=1}^c \frac{i+j+k-1}{i+j+k-2}\right)
   \frac{b (1)_c}{(b+1)_{c+a-1}}       
\binom{\frac{b+c-1}{2}}{\frac{b}{2}} \binom{\frac{a+b+c-1}{2}}{\frac{b}{2}} 2^{a-2} \\
\times \sum_{k=0}^{(a-2)/2} \left[ \left(\frac{c+2}{2}\right)_k           
      \left(\frac{1+b+c}{2}\right)_{k} 
      \left(\frac{c+2k+3}{2}\right)_{(a-2k-2)/2} \right. \\  
      \hspace{2cm} \left. \times 
      \left(\frac{b+c+2k+3}{2}\right)_{(a-2k-2)/2}
      \frac{(\frac{1}{2})_{(a-2k-2)/2}}{(1)_{(a-2k-2)/2}} \right] 
\end{multline}
in case that $a$ is even, where the Pochhammer symbol $(a)_{k}$ is defined by 
$(a)_{k}:=a(a+1) \dots (a+k-1)$ if  $k \ge 1$ and $(a)_0:=1$.
\end{theo}

\begin{figure}[h]
\centertexdraw{
\drawdim cm
\linewd 0.02

\rmove(0.8660254037844 -0.5)
\rmove(0.8660254037844 -0.5)
\rmove(0.8660254037844 -0.5)

\rlvec(0.8660254037844 0.5)
\rlvec(0.8660254037844 -0.5)
\rlvec(-0.8660254037844 -0.5)

\lfill f:0.5

\move(0 0)
\DreieckBreit
\rmove(0 1)
\DreieckSpitz
\rmove(0.8660254037844 -0.5)
\DreieckBreit
\rmove(0 1)
\DreieckSpitz
\rmove(0.8660254037844 -0.5)
\DreieckBreit
\rmove(0 1)
\DreieckSpitz
\rmove(0.8660254037844 -0.5)
\DreieckBreit

\rmove(-0.8660254037844 -0.5)
\rmove(-0.8660254037844 -0.5)
\rmove(-0.8660254037844 -0.5)
\rmove(0 -1)

\DreieckBreit
\rmove(0 1)
\DreieckSpitz
\rmove(0.8660254037844 -0.5)
\DreieckBreit
\rmove(0 1)
\DreieckSpitz
\rmove(0.8660254037844 -0.5)
\DreieckBreit
\rmove(0 1)
\DreieckSpitz
\rmove(0.8660254037844 -0.5)
\DreieckBreit
\rmove(0 1)
\DreieckSpitz
\rmove(0.8660254037844 -0.5)
\DreieckBreit

\rmove(-0.8660254037844 -0.5)
\rmove(-0.8660254037844 -0.5)
\rmove(-0.8660254037844 -0.5)
\rmove(-0.8660254037844 -0.5)
\rmove(0 -1)

\DreieckBreit
\rmove(0 1)
\DreieckSpitz
\rmove(0.8660254037844 -0.5)
\DreieckBreit
\rmove(0 1)
\DreieckSpitz
\rmove(0.8660254037844 -0.5)
\DreieckBreit
\rmove(0 1)
\DreieckSpitz
\rmove(0.8660254037844 -0.5)
\DreieckBreit
\rmove(0 1)
\DreieckSpitz
\rmove(0.8660254037844 -0.5)
\DreieckBreit
\rmove(0 1)
\DreieckSpitz
\rmove(0.8660254037844 -0.5)
\DreieckBreit

\rmove(-0.8660254037844 -0.5)
\rmove(-0.8660254037844 -0.5)
\rmove(-0.8660254037844 -0.5)
\rmove(-0.8660254037844 -0.5)
\rmove(-0.8660254037844 -0.5)
\rmove(0 -1)

\DreieckBreit
\rmove(0 1)
\DreieckSpitz
\rmove(0.8660254037844 -0.5)
\DreieckBreit
\rmove(0 1)
\DreieckSpitz
\rmove(0.8660254037844 -0.5)
\DreieckBreit
\rmove(0 1)
\DreieckSpitz
\rmove(0.8660254037844 -0.5)
\DreieckBreit
\rmove(0 1)
\DreieckSpitz
\rmove(0.8660254037844 -0.5)
\DreieckBreit
\rmove(0 1)
\DreieckSpitz
\rmove(0.8660254037844 -0.5)
\DreieckBreit
\rmove(0 1)
\DreieckSpitz
\rmove(0.8660254037844 -0.5)
\DreieckBreit

\rmove(-0.8660254037844 -0.5)
\rmove(-0.8660254037844 -0.5)
\rmove(-0.8660254037844 -0.5)
\rmove(-0.8660254037844 -0.5)
\rmove(-0.8660254037844 -0.5)
\rmove(-0.8660254037844 -0.5)

\DreieckSpitz
\rmove(0.8660254037844 -0.5)
\DreieckBreit
\rmove(0 1)
\DreieckSpitz
\rmove(0.8660254037844 -0.5)
\DreieckBreit
\rmove(0 1)
\DreieckSpitz
\rmove(0.8660254037844 -0.5)
\DreieckBreit
\rmove(0 1)
\DreieckSpitz
\rmove(0.8660254037844 -0.5)
\DreieckBreit
\rmove(0 1)
\DreieckSpitz
\rmove(0.8660254037844 -0.5)
\DreieckBreit
\rmove(0 1)
\DreieckSpitz
\rmove(0.8660254037844 -0.5)
\DreieckBreit
\rmove(0 1)
\DreieckSpitz
\rmove(0.8660254037844 -0.5)
\DreieckBreit

\rmove(-0.8660254037844 -0.5)
\rmove(-0.8660254037844 -0.5)
\rmove(-0.8660254037844 -0.5)
\rmove(-0.8660254037844 -0.5)
\rmove(-0.8660254037844 -0.5)
\rmove(-0.8660254037844 -0.5)

\DreieckSpitz
\rmove(0.8660254037844 -0.5)
\DreieckBreit
\rmove(0 1)
\DreieckSpitz
\rmove(0.8660254037844 -0.5)
\DreieckBreit
\rmove(0 1)
\DreieckSpitz
\rmove(0.8660254037844 -0.5)
\DreieckBreit
\rmove(0 1)
\DreieckSpitz
\rmove(0.8660254037844 -0.5)
\DreieckBreit
\rmove(0 1)
\DreieckSpitz
\rmove(0.8660254037844 -0.5)
\DreieckBreit
\rmove(0 1)
\DreieckSpitz
\rmove(0.8660254037844 -0.5)
\DreieckBreit
\rmove(0 1)
\DreieckSpitz
\rmove(0.8660254037844 -0.5)

\rmove(-0.8660254037844 -0.5)
\rmove(-0.8660254037844 -0.5)
\rmove(-0.8660254037844 -0.5)
\rmove(-0.8660254037844 -0.5)
\rmove(-0.8660254037844 -0.5)
\rmove(-0.8660254037844 -0.5)

\DreieckSpitz
\rmove(0.8660254037844 -0.5)
\DreieckBreit
\rmove(0 1)
\DreieckSpitz
\rmove(0.8660254037844 -0.5)
\DreieckBreit
\rmove(0 1)
\DreieckSpitz
\rmove(0.8660254037844 -0.5)
\DreieckBreit
\rmove(0 1)
\DreieckSpitz
\rmove(0.8660254037844 -0.5)
\DreieckBreit
\rmove(0 1)
\DreieckSpitz
\rmove(0.8660254037844 -0.5)
\DreieckBreit
\rmove(0 1)
\DreieckSpitz
\rmove(0.8660254037844 -0.5)

\rmove(-0.8660254037844 -0.5)
\rmove(-0.8660254037844 -0.5)
\rmove(-0.8660254037844 -0.5)
\rmove(-0.8660254037844 -0.5)
\rmove(-0.8660254037844 -0.5)

\DreieckSpitz
\rmove(0.8660254037844 -0.5)
\DreieckBreit
\rmove(0 1)
\DreieckSpitz
\rmove(0.8660254037844 -0.5)
\DreieckBreit
\rmove(0 1)
\DreieckSpitz
\rmove(0.8660254037844 -0.5)
\DreieckBreit
\rmove(0 1)
\DreieckSpitz
\rmove(0.8660254037844 -0.5)
\DreieckBreit
\rmove(0 1)
\DreieckSpitz
\rmove(0.8660254037844 -0.5)

\rmove(-0.8660254037844 -0.5)
\rmove(-0.8660254037844 -0.5)
\rmove(-0.8660254037844 -0.5)
\rmove(-0.8660254037844 -0.5)

\DreieckSpitz
\rmove(0.8660254037844 -0.5)
\DreieckBreit
\rmove(0 1)
\DreieckSpitz
\rmove(0.8660254037844 -0.5)
\DreieckBreit
\rmove(0 1)
\DreieckSpitz
\rmove(0.8660254037844 -0.5)
\DreieckBreit
\rmove(0 1)
\DreieckSpitz
\rmove(0.8660254037844 -0.5)
\htext (-3 -7.0){\small A hexagon with side lengths $a,b,c,a,b,c$, where the central 
rhombus is marked.}
}
\caption{}
\label{fig:3}
\end{figure}

The special case $a=b$ in Theorem~\ref{th:hp} was previously derived in 
\cite[Theorem 1, 2]{CiKrat}. (In order to see that the two sums in 
Theorem~\ref{th:hp} are equal to two sums in Theorem~1 and Theorem~2 from  
\cite{CiKrat} in this special case, one 
has to apply Bailey's transformation \eqref{bai} for balanced $_4 F _3$--series.)

The assumption about the parity of the side lengths $a$,$b$,$c$, i.e., that 
not all side lengths of the hexagon have the same parity, comes from the fact 
that this condition is necessary for the existence of a central rhombus in a rhombus 
tiling of a hexagon with side lengths $a$,$b$,$c$,$a$,$b$,$c$.  There is also a result
similar to Theorem~\ref{th:hp} for a hexagon whose side 
lengths have the same parity. In this case we choose for the fixed rhombus 
the horizontal rhombus whose lowest vertex is the centre of 
the hexagon. (We could have also chosen the horizontal rhombus whose uppermost vertex is 
the centre of the hexagon).  The number of these rhombus tilings is 
given by formulas quite similar to those in Theorem~\ref{th:hp},
and, furthermore, the proofs of these formulas are analogous to the 
proofs of the formulas in Theorem~\ref{th:hp}. We obtain the following:

\begin{theo} 
\label{th:hp2}
Let $a$,$b$,$c$ be positive integers with $a \equiv b \equiv c \hspace{2mm} (\text{mod} 
\hspace{1mm} 2)$. 
Then the number of rhombus tilings of hexagon with side lengths
$a$,$b$,$c$,$a$,$b$,$c$ which contain the horizontal rhombus whose lowest vertex is the centre  
of the hexagon is equal to
\begin{multline}
\label{number:odd2}
\left( \prod_{i=1}^a \prod_{j=1}^b \prod_{k=1}^c \frac{i+j+k-1}{i+j+k-2}\right)
   \frac{(1)_c}{(b+1)_{c+a-1}}  
\binom{\frac{b+c-2}{2}}{\frac{b-1}{2}} \binom{\frac{a+b+c-1}{2}}{\frac{b-1}{2}} 2^{a-1}\\ 
\times \left( \left(\frac{c+1}{2} \right)_{(a-1)/2} \left( \frac{b+c+2}{2} \right)_{(a-1)/2} 
       \frac{\left( \frac{1}{2} \right)_{(a-1)/2}}
            {(1)_{(a-1)/2}}  \right. \\
+ \sum_{k=1}^{(a-1)/2}  \left(\frac{c+2}{2}\right)_{k-1} 
      \left(\frac{b+c}{2}\right)_{k}  
      \left(\frac{c+2k+1}{2}\right)_{(a-2k+1)/2} \\ 
       \left. \times  
      \left(\frac{b+c+2k+2}{2}\right)_{(a-2k-1)/2}
      \frac{(\frac{1}{2})_{(a-2k-1)/2}}{(1)_{(a-2k-1)/2}}  \right)
\end{multline}
in case that $a$ is odd, and 
\begin{multline}
\label{number:even2}
\left( \prod_{i=1}^a \prod_{j=1}^b \prod_{k=1}^c \frac{i+j+k-1}{i+j+k-2}\right)
   \frac{ (1)_c}{(b)_{c+a}}       
\binom{\frac{b+c-2}{2}}{\frac{b-2}{2}} \binom{\frac{a+b+c-2}{2}}{\frac{b-2}{2}} 2^{a} \\
\times \left( \left(\frac{c+2}{2} \right)_{(a-2)/2} \left( \frac{b+c+2}{2} \right)_{a/2} 
       \frac{\left( \frac{1}{2} \right)_{a/2}}
            {(1)_{(a-2)/2}}  \right. \\
+ \sum_{k=1}^{a/2}  \left(\frac{c+1}{2}\right)_k           
      \left(\frac{b+c}{2}\right)_{k} 
      \left(\frac{c+2k+2}{2}\right)_{(a-2k)/2}  \\  
       \left. \times 
      \left(\frac{b+c+2k+2}{2}\right)_{(a-2k)/2}
      \frac{(\frac{1}{2})_{(a-2k)/2}}{(1)_{(a-2k)/2}} \right)
\end{multline}
in case that $a$ is even.
\end{theo}

\begin{figure}[h]
\centertexdraw{
\drawdim cm
\linewd 0.02

\rmove(0.8660254037844 -0.5)
\rmove(0.8660254037844 -0.5)
\rmove(0.8660254037844 -0.5)
\rmove(0 1)

\rlvec(0.8660254037844 0.5)
\rlvec(0.8660254037844 -0.5)
\rlvec(-0.8660254037844 -0.5)

\lfill f:0.5

\move(0 0)
\DreieckBreit
\rmove(0 1)
\DreieckSpitz
\rmove(0.8660254037844 -0.5)
\DreieckBreit
\rmove(0 1)
\DreieckSpitz
\rmove(0.8660254037844 -0.5)
\DreieckBreit
\rmove(0 1)
\DreieckSpitz
\rmove(0.8660254037844 -0.5)
\DreieckBreit

\rmove(-0.8660254037844 -0.5)
\rmove(-0.8660254037844 -0.5)
\rmove(-0.8660254037844 -0.5)
\rmove(0 -1)

\DreieckBreit
\rmove(0 1)
\DreieckSpitz
\rmove(0.8660254037844 -0.5)
\DreieckBreit
\rmove(0 1)
\DreieckSpitz
\rmove(0.8660254037844 -0.5)
\DreieckBreit
\rmove(0 1)
\DreieckSpitz
\rmove(0.8660254037844 -0.5)
\DreieckBreit
\rmove(0 1)
\DreieckSpitz
\rmove(0.8660254037844 -0.5)
\DreieckBreit

\rmove(-0.8660254037844 -0.5)
\rmove(-0.8660254037844 -0.5)
\rmove(-0.8660254037844 -0.5)
\rmove(-0.8660254037844 -0.5)
\rmove(0 -1)

\DreieckBreit
\rmove(0 1)
\DreieckSpitz
\rmove(0.8660254037844 -0.5)
\DreieckBreit
\rmove(0 1)
\DreieckSpitz
\rmove(0.8660254037844 -0.5)
\DreieckBreit
\rmove(0 1)
\DreieckSpitz
\rmove(0.8660254037844 -0.5)
\DreieckBreit
\rmove(0 1)
\DreieckSpitz
\rmove(0.8660254037844 -0.5)
\DreieckBreit
\rmove(0 1)
\DreieckSpitz
\rmove(0.8660254037844 -0.5)
\DreieckBreit

\rmove(-0.8660254037844 -0.5)
\rmove(-0.8660254037844 -0.5)
\rmove(-0.8660254037844 -0.5)
\rmove(-0.8660254037844 -0.5)
\rmove(-0.8660254037844 -0.5)
\rmove(0 -1)

%\DreieckBreit
\rmove(0 1)
\DreieckSpitz
\rmove(0.8660254037844 -0.5)
\DreieckBreit
\rmove(0 1)
\DreieckSpitz
\rmove(0.8660254037844 -0.5)
\DreieckBreit
\rmove(0 1)
\DreieckSpitz
\rmove(0.8660254037844 -0.5)
\DreieckBreit
\rmove(0 1)
\DreieckSpitz
\rmove(0.8660254037844 -0.5)
\DreieckBreit
\rmove(0 1)
\DreieckSpitz
\rmove(0.8660254037844 -0.5)
\DreieckBreit
\rmove(0 1)
\DreieckSpitz
\rmove(0.8660254037844 -0.5)
\DreieckBreit

\rmove(-0.8660254037844 -0.5)
\rmove(-0.8660254037844 -0.5)
\rmove(-0.8660254037844 -0.5)
\rmove(-0.8660254037844 -0.5)
\rmove(-0.8660254037844 -0.5)
%\rmove(-0.8660254037844 -0.5)

%\DreieckSpitz
%\rmove(0.8660254037844 -0.5)
%\DreieckBreit
%\rmove(0 1)
\DreieckSpitz
\rmove(0.8660254037844 -0.5)
\DreieckBreit
\rmove(0 1)
\DreieckSpitz
\rmove(0.8660254037844 -0.5)
\DreieckBreit
\rmove(0 1)
\DreieckSpitz
\rmove(0.8660254037844 -0.5)
\DreieckBreit
\rmove(0 1)
\DreieckSpitz
\rmove(0.8660254037844 -0.5)
\DreieckBreit
\rmove(0 1)
\DreieckSpitz
\rmove(0.8660254037844 -0.5)
\DreieckBreit
\rmove(0 1)
\DreieckSpitz
\rmove(0.8660254037844 -0.5)
\DreieckBreit

\rmove(-0.8660254037844 -0.5)
\rmove(-0.8660254037844 -0.5)
\rmove(-0.8660254037844 -0.5)
\rmove(-0.8660254037844 -0.5)
\rmove(-0.8660254037844 -0.5)
%\rmove(-0.8660254037844 -0.5)

%\DreieckSpitz
%\rmove(0.8660254037844 -0.5)
%\DreieckBreit
%\rmove(0 1)
\DreieckSpitz
\rmove(0.8660254037844 -0.5)
\DreieckBreit
\rmove(0 1)
\DreieckSpitz
\rmove(0.8660254037844 -0.5)
\DreieckBreit
\rmove(0 1)
\DreieckSpitz
\rmove(0.8660254037844 -0.5)
\DreieckBreit
\rmove(0 1)
\DreieckSpitz
\rmove(0.8660254037844 -0.5)
\DreieckBreit
\rmove(0 1)
\DreieckSpitz
\rmove(0.8660254037844 -0.5)
\DreieckBreit
\rmove(0 1)
\DreieckSpitz
\rmove(0.8660254037844 -0.5)

\rmove(-0.8660254037844 -0.5)
\rmove(-0.8660254037844 -0.5)
\rmove(-0.8660254037844 -0.5)
\rmove(-0.8660254037844 -0.5)
\rmove(-0.8660254037844 -0.5)
%\rmove(-0.8660254037844 -0.5)

%\DreieckSpitz
%\rmove(0.8660254037844 -0.5)
%\DreieckBreit
%\rmove(0 1)
\DreieckSpitz
\rmove(0.8660254037844 -0.5)
\DreieckBreit
\rmove(0 1)
\DreieckSpitz
\rmove(0.8660254037844 -0.5)
\DreieckBreit
\rmove(0 1)
\DreieckSpitz
\rmove(0.8660254037844 -0.5)
\DreieckBreit
\rmove(0 1)
\DreieckSpitz
\rmove(0.8660254037844 -0.5)
\DreieckBreit
\rmove(0 1)
\DreieckSpitz
\rmove(0.8660254037844 -0.5)

\rmove(-0.8660254037844 -0.5)
\rmove(-0.8660254037844 -0.5)
\rmove(-0.8660254037844 -0.5)
\rmove(-0.8660254037844 -0.5)
%\rmove(-0.8660254037844 -0.5)

%\DreieckSpitz
%\rmove(0.8660254037844 -0.5)
%\DreieckBreit
%\rmove(0 1)
\DreieckSpitz
\rmove(0.8660254037844 -0.5)
\DreieckBreit
\rmove(0 1)
\DreieckSpitz
\rmove(0.8660254037844 -0.5)
\DreieckBreit
\rmove(0 1)
\DreieckSpitz
\rmove(0.8660254037844 -0.5)
\DreieckBreit
\rmove(0 1)
\DreieckSpitz
\rmove(0.8660254037844 -0.5)

%\rmove(-0.8660254037844 -0.5)
%\rmove(-0.8660254037844 -0.5)
%\rmove(-0.8660254037844 -0.5)
%\rmove(-0.8660254037844 -0.5)

%\DreieckSpitz
%\rmove(0.8660254037844 -0.5)
%\DreieckBreit
%\rmove(0 1)
%\DreieckSpitz
%\rmove(0.8660254037844 -0.5)
%\DreieckBreit
%\rmove(0 1)
%\DreieckSpitz
%\rmove(0.8660254037844 -0.5)
%\DreieckBreit
%\rmove(0 1)
%\DreieckSpitz
%\rmove(0.8660254037844 -0.5)
\htext (-3 -5.5){\small A hexagon with side lengths $a,b,c,a,b,c$, }
\htext (-3 -6){\small where the `almost central' 
rhombus above the center is marked.}
}
\caption{}
\label{fig:20}
\end{figure}

The special case $a=b$ in Theorem~\ref{th:hp2} was previously derived 
in \cite[Theorem 1, 2]{FulKrat1}. (In order to see that the two sums 
in Theorem~\ref{th:hp2} are equal to the two sums  
in Theorem~1 and Theorem~2 from \cite{FulKrat1} for this special case, 
one first has to apply a certain contiguous relation on the balanced $_4 F _3$--series in 
Theorem~\ref{th:hp2}. In each case we obtain the sum of a balanced 
$_3 F _2$--series, which can be evaluated by Saalsch\"utz's summation 
formula (see \cite[(2.3.1.3)]{slater}), and another balanced $_4 F _3$--series.
An application of Bailey's transformation to the latter  
$_4 F _3$--series  finally shows the equivalence of the 
result in Theorem~\ref{th:hp2} and the result in  \cite[Theorem 1, 2]{FulKrat1} for 
the special case.)

\bigskip
   
These  enumeration results for rhombus tilings with a fixed rhombus are not only interesting 
because they add more results to the growing set of results on the enumeration 
of rhombus tilings with special properties. These enumeration also contribute to  
the interesting question of what a `typical' rhombus tiling of a hexagon 
looks like. Cohen, Larsen und Propp address this question in 
\cite{propp}. We are able to deduce the following 
two theorems from Theorem~\ref{th:hp} and Theorem~\ref{th:hp2}.
They confirm Conjecture~1 in \cite{propp} for a special case.

\begin{theo}
\label{asym:th}
Let $\alpha$, $\beta$ and $\gamma$ be non-negative real numbers. Then the 
probability that a rhombus tiling of the hexagon with side lengths $a,b,c,a,b,c$, where
$a \sim \alpha N$, $b \sim \beta N$, $c \sim \gamma N$ and  $a \equiv b \not\equiv c \hspace{2mm} (\text{mod} 
\hspace{1mm} 2)$, contains the rhombus in the centre 
is asymptotically 
\begin{equation}
\label{asym}
\frac{2}{\pi} 
\arcsin\left( \sqrt{\frac{\alpha \beta}{(\beta + \gamma)(\alpha + \gamma)}} \right)
\end{equation}
as $N$ tends to infinity.
\end{theo} 

\begin{theo}
\label{asym:th2}
Let $\alpha$, $\beta$ and $\gamma$ be non-negative real numbers. Then the 
probability that a rhombus tiling of the hexagon with side lengths $a,b,c,a,b,c$, where
$a \sim \alpha N$, $b \sim \beta N$, $c \sim \gamma N$ and $a \equiv b \equiv c \hspace{2mm} (\text{mod} 
\hspace{1mm} 2)$, contains the `almost central' rhombus 
above the centre  
is asymptotically 
\begin{equation}
\label{asym2}
\frac{2}{\pi} 
\arcsin\left( \sqrt{\frac{\alpha \beta}{(\beta + \gamma)(\alpha + \gamma)}} \right)
\end{equation}
as $N$ tends to infinity.
\end{theo}

Roughly speaking, 
Conjecture 1 in \cite{propp} predicts the following: Fix an arbitrary point $(x,y)$ in 
the hexagon with side lengths $a$,$b$,$c$,$a$,$b$,$c$. 
Then the probability that a random rhombus tiling 
of a hexagon contains a vertical rhombus at this point tends to
$\mathcal{P}_{\alpha,\beta,\gamma}(x,y)$ as the hexagon becomes large. Here, 
$\mathcal{P}_{\alpha,\beta,\gamma}$ is a certain function (defined 
in \cite[Theorem 1]{propp}) that depends on the proportions $\alpha$,$\beta$,$\gamma$ 
of the side lengths of the hexagon (see Theorem~\ref{asym:th}). 
Theorem~\ref{asym:th} and  Theorem~\ref{asym:th2} confirm this conjecture for the centre of the hexagon.
(The reader should know that Cohen, Larsen and Propp use another coordinate system in 
\cite{propp} than we do and therefore they are in the position to speak of vertical rhombi.
But of course, their result can easily be translated into our coordinate system where we fix
horizontal rhombi instead of vertical rhombi.) 

We want to direct the reader's attention to a side result of the present work, given 
in Lemma~\ref{dreisum}. It expresses, for {\em any} fixed rhombus, the number of 
rhombus tilings containing that fixed rhombus as a triple sum. I am currently 
pursuing an asymptotic analysis of this triple sum with the ultimate goal of 
effectively proving Cohn, Larsen and Propp's conjecture.

\bigskip

In order to prove Theorem~\ref{th:hp} and Theorem~\ref{th:hp2}  we make use of a 
bijection between rhombus tilings which contain the central rhombus, 
respectively rhombus tilings which contain the `almost central' rhombus above the 
centre,  and 
non-intersecting lattice paths. This bijection actually works for arbitrary rhombi, i.e., the 
fixed rhombus need not be placed in the centre or next to the centre. The 
description of this bijection  is the subject of Section~\ref{bijection}. 
Since the number of non-intersecting 
lattice paths is given by a determinant due to Lindstr\"om, Gessel and Viennot
(see Theorem~\ref{gessel-viennot}), this bijection 
reduces the problem to evaluating a certain determinant with binomial 
entries, see Lemma~\ref{det}. 
In case that the fixed rhombus is situated in the centre or next to the centre 
we are able to 
evaluate the determinant, see Lemma~\ref{detlem} and Lemma~\ref{detlem2} in Section~\ref{comp}. We partly make us of a method 
that has already produced evaluations of other binomial determinants
(see e.g. \cite{FulKrat1}, \cite{kratdet}). In the course of evaluating the determinant 
we need an alternative expression for the number of rhombus tilings that we are interested in,  in form of the aforementioned  triple sum. This is given in 
Lemma~\ref{dreisum} in Section~\ref{dreisumsec}. Finally, in Section~\ref{asymsec} we provide the proofs of Theorem~\ref{asym:th} and Theorem~\ref{asym:th2}.
 
 \section{From rhombus tilings to non-intersecting lattice paths and determinants}
\label{bijection}

As mentioned before, the subject of this section is the 
bijection between non-intersecting lattice paths and rhombus tilings. 
In our context, non-intersecting lattice paths are 
disjoint paths in the lattice $\mathbb{Z}^2$ with steps in 
direction $(1,0)$ and $(0,-1)$. Figure~\ref{git} shows an example of such a family of 
non-intersecting lattice paths. It consists of three paths, each one connecting 
an initial point $A_i$ with a destination 
point $E_i$, $i=1,2,3$.

Figure~\ref{tilgit} and Figure~\ref{git} illustrate the bijection between rhombus tilings of 
a hexagon with side lengths $a$,$b$,$c$,$a$,$b$,$c$ and non-intersecting lattice 
paths with initial  points 
\begin{equation}
\label{start}
A_i=(i-1,c+i-1)
\end{equation}
and destination points
\begin{equation}
\label{end}
E_i=(b+i-1,i-1),
\end{equation}
$i=1,2,\dots,a$. 
Figure~\ref{tilgit} shows the rhombus tiling from Figure~\ref{fig:2} and an 
indication of the corresponding non-intersecting lattice paths. The initial points $A'_i$ and 
the destination points $E'_i$ of the paths in this figure are 
the centres of the sides of 
the rhombi that form the two sides of the hexagon with side lengths $a$. 
Roughly speaking these paths just 
describe `the way down' from $A'_i$ to $E'_i$ on the three dimensional 
pile of cubes --- elsewhere called {\em plane partitions} --- which the rhombus tiling, when interpreted as three-dimensional object, gives.

\begin{figure}
\centertexdraw{
\drawdim cm
% Tiling 
\linewd 0.03
\RhombusPos
\rmove(0.8660254037844 -0.5)
\RhombusNeg
\rmove(0 1)
\RhombusVert
\rmove(0.8660254037844 0.5)
\RhombusVert

\rmove(-0.8660254037844 -0.5)
\rmove(-0.8660254037844 -0.5)
\rmove(0 -1)

\RhombusPos
\rmove(0.8660254037844 -0.5)
\RhombusNeg
\rmove(0.8660254037844 0.5)
\rmove(0 1)
\RhombusPos
\rmove(0.8660254037844 0.5)
\RhombusPos

\rmove(-0.8660254037844 -0.5)
\rmove(-0.8660254037844 -0.5)
\rmove(-0.8660254037844 -0.5)
\rmove(0 -1)

\RhombusPos
\rmove(0.8660254037844 -0.5)
\RhombusNeg
\rmove(0.8660254037844 0.5)
\RhombusNeg
\rmove(0 1)
\RhombusVert
\rmove(0.8660254037844 0.5)
\RhombusVert
\rmove(0.8660254037844 0.5)
\RhombusNeg

\rmove(-0.8660254037844 -0.5)
\rmove(-0.8660254037844 -0.5)
\rmove(-0.8660254037844 -0.5)
\rmove(-0.8660254037844 -0.5)
\rmove(0 -1)

\RhombusVert
\rmove(0.8660254037844 0.5)
\rmove(0.8660254037844 -0.5)
\RhombusNeg
\rmove(0.8660254037844 0.5)
\rmove(0 1)
\RhombusPos
\rmove(0.8660254037844 -0.5)
\RhombusNeg
\rmove(0 1)
\RhombusVert
\rmove(0.8660254037844 0.5)
\RhombusNeg

\rmove(-0.8660254037844 -0.5)
\rmove(-0.8660254037844 -0.5)
\rmove(-0.8660254037844 -0.5)
\rmove(-0.8660254037844 -0.5)
\rmove(-0.8660254037844 -0.5)
\rmove(0 -1)

\RhombusNeg
\rmove(0.8660254037844 0.5)
\RhombusPos
\rmove(0.8660254037844 -0.5)
\RhombusNeg
\rmove(0.8660254037844 0.5)
\RhombusNeg
\rmove(0 1)
%Mitte
\rlvec(0.8660254037844 0.5)
\rlvec(0.8660254037844 -0.5)
\rlvec(-0.8660254037844 -0.5)
\lfill f:0.5
\rlvec(-0.8660254037844  0.5)
\rmove(0.8660254037844 0.5)
\rmove(0.8660254037844 0.5)
\RhombusPos
\rmove(0.8660254037844 0.5)
\RhombusNeg
\rmove(0.8660254037844 -0.5)
\RhombusNeg

\rmove(-0.8660254037844 -0.5)
\rmove(-0.8660254037844 -0.5)
\rmove(-0.8660254037844 -0.5)
\rmove(-0.8660254037844 -0.5)
\rmove(-0.8660254037844 -0.5)
\rmove(-0.8660254037844 -0.5)

\RhombusVert
\rmove(0.8660254037844 0.5)
\rmove(0.8660254037844 -0.5)
\RhombusNeg
\rmove(0.8660254037844 0.5)
\RhombusNeg
\rmove(0 1)
\RhombusVert
\rmove(0.8660254037844 0.5)
\RhombusVert
\rmove(0.8660254037844 0.5)
\RhombusNeg
\rmove(0.8660254037844 -0.5)
\RhombusNeg

\rmove(-0.8660254037844 -0.5)
\rmove(-0.8660254037844 -0.5)
\rmove(-0.8660254037844 -0.5)
\rmove(-0.8660254037844 -0.5)
\rmove(-0.8660254037844 -0.5)

\RhombusVert
\rmove(0.8660254037844 0.5)
\rmove(0.8660254037844 -0.5)
\RhombusNeg
\rmove(0.8660254037844 0.5)
\rmove(0 1)
\RhombusPos
\rmove(0.8660254037844 0.5)
\RhombusPos
\rmove(0.8660254037844 -0.5)
\RhombusNeg

\rmove(-0.8660254037844 -0.5)
\rmove(-0.8660254037844 -0.5)
\rmove(-0.8660254037844 -0.5)
\rmove(-0.8660254037844 -0.5)

\RhombusVert 
\rmove(0.8660254037844 0.5)
\rmove(0.8660254037844 0.5)
\RhombusPos
\rmove(0.8660254037844 0.5)
\RhombusVert

\rmove(-0.8660254037844 -0.5)
\rmove(-0.8660254037844 -0.5)
\rmove(0 -1)

\RhombusVert
\rmove(0.8660254037844 0.5)
\RhombusVert
\rmove(0.8660254037844 0.5)
\RhombusNeg
\rmove(0.8660254037844 0.5)
\RhombusPos

%Gitterpunktweg
\linewd 0.06
\lpatt (0.06 0.1)
\move(0 0)
\rmove(0.4330127018922 0.25)

\rlvec(0 -1)
\rlvec(0 -1)
\rlvec(0 -1)
\rlvec(0.8660254037844 -0.5)
\rlvec(0 -1)
\rlvec(0.8660254037844 -0.5)
\rlvec(0.8660254037844 -0.5)
\rlvec(0.8660254037844 -0.5)
\rlvec(0.8660254037844 -0.5)

\move(0 0)
\rmove(0.8660254037844 0.5)
\rmove(0.4330127018922 0.25)

\rlvec(0.8660254037844 -0.5)
\rlvec(0 -1)
\rlvec(0.8660254037844 -0.5)
\rlvec(0 -1)
\rlvec(0.8660254037844 -0.5)
\rlvec(0.8660254037844 -0.5)
\rlvec(0 -1)
\rlvec(0 -1)
\rlvec(0.8660254037844 -0.5)

\move(0 0)
\rmove(0.8660254037844 0.5)
\rmove(0.8660254037844 0.5)
\rmove(0.4330127018922 0.25)

\rlvec(0.8660254037844 -0.5)
\rlvec(0 -1)
\rlvec(0.8660254037844 -0.5)
\rlvec(0.8660254037844 -0.5)
\rlvec(0 -1)
\rlvec(0.8660254037844 -0.5)
\rlvec(0 -1)
\rlvec(0.8660254037844 -0.5)
\rlvec(0 -1) 
\htext(0 0.4){$A'_1$}
\htext(0.9 0.9){$A'_2$}
\htext(1.8 1.4){$A'_3$}
\htext(4.7 -6.7){$E'_1$}
\htext(5.6 -6.2){$E'_2$}
\htext(6.5 -5.7){$E'_3$}
\htext(-1 -7.5){\small A rhombus tiling of a hexagon with an indication}
\htext(-1 -8){\small of the corresponding family of 
non-intersecting lattice paths.} }
\caption{}
\label{tilgit}
\end{figure}  

In order to obtain the corresponding family of non-intersecting lattice 
paths in  Figure~\ref{git}, we only have to change the $120^\circ$ angles 
of the paths in Figure~\ref{tilgit} into right angles.

\begin{figure}
\begin{picture} (200,250)(0,-40)
\Einheit 1cm
\PfadDicke{1.3pt}
\Koordinatenachsen(8,7)(0,0)
\Pfad(0,4),  2 2 2 1 2 1 1 1 1 \endPfad
\Pfad(1,5),  1 2 1 2 1 1 2 2 1 \endPfad
\Pfad(2,6),  1 2 1 1 2 1 2 1 2 \endPfad
\NormalPunkt(0,4)
\NormalPunkt(1,5)
\NormalPunkt(2,6)
\NormalPunkt(5,0)
\NormalPunkt(6,1)
\NormalPunkt(7,2)
\Kreis(3,3)
\Kreis(4,3)
\put(-15,110) {$A_1$}
\put(11,140) {$A_2$}
\put(40,168) {$A_3$}
\put(135,-15) {$E_1$}
\put(164,13) {$E_2$}
\put(193,42) {$E_3$}
\put(60,70) {$(x-1,y)$}
\put(108,70) {$(x,y)$}
\put(-100,-45){\small The family of non-intersecting lattice paths corresponding to the rhombus tiling in 
Figure~\ref{tilgit}.}
\end{picture}
\caption{}
\label{git}
\end{figure}

Therefore we already have a bijection between all rhombus tilings of a hexagon with side lengths 
$a$,$b$,$c$,$a$,$b$,$c$ and all 
non-intersecting lattice paths starting in $A_i$ and stopping in $E_i$. 
However, we are interested in special rhombus tilings --- namely in those 
with one fixed rhombus. Because of that we have to look for a simple 
description of the non-intersecting lattice paths that correspond to the 
rhombus tilings with a fixed rhombus.
Studying Figure~\ref{tilgit} and Figure~\ref{git} again, we see that the marked rhombus 
in Figure~\ref{tilgit} corresponds to the edge $(x-1,y) \rightarrow (x,y)$ in
Figure~\ref{git} (in our example $x=4$ and $y=3$). Therefore the rhombus tilings which 
contain a fixed horizontal rhombus correspond to 
the families of non-intersecting lattice paths which contain a fixed horizontal edge. 
Thus, we have found a bijection between rhombus tilings of a hexagon 
with side lengths $a$,$b$,$c$,$a$,$b$,$c$ which contain a fixed 
horizontal rhombus and non-intersecting lattice paths with initial 
points $A_i$ and destination points $E_i$, $i=1,2,\dots,a$, which contain a fixed horizontal edge.
As we will see in a moment, the latter can be counted with the help of the following 
main theorem on non-intersecting lattice paths: 

\begin{theo}[Lindstr\"om-Gessel-Viennot]
\label{gessel-viennot}
Let $A_1,A_2, \dots , A_r$  and $E_1,E_2, \dots , E_r$ be points in $\mathbb{Z}^2$. 
Then the determinant
\begin{equation}
\label{gv:det}
\det_{1 \le i , j \le r} 
    \left( \left|P\left(A_i \rightarrow E_j\right) \right| \right)
\end{equation}
is equal to 
\begin{equation*}
\sum_{\sigma \in \mathcal{S}_r} \sgn \sigma \cdot 
\left|P^+\left(A \rightarrow E_\sigma \right) \right|,
\end{equation*} 
where $\left|P\left(A_i \rightarrow E_j \right) \right|$ denotes the number of 
lattice paths connecting $A_i$ to $E_j$, and 
$\left|P^+\left(A \rightarrow E_\sigma \right) \right|$ denotes the number 
of families $(P_1,P_2,\dots,P_r)$ of non-intersecting lattice 
paths, the $i$th path $P_i$ connecting $A_i$ to $E_{\sigma(i)}$, $1 \le i \le r$.   
\end{theo}

This theorem can be found in \cite[Theorem 1]{gesvie} or, alternatively,
in \cite[Lemma 1]{lin}.

\begin{rem}
\label{gitterweg}   \rm
Usually the theorem is applied to the following situation:
There exists only one permutation $\sigma$ (normally the identity permutation) 
such that the set of families
 of non-intersecting lattice 
paths connecting $A_i$ to $E_{\sigma(i)}$, $1 \le i \le r$, is not empty
(e.g., this is the case for the points $A_i$ and $E_i$ defined in 
\eqref{start} and \eqref{end}).
Then Theorem~\ref{gessel-viennot} solves the enumeration problem of 
counting non-intersecting lattice paths --- by means of the determinant \eqref{gv:det}  --- totally. This is because 
$\left|P\left(A \rightarrow E \right) \right|$ is easy to determine: 
Let $A=(a_1,a_2)$ and $E=(e_1,e_2)$ be two points in $\mathbb{Z}^2$, such 
that $A$ is located in the north-west of $E$ (i.e., $e_1 \ge a_1$ and 
$a_2 \ge e_2$). Then the number of lattice paths with 
steps in direction $(1,0)$ and $(0,-1)$ is 
\begin{equation}
\label{number:lp}
| P(A \rightarrow E)| = \binom{e_1-a_1+a_2-e_2}{e_1-a_1}
=\binom{e_1-a_1+a_2-e_2}{a_2-e_2},
\end{equation}
since each lattice paths corresponds to a choice of $e_1-a_1$ 
steps in direction $(1,0)$ out of the total number of 
$e_1-a_1+a_2-e_2$ steps.
\end{rem}

Remark~\ref{gitterweg} shows that Theorem~\ref{gessel-viennot} can be useful in  the enumeration of non-intersecting 
lattice paths with given initial and destination points. 
Since our enumeration problem of non-intersecting lattice paths not 
only involves fixed initial and destination points but also the additional condition about 
the fixed edge $(x-1,y) \rightarrow (x,y)$, Theorem~\ref{gessel-viennot}
seems useless in our situation at first glance. The following `trick' 
remedies this matter: We add the following pair of initial and 
destination points:
\begin{equation}
\label{zus}
A_{a+1}=(x,y) \qquad \text{and} \qquad E_{a+1}=(x-1,y)
\end{equation}
(see Figure~\ref{git}). Then, as is not difficult to see, the number of rhombus tilings with a fixed rhombus 
`at' $(x,y)$ equals the number of families of non-intersecting lattice paths with 
initial points $A_i$ and destination points $E_i$, $i=1,2, \dots , a+1$.

Now we are ready to apply Theorem~\ref{gessel-viennot} to the points $A_i$ and $E_i$, $i=1,2,
\dots , a+1$, defined in \eqref{start}, \eqref{end} and 
\eqref{zus}. In order to figure out what the determinant \eqref{gv:det} actually 
gives, we have to find the  
permutations $\sigma$ for which there exist non-intersecting 
lattice paths connecting $A_i$ to $E_{\sigma(i)}$ with $1 \le i \le a+1$. 
It is quite easy to see that this is only accomplished by
transpositions of the form $(i, a+1)$, where $i \not= a+1$.
The fact that all transpositions have the same sign --- namely $-1$ --- 
implies that the determinant in \eqref{gv:det} applied to our special 
points  gives the number of families of non-intersecting lattice paths with initial 
points $A_i$ and destination points $E_i$, $1 \le i \le a+1$, with negative sign. 

We use \eqref{number:lp} to compute $\left|P \left(A_i \rightarrow E_j \right) \right|$ in 
\eqref{gv:det},  and finally obtain the following: 

\begin{lem}
\label{det}
Let $a$,$b$ and $c$ be positive integers and $(x,y)$ be an integer point such 
that $0 \le x \le b+a-1$ and $1 \le y \le c+a-1$. Then the number of rhombus tilings 
of a hexagon with side lengths $a$,$b$,$c$,$a$,$b$,$c$ which contain the fixed 
horizontal rhombus that corresponds to the point $(x,y)$ in the bijection described 
above equals

\unitlength1cm
\begin{picture}(0,3.5)
\put(-1,1.5){\parbox{3cm}{
\begin{equation} 
\label{matrix:allg:1}
-\det_{1 \le i, j \le a+1} \left(
\begin{array}{ccc|c}
      &                    & &   \\
      &\binom{b+c}{c-i+j}  & & \binom{b-x+y}{y-i+1} \\
      &                    & &   \\ \hline
      &\binom{c+x-y-1}{x-j}& & 0 
\end{array}
\right)
\end{equation}}}
\put(5.7,1.8){\makebox(2,1)[t]{ \scriptsize $1 \le j \le a$}}
\put(8,1.8){\makebox(2,1)[t]{ \scriptsize $j=a+1$}}
\put(9.7,0.9){\makebox(2,1)[t]{ \scriptsize $1 \le i \le a$}}
\put(9.7,-0.1){\makebox(2,1)[t]{ \scriptsize $i=a+1$}}
\put(10.7, 0.5){\makebox(2,1)[t]{.}}
\end{picture}
\end{lem}

In Lemma~\ref{det} we refer to a correspondence between the integer points $(x,y)$, 
$0 \le x \le b+a-1$ and $1 \le y \le c+a-1$, and the horizontal rhombi we use 
for the rhombus tilings of a hexagon with side lengths $a$,$b$,$c$,$a$,$b$,$c$. 
This correspondence is implicitly given by the bijection between 
rhombus tilings and non-intersecting lattice paths described above. The following 
remark makes the bijection between points and rhombi more 
explicit.

\begin{rem} 
\label{cor}
\rm 
Again we consider a hexagon with side lengths $a$,$b$,$c$,$a$,$b$,$c$. 
We introduce the following oblique angled coordinate system: Its 
origin is located in one of the two vertices, where sides of lengths $b$ and 
$c$ meet, and the axes are induced by those two sides (see Figure~\ref{fig:10}). 
The units are chosen 
such that the (Euclidean) side lengths of the considered hexagon 
are $a$,$b$,$c$,$a$,$b$,$c$  in this coordinate system, too.
(That is to say, the two triangles in Figure~\ref{fig:10} with 
vertices in the origin form  the unit `square').
Thus, in this coordinate system, the points in Figure~\ref{tilgit} have coordinates 
$A'_1=(1/2,9/2)$, $A'_2=(3/2,11/2)$, $A'_3=(5/2,13/2)$, 
$E'_1=(11/2,1/2)$, $E'_2=(13/2,3/2)$ and  $E'_3=(15/2,5/2)$.
Phrased differently, the coordinate system  was chosen such that 
$A'_i=A_i+(1/2,1/2)$ and $E'_i=E_i+(1/2,1/2)$, where $A'_i$ and $E'_i$ 
are the initial and destination points of the paths in 
Figure~\ref{tilgit}, and $A_i$ and $E_i$ are defined as in 
\eqref{start}, \eqref{end} and \eqref{zus}, $1 \le i \le a+1$. From this 
point of view, the family of paths in Figure~\ref{tilgit} is the family of 
paths in Figure~\ref{git} translated by $(1/2,1/2)$ and drawn in the 
oblique angled coordinate system. Furthermore, we see that the integer 
point $(x,y)$ in the oblique angled coordinate system  is just the lowest vertex of 
the corresponding rhombus under the aforementioned bijection between points and rhombi
(see Figure~\ref{tilgit} and Figure~\ref{git}).
\end{rem}

\begin{figure}
\centertexdraw{
\drawdim cm
\arrowheadtype t:F
\linewd 0.05
\move(0 -3)
\ravec(0 5)
\move(0 -3)
\rlvec(0.8660254037844 -0.5)
\rlvec(0 1)
\rlvec(-0.866025403784 0.5)
\rlvec(0 -1)
\rlvec(0.8660254037844 -0.5)
\rlvec(0.8660254037844 -0.5)
\rlvec(0.8660254037844 -0.5)
\rlvec(0.8660254037844 -0.5)
\rlvec(0.8660254037844 -0.5)
\ravec(0.8660254037844 -0.5)
\linewd 0.01
\move(0 0)
\DreieckBreit
\rmove(0 1)
\DreieckSpitz
\rmove(0.8660254037844 -0.5)
\DreieckBreit
\rmove(0 1)
\DreieckSpitz
\rmove(0.8660254037844 -0.5)
\DreieckBreit
\rmove(0 1)
\DreieckSpitz
\rmove(0.8660254037844 -0.5)
\DreieckBreit
\rmove(-0.8660254037844 -0.5)
\rmove(-0.8660254037844 -0.5)
\rmove(-0.8660254037844 -0.5)
\rmove(0 -1)
\DreieckBreit
\rmove(0 1)
\DreieckSpitz
\rmove(0.8660254037844 -0.5)
\DreieckBreit
\rmove(0 1)
\DreieckSpitz
\rmove(0.8660254037844 -0.5)
\DreieckBreit
\rmove(0 1)
\DreieckSpitz
\rmove(0.8660254037844 -0.5)
\DreieckBreit
\rmove(0 1)
\DreieckSpitz
\rmove(0.8660254037844 -0.5)
\DreieckBreit
\rmove(-0.8660254037844 -0.5)
\rmove(-0.8660254037844 -0.5)
\rmove(-0.8660254037844 -0.5)
\rmove(-0.8660254037844 -0.5)
\rmove(0 -1)
\DreieckBreit
\rmove(0 1)
\DreieckSpitz
\rmove(0.8660254037844 -0.5)
\DreieckBreit
\rmove(0 1)
\DreieckSpitz
\rmove(0.8660254037844 -0.5)
\DreieckBreit
\rmove(0 1)
\DreieckSpitz
\rmove(0.8660254037844 -0.5)
\DreieckBreit
\rmove(0 1)
\DreieckSpitz
\rmove(0.8660254037844 -0.5)
\DreieckBreit
\rmove(0 1)
\DreieckSpitz
\rmove(0.8660254037844 -0.5)
\DreieckBreit
\rmove(-0.8660254037844 -0.5)
\rmove(-0.8660254037844 -0.5)
\rmove(-0.8660254037844 -0.5)
\rmove(-0.8660254037844 -0.5)
\rmove(-0.8660254037844 -0.5)
\rmove(0 -1)
\DreieckBreit
\rmove(0 1)
\DreieckSpitz
\rmove(0.8660254037844 -0.5)
\DreieckBreit
\rmove(0 1)
\DreieckSpitz
\rmove(0.8660254037844 -0.5)
\DreieckBreit
\rmove(0 1)
\DreieckSpitz
\rmove(0.8660254037844 -0.5)
\DreieckBreit
\rmove(0 1)
\DreieckSpitz
\rmove(0.8660254037844 -0.5)
\DreieckBreit
\rmove(0 1)
\DreieckSpitz
\rmove(0.8660254037844 -0.5)
\DreieckBreit
\rmove(0 1)
\DreieckSpitz
\rmove(0.8660254037844 -0.5)
\DreieckBreit
\rmove(-0.8660254037844 -0.5)
\rmove(-0.8660254037844 -0.5)
\rmove(-0.8660254037844 -0.5)
\rmove(-0.8660254037844 -0.5)
\rmove(-0.8660254037844 -0.5)
\rmove(-0.8660254037844 -0.5)
\DreieckSpitz
\rmove(0.8660254037844 -0.5)
\DreieckBreit
\rmove(0 1)
\DreieckSpitz
\rmove(0.8660254037844 -0.5)
\DreieckBreit
\rmove(0 1)
\DreieckSpitz
\rmove(0.8660254037844 -0.5)
\DreieckBreit
\rmove(0 1)
\DreieckSpitz
\rmove(0.8660254037844 -0.5)
\DreieckBreit
\rmove(0 1)
\DreieckSpitz
\rmove(0.8660254037844 -0.5)
\DreieckBreit
\rmove(0 1)
\DreieckSpitz
\rmove(0.8660254037844 -0.5)
\DreieckBreit
\rmove(0 1)
\DreieckSpitz
\rmove(0.8660254037844 -0.5)
\DreieckBreit
\rmove(-0.8660254037844 -0.5)
\rmove(-0.8660254037844 -0.5)
\rmove(-0.8660254037844 -0.5)
\rmove(-0.8660254037844 -0.5)
\rmove(-0.8660254037844 -0.5)
\rmove(-0.8660254037844 -0.5)
\DreieckSpitz
\rmove(0.8660254037844 -0.5)
\DreieckBreit
\rmove(0 1)
\DreieckSpitz
\rmove(0.8660254037844 -0.5)
\DreieckBreit
\rmove(0 1)
\DreieckSpitz
\rmove(0.8660254037844 -0.5)
\DreieckBreit
\rmove(0 1)
\DreieckSpitz
\rmove(0.8660254037844 -0.5)
\DreieckBreit
\rmove(0 1)
\DreieckSpitz
\rmove(0.8660254037844 -0.5)
\DreieckBreit
\rmove(0 1)
\DreieckSpitz
\rmove(0.8660254037844 -0.5)
\DreieckBreit
\rmove(0 1)
\DreieckSpitz
\rmove(0.8660254037844 -0.5)
\rmove(-0.8660254037844 -0.5)
\rmove(-0.8660254037844 -0.5)
\rmove(-0.8660254037844 -0.5)
\rmove(-0.8660254037844 -0.5)
\rmove(-0.8660254037844 -0.5)
\rmove(-0.8660254037844 -0.5)
\DreieckSpitz
\rmove(0.8660254037844 -0.5)
\DreieckBreit
\rmove(0 1)
\DreieckSpitz
\rmove(0.8660254037844 -0.5)
\DreieckBreit
\rmove(0 1)
\DreieckSpitz
\rmove(0.8660254037844 -0.5)
\DreieckBreit
\rmove(0 1)
\DreieckSpitz
\rmove(0.8660254037844 -0.5)
\DreieckBreit
\rmove(0 1)
\DreieckSpitz
\rmove(0.8660254037844 -0.5)
\DreieckBreit
\rmove(0 1)
\DreieckSpitz
\rmove(0.8660254037844 -0.5)
\rmove(-0.8660254037844 -0.5)
\rmove(-0.8660254037844 -0.5)
\rmove(-0.8660254037844 -0.5)
\rmove(-0.8660254037844 -0.5)
\rmove(-0.8660254037844 -0.5)
\DreieckSpitz
\rmove(0.8660254037844 -0.5)
\DreieckBreit
\rmove(0 1)
\DreieckSpitz
\rmove(0.8660254037844 -0.5)
\DreieckBreit
\rmove(0 1)
\DreieckSpitz
\rmove(0.8660254037844 -0.5)
\DreieckBreit
\rmove(0 1)
\DreieckSpitz
\rmove(0.8660254037844 -0.5)
\DreieckBreit
\rmove(0 1)
\DreieckSpitz
\rmove(0.8660254037844 -0.5)
\rmove(-0.8660254037844 -0.5)
\rmove(-0.8660254037844 -0.5)
\rmove(-0.8660254037844 -0.5)
\rmove(-0.8660254037844 -0.5)
\DreieckSpitz
\rmove(0.8660254037844 -0.5)
\DreieckBreit
\rmove(0 1)
\DreieckSpitz
\rmove(0.8660254037844 -0.5)
\DreieckBreit
\rmove(0 1)
\DreieckSpitz
\rmove(0.8660254037844 -0.5)
\DreieckBreit
\rmove(0 1)
\DreieckSpitz
\rmove(0.8660254037844 -0.5)
\move(-0.5 -1)
\htext{$c$}
\move(1.299038105677 1.75)
\rmove(0  0.5)
\htext{$a$}
\rmove(0 -0.5)
\rmove(1.299038105677 0.75)
\rmove(1.299038105677 -0.75)
\rmove(0.8660254037844 -0.5)
\rmove(0 0.5)
\htext{$b$}
\rmove(0 -0.5)
\rmove(1.299038105677 -0.75)
\rmove(0.8660254037844 -0.5)
\rmove(0 -2)
\rmove(0.5 0)
\htext{$c$}
\rmove(-0.5 0)
\rmove(0 -2)
\rmove(-1.299038105677 -0.75)
\rmove(0 -0.5)
\htext{$a$}
\rmove(0.5 0)
\rmove(-1.299038105677 -0.75)
\rmove(-1.299038105677 0.75)
\rmove(-0.8660254037844 0.5)
\rmove(-0.8660254037844 0.5)
\rmove(0 -0.5)
\htext{$b$}
\htext(-1 -7){The oblique angled coordinate system.}
}
\caption{}
\label{fig:10}
\end{figure}

\section{From the determinant to a triple sum}
\label{dreisumsec}

The aim of this section is the derivation of a triple sum 
that is equal to \eqref{matrix:allg:1} and therefore gives the number of all
rhombus tilings of a hexagon with side lengths $a$,$b$,$c$,$a$,$b$,$c$ which contain a fixed 
rhombus with lowest vertex $(x,y)$: 

\begin{lem} 
\label{dreisum}
Let $a$, $b$ and $c$ be positive integers, and let $(x,y)$ be an integer point such that
$0 \le x \le b+a-1$ and $1 \le y \le c+a-1$. Then the number of rhombus tilings of a hexagon 
with side lengths 
$a$,$b$,$c$,$a$,$b$,$c$ which contain the  fixed horizontal rhombus with lowest vertex 
$(x,y)$ (in the oblique angled coordinate system; see Remark~\ref{cor}) equals
\begin{multline}
\label{drei}
\left( \prod_{i=1}^{a} \prod_{j=1}^b \prod_{k=1}^c \frac{i+j+k-1}{i+j+k-2}  \right)
\frac{(1)_c}{(b+1)_c}\\
\times
\sum_{n=1}^{a} \sum_{m=1}^{a} \sum_{s=1}^{m}  \left[ (-1)^{n+s} 
\binom{c+x-y+n-2}{x-1} \binom{b-x+y+s-1}{b-x+s-1} \frac{(b+1)_{s-1}}{(b+c+1)_{s-1}}  \right. \\
\left.
\binom{m-1}{s-1} \frac{(c+1)_{n-1}}{(n-1)!} \frac{(b+c+n)_{m-n}}{(m-n)!}  \right] .
\end{multline}
\end{lem}

{\bf Outline of the proof of Lemma~\ref{dreisum}.} We derive the triple sum by starting from \eqref{matrix:allg:1}.
The matrix underlying the determinant in \eqref{matrix:allg:1} 
has a `homogeneous' definition, except for 
the last row and the last column. Our proof starts with  
some elementary row and column operations that transform the 
`homogeneous' submatrix into a matrix of triangular form (see \eqref{matrix:allg:8}). 
Next we expand the determinant along the exceptional row and 
along the exceptional column (see \eqref{expand}). The result is a triple sum with 
a summand that involves another determinant. But this determinant is 
the determinant of the aforementioned 
triangular matrix with one row and one column deleted (which 
row and which column depends on the summation index of the triple sum, see \eqref{matrix:2}) 
and therefore we are able to compute it: 
The triangle property allows an expansion of this remaining determinant along the 
first $(n-1)$ columns and the last $(a-m)$ rows, where $m$ denotes 
the missing row and $n$ denotes the missing column (see \eqref{expand:2}). Finally the 
remaining $(m-n) \times (m-n)$ determinants can be reduced to Vandermonde's determinant (see 
\eqref{matrix:3}).  

\medskip

{\bf Proof of Lemma~\ref{dreisum} -- The details.}
As mentioned before, we first describe some elementary row and column operations 
that transform the `homogeneous' $a \times a$ submatrix $\left(\binom{b+c}{c-i+j}\right)_{1 \le i,j \le a}$ 
of the matrix underlying the determinant in \eqref{matrix:allg:1} into an upper 
triangular matrix. 

We begin with some elementary column operations:
By using the symmetry of the binomial coefficient (i.e., the identity
$\binom{n}{k}=\binom{n}{n-k}$) in the last row of the 
matrix underlying the determinant in \eqref{matrix:allg:1}  we observe that 
the determinant in \eqref{matrix:allg:1} is equal to

\unitlength1cm
\begin{picture}(20,3.5)
\put(-0.5,1.5){\parbox{3cm}{
\begin{equation} 
\label{matrix:allg:2}
-\det_{1 \le i, j \le a+1} \left(
\begin{array}{ccc|c}
      &                    & &   \\
      &\binom{b+c}{c-i+j}  & & \binom{b-x+y}{y-i+1} \\
      &                    & &   \\ \hline
      &\binom{c+x-y-1}{c-y+j-1}& & 0 
\end{array}
\right)
\end{equation}}}
\put(6.2,1.8){\makebox(2,1)[t]{ \scriptsize $1 \le j \le a$}}
\put(8.5,1.8){\makebox(2,1)[t]{ \scriptsize $j=a+1$}}
\put(10.2,0.9){\makebox(2,1)[t]{ \scriptsize $1 \le i \le a$}}
\put(10.2,-0.1){\makebox(2,1)[t]{ \scriptsize $i=a+1$}}
\put(11.2, 0.5){\makebox(2,1)[t]{.}}
\end{picture}

We add the $(j-1)$th column to the $j$th column, $j=a,a-1,\dots,2$, in that order. 
The entries of the changed matrix read as 
\begin{equation*}
\binom{b+c}{c-i+j} + \binom{b+c}{c-i+j-1} = \binom{b+c+1}{c-i+j}
\end{equation*}
for $i=1,2, \dots,a$ and $j=2,3, \dots, a$, and 
\begin{equation*}
\binom{c+x-y-1}{c-y+j-1}+\binom{c+x-y-1}{c-y+j-2}=
\binom{c+x-y}{c-y+j-1}
\end{equation*}
for $i=a+1$ and $j=2,3, \dots, a$. The other entries do not change. Thus, 
the following determinant is equal to the determinant in \eqref{matrix:allg:2}:

\unitlength1cm
\begin{picture}(0,3.5)
\put(-0.5,1.5){\parbox{3cm}{
\begin{equation}
\label{allg:3}
-\det_{1 \le i, j \le a+1} \left(
\begin{array}{c|ccc|c}
                   &     &                      & &   \\
 \binom{b+c}{c-i+1}&     &\binom{b+c+1}{c-i+j}  & & \binom{b-x+y}{y-i+1} \\
                       &     &                      & &   \\ \hline
 \binom{c+x-y-1}{c-y} &     &\binom{c+x-y}{c-y+j-1}& & 0 
\end{array}
\right)
\end{equation}}}
\put(5,1.8){\makebox(2,1)[t]{ \scriptsize $j=1$}}
\put(7.3,1.8){\makebox(2,1)[t]{ \scriptsize $2 \le j \le a$}}
\put(9.5,1.8){\makebox(2,1)[t]{ \scriptsize $j=a+1$}}
\put(11.3,0.9){\makebox(2,1)[t]{ \scriptsize $1 \le i \le a$}}
\put(11.3,-0.1){\makebox(2,1)[t]{ \scriptsize $i=a+1$}}
\put(12, 0.5){\makebox(2,1)[t]{.}}
\end{picture}

Next we repeat the procedure, i.e., we add the $(j-1)$th column to the $j$th column, 
$j=a,a-1, \dots,3$, in that order. Thus, we obtain that 
the following determinant  is equal to the determinant in \eqref{allg:3}:

\unitlength1cm
\begin{picture}(0,3.5)
\put(-1,1.5){\parbox{3cm}{
\begin{equation*} 
-\det_{1 \le i, j \le a+1} \left(
\begin{array}{c|c|ccc|c}
                   &  &    &                      & &   \\
 \binom{b+c}{c-i+1}& \binom{b+c+1}{c-i+2} &    &\binom{b+c+2}{c-i+j}  & & \binom{b-x+y}{y-i+1} \\
                    &   &     &                      & &   \\ \hline
 \binom{c+x-y-1}{c-y}& \binom{c+x-y}{c-y+1}&    &\binom{c+x-y+1}{c-y+j-1}& & 0 
\end{array}
\right)
\end{equation*}}}
\put(3.5,1.8){\makebox(2,1)[t]{ \scriptsize $j=1$}}
\put(5.3,1.8){\makebox(2,1)[t]{ \scriptsize $j=2$}}
\put(7.7,1.8){\makebox(2,1)[t]{ \scriptsize $3 \le j \le a$}}
\put(10,1.8){\makebox(2,1)[t]{ \scriptsize $j=a+1$}}
\put(11.8,0.9){\makebox(2,1)[t]{ \scriptsize $1 \le i \le a$}}
\put(11.8,-0.1){\makebox(2,1)[t]{ \scriptsize $i=a+1$}}
\put(12.5, 0.5){\makebox(2,1)[t]{.}}
\end{picture}

We repeat this procedure of adding successive  columns from right to left each 
time stopping one column earlier than before, as long as possible. This means that the procedure
is performed $a-1$ times including the two steps described in detail. 
Since the upper 
parameter of the binomial entry increases by one every time it is involved and 
every entry in the $j$th column participates $j-1$ times exactly,
$1 \le j \le a$, this 
procedure yields the following determinant:

\unitlength1cm
\begin{picture}(0,3.5)
\put(-0.5,1.5){\parbox{3cm}{
\begin{equation} 
\label{matrix:allg:3}
-\det_{1 \le i, j \le a+1} \left(
\begin{array}{ccc|c}
      &                    & &   \\
      &\binom{b+c+j-1}{c-i+j}  & & \binom{b-x+y}{y-i+1} \\
      &                    & &   \\ \hline
      &\binom{c+x-y+j-2}{c-y+j-1}& & 0 
\end{array}
\right)
\end{equation}}}
\put(6.2,1.8){\makebox(2,1)[t]{ \scriptsize $1 \le j \le a$}}
\put(8.5,1.8){\makebox(2,1)[t]{ \scriptsize $j=a+1$}}
\put(10.5,0.9){\makebox(2,1)[t]{ \scriptsize $1 \le i \le a$}}
\put(10.5,-0.1){\makebox(2,1)[t]{ \scriptsize $i=a+1$}}
\put(11.4, 0.5){\makebox(2,1)[t]{.}}
\end{picture}

Now we apply some elementary row operations to  \eqref{matrix:allg:3}. In fact these row 
operations are analogous to the column operations we just applied to 
\eqref{matrix:allg:2}. In order to do so, we first use the symmetry of the 
binomial coefficient for the first $a$ rows of the  
determinant in \eqref{matrix:allg:3} and observe that the determinant

\unitlength1cm
\begin{picture}(0,3.5)
\put(-0.5,1.5){\parbox{3cm}{
\begin{equation} 
\label{matrix:allg:4}
-\det_{1 \le i, j \le a+1} \left(
\begin{array}{ccc|c}
      &                    & &   \\
      &\binom{b+c+j-1}{b+i-1}  & & \binom{b-x+y}{b+i-x-1} \\
      &                    & &   \\ \hline
      &\binom{c+x-y+j-2}{c-y+j-1}& & 0 
\end{array}
\right)
\end{equation}}}
\put(6.2,1.8){\makebox(2,1)[t]{ \scriptsize $1 \le j \le a$}}
\put(8.5,1.8){\makebox(2,1)[t]{ \scriptsize $j=a+1$}}
\put(10.7,0.9){\makebox(2,1)[t]{ \scriptsize $1 \le i \le a$}}
\put(10.7,-0.1){\makebox(2,1)[t]{ \scriptsize $i=a+1$}}
\end{picture}

{ \parindent0cm
is equal to the determinant in \eqref{matrix:allg:3}.
As announced, we now add the $(i-1)$th row of the 
determinant in \eqref{matrix:allg:4} to the $i$th row, starting at $i=a$ and 
stopping at $i=2$. Next we do the same with the resulting determinant, starting 
at $i=a$ but stopping at $i=3$. After repeating this procedure 
$a-1$ times we obtain that the following determinant is equal to the determinant in 
\eqref{matrix:allg:4}:} 

\unitlength1cm
\begin{picture}(0,3.5)
\put(-0.5,1.5){\parbox{3cm}{
\begin{equation} 
\label{matrix:allg:5}
-\det_{1 \le i, j \le a+1} \left(
\begin{array}{ccc|c}
      &                    & &   \\
      &\binom{b+c+i+j-2}{b+i-1}  & & \binom{b-x+y+i-1}{b+i-x-1} \\
      &                    & &   \\ \hline
      &\binom{c+x-y+j-2}{c-y+j-1}& & 0 
\end{array}
\right)
\end{equation}}}
\put(6.2,1.8){\makebox(2,1)[t]{ \scriptsize $1 \le j \le a$}}
\put(8.5,1.8){\makebox(2,1)[t]{ \scriptsize $j=a+1$}}
\put(11,0.9){\makebox(2,1)[t]{ \scriptsize $1 \le i \le a$}}
\put(11,-0.1){\makebox(2,1)[t]{ \scriptsize $i=a+1$}}
\put(11.9, 0.5){\makebox(2,1)[t]{.}}
\end{picture}

Now we take the factor $(-1)^{b+c+i-1}(b+c+i-1)!/(b+i-1)!$ out 
of the $i$th row of the determinant in 
\eqref{matrix:allg:5}, $i=1,2, \dots, a$. This yields

\unitlength1cm
\begin{picture}(0,4.5)
\put(-0.5,2.5){\parbox{3cm}{
\begin{multline} 
\label{matrix:allg:6}
- \prod_{i=1}^{a} (-1)^{b+c+i-1} \frac{(b+c+i-1)!}{(b+i-1)!} \\
\times \det_{1 \le i, j \le a+1} \left(
\begin{array}{ccc|c}
      &                    & &   \\
      &\frac{(j-1)!}{(c+j-1)!} \binom{-j}{b+c+i-1}  & & (-1)^{b+c+i-1}\frac{(b+i-1)!}{(b+c+i-1)!}\binom{b-x+y+i-1}{b+i-x-1} \\
      &                    & &   \\ \hline
      &\binom{c+x-y+j-2}{c-y+j-1}& & 0       
\end{array}
\right) \\
\end{multline}}}
\put(4.7,2.4){\makebox(2,1)[t]{ \scriptsize $1 \le j \le a$}}
\put(8.9,2.4){\makebox(2,1)[t]{ \scriptsize $j=a+1$}}
\put(13.2,1.3){\makebox(2,1)[t]{ \scriptsize $1 \le i \le a$}}
\put(13.2,0.3){\makebox(2,1)[t]{ \scriptsize $i=a+1$}}
\put(13.9, 1){\makebox(2,1)[t]{.}}
\end{picture}

Finally we want to apply the elementary row operations we just applied to the 
determinant in \eqref{matrix:allg:4} once again. 
An analysis of the elementary row operations that  we applied to the  
determinant in \eqref{matrix:allg:4} yields that the  
determinant in \eqref{matrix:allg:5} was produced from 
the determinant in \eqref{matrix:allg:4} by replacing the $i$th row by 
\begin{equation*}
\sum_{s=1}^{i} \binom{i-1}{s-1} A(s),
\end{equation*}
$1 \le i \le a$, where $A(s)$ denotes the $s$th row of the 
determinant in \eqref{matrix:allg:4}.
We perform this replacement of the entries in the determinant 
in \eqref{matrix:allg:6} and obtain 

\unitlength1cm
\begin{picture}(0,4.5)
\put(-0.5,2.5){\parbox{3cm}{
\begin{multline} 
\label{matrix:allg:7}
- \prod_{i=1}^{a} (-1)^{b+c+i-1} \frac{(b+c+i-1)!}{(b+i-1)!} \\
\times \det_{1 \le i, j \le a+1} \left(
\begin{array}{c|c}
      \frac{(j-1)!}{(c+j-1)!}\binom{i-j-1}{b+c+i-1}  & 
      \sum_{s=1}^{i} \binom{i-1}{s-1} 
      (-1)^{b+c+s-1}\frac{(b+s-1)!}{(b+c+s-1)!}\binom{b-x+y+s-1}{b+s-x-1}  \\ \hline
      \binom{c+x-y+j-2}{c-y+j-1}& 0 
\end{array}
\right). \\
\end{multline}}}
\put(3.2,0.4){\makebox(2,1)[t]{ \scriptsize $1 \le j \le a$}}
\put(8.4,0.4){\makebox(2,1)[t]{ \scriptsize $j=a+1$}}
\put(13.5,1.7){\makebox(2,1)[t]{ \scriptsize $1 \le i \le a$}}
\put(13.5,0.8){\makebox(2,1)[t]{ \scriptsize $i=a+1$}}
\end{picture}

The entry in the $i$th row and $j$th column, $1 \le i \le a$ and 
$1 \le j \le a$,  develops from the corresponding entry in \eqref{matrix:allg:6}
by using Vandermonde's summation formula: 
\begin{align*}
\sum_{s=1}^{i} \binom{i-1}{s-1} \frac{(j-1)!}{(c+j-1)!} \binom{-j}{b+c+s-1}=
\frac{(j-1)!}{(c+j-1)!} \binom{i-j-1}{b+c+i-1}.
\end{align*}
Next we apply the elementary identity
\begin{equation}
\label{bin:id}
\binom{n}{k} = \binom{-n+k-1}{k} (-1)^k
\end{equation}
to this entry in 
the $i$th row and $j$th column of the matrix underlying the determinant 
in \eqref{matrix:allg:7}, $1 \le i \le a$ and $1 \le j \le a$, and then use the
symmetry of the binomial coefficient. 
Finally we take the factor $(-1)^{b+c+i-1}$ out of the $i$th row, $1 \le i \le a$, 
and obtain  that the expression in \eqref{matrix:allg:7} is equal 
to 

\unitlength1cm
\begin{picture}(0,4.5)
\put(-0.5,2.5){\parbox{3cm}{
\begin{multline} 
\label{matrix:allg:8}
- \prod_{i=1}^{a} \frac{(b+c+i-1)!}{(b+i-1)!} \\
\times \det_{1 \le i, j \le a+1} \left(
\begin{array}{c|c}
      \frac{(j-1)!}{(c+j-1)!}\binom{b+c+j-1}{j-i}  & 
      \sum_{s=1}^{i} \binom{i-1}{s-1} 
      (-1)^{i-s}\frac{(b+s-1)!}{(b+c+s-1)!}\binom{b-x+y+s-1}{b+s-x-1}  \\ \hline
      \binom{c+x-y+j-2}{c-y+j-1}& 0 
\end{array}
\right). \\
\end{multline}}}
\put(3.2,2.1){\makebox(2,1)[t]{ \scriptsize $1 \le j \le a$}}
\put(8.4,2.1){\makebox(2,1)[t]{ \scriptsize $j=a+1$}}
\put(13.2,1.7){\makebox(2,1)[t]{ \scriptsize $1 \le i \le a$}}
\put(13.2,0.8){\makebox(2,1)[t]{ \scriptsize $i=a+1$}}
\end{picture}

With pleasure we discover that the $a \times a$ submatrix induced by the first $a$ rows and 
first $a$ columns  of the matrix underlying the determinant 
in \eqref{matrix:allg:8} is an upper triangular matrix. This is because  the binomial coefficient 
$\binom{\alpha}{k}$ is defined to be zero if $k < 0$ for any indeterminante $\alpha$. 
Therefore the determinant is of that form we were looking 
for at the beginning of the proof. 

Next we expand \eqref{matrix:allg:8} along the last row and then along the last 
column, to obtain
\begin{multline}
\label{expand} 
-\prod_{i=1}^{a} \frac{(b+c+i-1)!}{(b+i-1)!} \sum_{n=1}^{a} (-1)^{a+1+n} 
  \binom{c+x-y+n-2}{c-y+n-1} \\
  \quad \times \sum_{m=1}^{a} (-1)^{a+m}
\sum_{s=1}^{m} \binom{m-1}{s-1}
(-1)^{m-s} \frac{(b+s-1)!}{(b+c+s-1)!} \binom{b-x+y+s-1}{b+s-x-1} \\
\quad \times 
\underset {i \not= m, j \not= n} {\det_{1 \le i, j \le a} }
\left( \frac{(j-1)!}{(c+j-1)!} \binom{b+c+j-1}{j-i} \right).
\end{multline}

Then we take the factor $(j-1)!/(c+j-1)!$ out of the $j$th column of the remaining determinant. 
This, together with some other manipulations, gives
\begin{multline}
\label{kA}
 \prod_{i=1}^{a}\frac{(b+c+i-1)!(i-1)!}{(b+i-1)!(c+i-1)!}  \\
\times \sum_{n=1}^{a} \sum_{m=1}^{a} \sum_{s=1}^{m} (-1)^{n+s}
\binom{c+x-y+n-2}{x-1} \binom{b-x+y+s-1}{b+s-x-1} \binom{m-1}{s-1}  \\
 \times  \frac{(b+s-1)!}{(b+c+s-1)!}
\frac{(c+n-1)!}{(n-1)!} 
\underset {i \not= m, j \not= n }{\det_{1 \le i, j \le a}}
\left( \binom{b+c+j-1}{j-i} \right).
\end{multline}
  
Now we have to compute 
\begin{equation} 
\label{matrix:2}
\underset{i \not= m, j \not= n}{\det_{1 \le i, j \le a}}
\left( \binom{b+c+j-1}{j-i} \right).
\end{equation}
This is a determinant of an upper triangular matrix, where the $m$th row and $n$th column was 
deleted. Therefore it is only different from zero if $n \le m$. 
So, let us assume $n \le m$. Expansion of the determinant along the first 
$(n-1)$ columns and along the last $(a-m)$ rows yields  
\begin{equation}
\label{expand:2}
\underset{n \le i \le m-1}{\det_{n+1 \le j \le m}}
\left( \binom{b+c+j-1}{j-i} \right)=
\underset{n \le i \le  m-1}{\det_{n \le j \le m-1}}
\left( \binom{b+c+j}{j-i+1} \right)
\end{equation}
for the determinant in \eqref{matrix:2}.

We are going  to reduce this determinant to Vandermonde's determinant. 
In order to do so, we take $(b+c+j)!/(j-n+1)!$ out of the $j$th column, $1 \le j \le a$, 
and $1/(b+c+i-1)!$ out of the $i$th row, $1 \le i \le a$. Thus, we obtain
that the determinant in \eqref{expand:2} is equal to
\begin{equation}  
\label{matrix:3}
(b+c+n)_{m-n}  \prod_{j=1}^{m-n} \frac{1}{j!} 
\underset {1 \le i \le m-n}{\det_{1 \le j \le m-n}}
\left( (j-i+2)_{i-1} \right).
\end{equation}
The entries of this determinant are 
monic polynomials in $j$ of degree $i-1$, where $i$ and $j$ denote as usual the index of 
the row and the column of the entry. It is now straightforward to reduce this 
determinant by appropriate row operations to Vandermonde's determinant,
\begin{equation}
\label{vand}
\underset {1 \le j \le m-n}{\det_{1 \le i \le m-n}}
\left( (j-i+2)_{i-1} \right)  =
\underset{1 \le j \le m-n}{\det_{1 \le i \le m-n}}
\left( j^{i-1} \right)=\prod_{i=1}^{m-n-1} i!.
\end{equation}
The last equation was obtained by applying the well-known formula for 
Vandermonde's determinant to 
$\det_{1 \le i \le m-n, 1 \le j \le m-n} \left( j^{i-1} \right)$. 

Using \eqref{vand} in \eqref{matrix:3}, we obtain that the determinant in
\eqref{matrix:2} equals
\begin{equation*}
(b+c+n)_{m-n}  \prod_{j=1}^{m-n} \frac{1}{j!} 
\prod_{i=1}^{m-n-1} i! =
 \frac{(b+c+n)_{m-n} }{(m-n)!}.
\end{equation*}
Accordingly, we replace the determinant in \eqref{kA}  by its value $(b+c+n)_{m-n}/(m-n)!$.
Thus, we obtain  the following triple sum for the determinant in \eqref{matrix:allg:1}: 
\begin{multline*}
\left( \prod_{i=2}^{a} \frac{(b+c+i-1)!(i-1)!}{(b+i-1)!(c+i-1)!} \right) \\
 \quad \times
\sum_{n=1}^{a} \sum_{m=1}^{a} \sum_{s=1}^{m}  \left[ (-1)^{n+s} 
\binom{c+x-y+n-2}{x-1} \binom{b-x+y+s-1}{b-x+s-1} \binom{m-1}{s-1}  \right. \\
 \quad \times \left.
\frac{(b+1)_{s-1}}{(b+c+1)_{s-1}}\frac{(c+1)_{n-1}}{(n-1)!} \frac{(b+c+n)_{m-n}}{(m-n)!}  \right].
\end{multline*}
Therefore Lemma~\ref{dreisum} is finally proved  since 
\begin{equation*}
\prod_{i=2}^{a} \frac{(b+c+i-1)!(i-1)!}{(b+i-1)!(c+i-1)!} =
\left( \prod_{i=1}^{a} \prod_{j=1}^b \prod_{k=1}^c \frac{i+j+k-1}{i+j+k-2} \right)
\frac{(1)_c}{(b+1)_c}.
\end{equation*}

\section{Evaluation of the determinants}
\label{comp}

In this section we compute the determinant from Lemma~\ref{det}  
in case that the fixed rhombus is placed in the centre 
(see Lemma~\ref{detlem}) and in case that the lowest vertex of 
the fixed rhombus is placed in the centre (see Lemma~\ref{detlem2}).
The combination of these two lemmas and Lemma~\ref{det} 
then establish Theorem~\ref{th:hp} and Theorem~\ref{th:hp2}.

In order to do so, we first have to figure out which integer point $(x,y)$ 
corresponds to the central rhombus, respectively the `almost central' rhombus above the centre, 
via the bijection described in 
Section~\ref{bijection}. In the oblique angled coordinate system 
introduced in Remark~\ref{cor} the point $((a+b)/2,(a+c)/2)$ is the 
centre of the considered hexagon.
Because of that, $((a+b)/2,(a+c-1)/2)$ is the lowest vertex of the 
central rhombus and therefore the integer point that corresponds to the 
central rhombus (see Remark~\ref{cor}).  Accordingly, $((a+b)/2,(a+c)/2)$ is the lowest vertex of 
the `almost central' rhombus above the centre.

The evaluations of both determinants, namely the determinant  corresponding to the case 
that the  fixed rhombus is placed in the centre, and the determinant corresponding to the 
case that the fixed rhombus is 
placed next to the centre on the other hand, are quite similar. Therefore we concentrate on the first case, 
see the following lemma. 
Lemma~\ref{detlem2} at the end of this section  is devoted to 
the second case, but there I only explain the major differences to the first case.   

\begin{lem} 
\label{detlem}
Let $a$,$b$,$c$ be integers. Then the determinant 

\unitlength1cm
\begin{picture}(0,5)
\put(-1,2.5){\parbox{5cm}{
\begin{equation*} 
-\det_{1 \le i, j \le a+1} \left(
\begin{array}{ccc|c}
      &                    & &   \\
      &\displaystyle \binom{b+c}{c-i+j}  & & \displaystyle \binom{\frac{b+c-1}{2}}{\frac{c}{2}-i+\frac{a+1}{2}} \\
      &                    & &   \\ \hline
      &\displaystyle \binom{\frac{b+c-1}{2}}{\frac{c}{2}-\frac{a+1}{2}+j}& & 0 
\end{array}
\right)
\end{equation*}}}
\put(5.2,3.3){\makebox(2,1)[t]{ \scriptsize $1 \le j \le a$}}
\put(8.2,3.3){\makebox(2,1)[t]{ \scriptsize $j=a+1$}}
\put(11,2.1){\makebox(2,1)[t]{ \scriptsize $1 \le i \le a$}}
\put(11,0.8){\makebox(2,1)[t]{ \scriptsize $i=a+1$}}
\end{picture}

{\parindent0cm is equal to  }

\begin{align*}
&\left( \prod_{i=1}^a \prod_{j=1}^b \prod_{k=1}^c \frac{i+j+k-1}{i+j+k-2}\right)
   \frac{(1)_c}{(b)_c}  \frac{1}{(b+c+1)_{a-1}} 
\binom{\frac{b+c-1}{2}}{\frac{b-1}{2}} \binom{\frac{a+b+c-2}{2}}{\frac{b-1}{2}} 2^{a-1}\\ 
&\times \sum_{k=0}^{(a-1)/2} \left[ \left(\frac{c+1}{2}\right)_k 
      \left(\frac{1+b+c}{2}\right)_{k}  
      \left(\frac{c+2k+2}{2}\right)_{(a-2k-1)/2} \right. \\ 
      & \hspace{2cm} \left. \times  
      \left(\frac{b+c+2k+3}{2}\right)_{(a-2k-1)/2}
      \frac{(\frac{1}{2})_{(a-2k-1)/2}}{(1)_{(a-2k-1)/2}} \right] 
\end{align*}
in case that $a$ is odd, and 
\begin{align*}
&\left( \prod_{i=1}^a \prod_{j=1}^b \prod_{k=1}^c \frac{i+j+k-1}{i+j+k-2}\right)
   \frac{(1)_c}{(b)_c}  \frac{b}{(b+c+1)_{a-1}}      
\binom{\frac{b+c-1}{2}}{\frac{b}{2}} \binom{\frac{a+b+c-1}{2}}{\frac{b}{2}} 2^{a-2} \\
&\times \sum_{k=0}^{(a-2)/2} \left[ \left(\frac{c+2}{2}\right)_k           
      \left(\frac{1+b+c}{2}\right)_{k} 
      \left(\frac{c+2k+3}{2}\right)_{(a-2k-2)/2} \right. \\  
      & \hspace{2cm} \left. \times 
      \left(\frac{b+c+2k+3}{2}\right)_{(a-2k-2)/2}
      \frac{(\frac{1}{2})_{(a-2k-2)/2}}{(1)_{(a-2k-2)/2}} \right]  
\end{align*}
in case that $a$ is even.
\end{lem} 

{\bf Outline of the proof of Lemma~\ref{detlem}.} The following procedure is used for 
computing the determinant: The determinant 
depends on the side lengths $a$, $b$ and $c$ of the hexagon. More precisely: The 
dimension of the matrix underlying the determinant is $a+1$, and the 
entries are certain binomial coefficients that depend on $b$, $c$ and $a$.
First, in Step~1, we reduce the 
problem to the computation of a  determinant with entries that are polynomials in $b$ and $c$ if we fix $a$, see 
\eqref{matrix}, and 
therefore a determinant that is itself a polynomial in those two variables.  
The comparison of the determinant in Lemma~\ref{detlem} and 
the determinant in \eqref{matrix}, the latter  will be denoted by $\det(D_a(b,c))$,  
yields that we have to show that 
the determinant $\det(D_a(b,c))$  is a product of certain linear factors and 
an irreducible  polynomial (over $\mathbb{Z}$) in $b$ and $c$, see \eqref{erg:polodd}, 
respectively \eqref{erg:poleven}. 

In Step 2 of our proof we show that every linear factor on the right-hand side of  \eqref{erg:polodd}, 
respectively \eqref{erg:poleven}, 
is indeed a linear factor of  $\det(D_a(b,c))$. I explain the used
procedure by an example: We have to show, e.g., that the linear factor $(c+1)$ is a factor of 
$\det(D_a(b,c))$ in case that $a$ is even (see 
\eqref{erg:poleven}). It is 
a fundamental algebraic fact that, in order to show that, 
it suffices to show that $\det(D_a(b,-1))=0$ if $a$ is even. (Here we use the fact that 
$\det(D_a(b,c))$ is a polynomial in $c$ if we fix $a$). 
Clearly,  a determinant of 
a matrix is equal to zero if and only if there exists a linear combination 
of rows, or, equivalently, of columns. Therefore we find a linear combination 
of rows of $D_a(b,-1)$ which vanishes, and then prove it, see \eqref{c+i:Lk}. Proving in this case means to 
establish hypergeometric identities
\footnote{At this point it is worth mentioning that I used C. Krattenthaler's 
{\it Mathematica} package HYP \cite{hyp} to handle most of the hypergeometric identities within this paper.}.

At this point it is worth mentioning that the procedure described so far gives 
a complete solution for determinants that factorise completely into linear factors.

In Step 3 of the proof we finally compute the irreducible polynomial,  which will be denoted by $P_a(b,c)$.
In order to do so, we first look for special values of $b$, where $P_a(b,c)$ is `nice'. 
Indeed, we discovered that $P_a(b,c)$ factors completely into linear factors for $b=-c-k$, 
$k=1,3,\dots,2 \lfloor(a-1)/2 \rfloor +1$. 
We work out these evaluations of $P_a(b,c)$ in \eqref{auswert2} and 
\eqref{auswert1}, and subsequently prove them by making use 
of the triple sum derived in Section~\ref{dreisumsec}. As 
it turns out, the degree of $P_a(b,c)$ as a polynomial in $b$ is exactly 
$ \lfloor(a-1)/2 \rfloor$. Thus, the above $ \lfloor(a-1)/2 \rfloor +1 $ evaluations suffice to 
compute $P_a(b,c)$ by using Lagrange interpolation, see \eqref{ir:1} and 
\eqref{ir:2}.

{\bf Proof of Lemma~\ref{detlem} -- The details.}

{\bf Step 1: From our determinant to a determinant with polynomial entries.}

In Section~\ref{bijection} it was finally shown that the number of rhombus tilings which contain  
a fixed horizontal rhombus with lowest vertex $(x,y)$ is given by the following $(a+1) \times (a+1)$ determinant:

\unitlength1cm
\begin{picture}(0,3)
\put(-0.5,1.5){\parbox{3cm}{
\begin{equation} 
\label{m:a:1}
-\det_{1 \le i, j \le a+1} \left(
\begin{array}{ccc|c}
      &                    & &   \\
      &\binom{b+c}{c-i+j}  & & \binom{b-x+y}{y-i+1} \\
      &                    & &   \\ \hline
      &\binom{c+x-y-1}{x-j}& & 0 
\end{array}
\right)
\end{equation}}}
\put(6.2,1.8){\makebox(2,1)[t]{ \scriptsize $1 \le j \le a$}}
\put(8.5,1.8){\makebox(2,1)[t]{ \scriptsize $j=a+1$}}
\put(10.3,0.9){\makebox(2,1)[t]{ \scriptsize $1 \le i \le a$}}
\put(10.3,-0.1){\makebox(2,1)[t]{ \scriptsize $i=a+1$}}
\put(11,0.7){\makebox(2,1)[t]{.}}
\end{picture}

With displeasure  we observe that the matrix underlying this determinant
has an exceptional row and an 
exceptional column, namely the $(a+1)$th  in both cases. In the following I 
describe some elementary row and column operations which lead 
to a matrix having the same determinant but a  `homogeneous' definition:

1. First we want the entries in the $(a+1)$th row to be zero ---  except 
for the first entry. Since $\binom{c+x-y-1}{x-j}$ is the $j$th entry in the 
$(a+1)$th row, we have to subtract 
$\binom{c+x-y-1}{x-j}/\binom{c+x-y-1}{x-j+1}=(x-j+1)/(c-y+j-1)$ times 
the $(j-1)$th column from the $j$th, starting at $j=a$ and stopping at 
$j=2$. We obtain a new matrix with the desired behaviour in the 
$(a+1)$th row. The $(i,j)$ entry of this new matrix is  
\begin{equation*}
\binom{b+c}{c-i+j}-\frac{x-j+1}{c-y+j-1} \binom{b+c}{c-i+j-1},
\end{equation*}
$i,j=2,3, \dots, a$, it is $\binom{c+x-y-1}{x-1}$
for $i=a+1$ and $j=1$, and it is $0$ for $i=a+1$ and $j \not= 1$. The other 
entries do not change.

2. Next we want to do the same for the $(a+1)$th column: Accordingly, since 
$\binom{b-x+y}{y-i+1}$ is the $i$th entry in the $(a+1)$th column, we 
subtract $\binom{b-x+y}{y-i+1}/\binom{b-x+y}{y-i+2}=(y-i+2)/(b+i-x-1)$ times 
the $(j-1)$th row from the $j$th row,  again starting at the bottom --- that is to say, 
$j=a,a-1,\dots, 2$.  Therefore the $(i,j)$ entry of our new matrix is given as 
follows:
\begin{equation}
\label{entry}
\begin{split}
&\binom{b+c}{c-i+j}-\frac{x-j+1}{c-y+j-1} \binom{b+c}{c-i+j-1}\\
-\frac{y-i+2}{b+i-x-1} & \left(  \binom{b+c}{c-(i-1)+j} - 
\frac{x-j+1}{c-y+j-1} \binom{b+c}{c-(i-1)+j-1}  \right)
\end{split}
\end{equation}
if $i,j=2,3 \dots, a$, it is $\binom{b-x+y}{y}$ for $j=a+1$ and $i=1$, and it is $0$ for $j=a+1$ and
$i \not= 1$. Again, the other entries do not change.

3. In summary, we obtain the following determinant,
\begin{equation*}
\det_{ 1 \le i,j \le a+1} \left(
\begin{array}{c c c c c}
*  & * & * & * & \binom{b-x+y}{y}  \\ 
* &   &  &  & 0 \\  
*  & &a_{ij} & & \vdots \\
*  & & & & \vdots \\
\binom{c+x-y-1}{x-1} & 0 & \dots &  \dots & 0
\end{array} 
\right),
\end{equation*}
where $*$ describes an entry which is of no interest for us, and where $a_{ij}$ denotes 
the expression in \eqref{entry}. This form suggests an 
expansion of the determinant along the last row and along the last column. Therefore
\eqref{m:a:1} is equal to
\begin{equation} 
\label{factor1}
-\binom{c+x-y-1}{x-1} \binom{b-x+y}{y} \det_{2 \le i,j \le a} (a_{ij}).
\end{equation}

We have reduced our problem to computing $\det_{2 \le i,j \le a} (a_{ij})$ and therefore 
the job to compute the determinant of a `homogeneous' matrix, even if the entries 
are now more complex. The advantage is that we can now take several factors out of the 
determinant so that the remaining entries are polynomials: 
A simple calculation shows that
\begin{equation*}
a_{ij}= \frac{(b+c)!(b+i-j+2)_{j-2}(c-i+j+2)_{a-j}}{(b+i-1)!(c-j+a+1)!(c-y+j-1)(b+i-x-1)}
\times H,
\end{equation*}
where $H$ is the polynomial
\begin{equation}
\label{H}
\begin{split}
H= (b+i-j+1)(c-i+j+1)(c-y+j-1)(b+i-x-1) \\
       -(c-i+j)(c-i+j+1)(x-j+1)(b+i-x-1) \\
       -(b+i-j)(b+i-j+1)(y-i+2)(c-y+j-1)\\
       +(b+i-j+1)(c-i+j+1)(x-j+1)(y-i+2).  
\end{split}
\end{equation}
We take $(b+c)!/((b+i-1)!(b+i-x-1))$ out of the $i$th row, and 
$1/((c-j+a-1)!(c-y+j-1))$ out of the $j$th column of the matrix $(a_{ij})_{2 \le i,j \le a}$. 
This gives
\begin{equation}
\label{factor2}
\begin{split}
\det_{2 \le i,j \le a} (a_{ij}) = 
\left(\prod_{i=2}^{a} \frac{(b+c)!}{(b+i-1)!(c-i+a+1)!(c-y+i-1)(b+i-x-1)} \right) \\
\times \det_{2 \le i,j \le a} \left( (b+i-j+2)_{j-2} (c-i+j+2)_{a-j} \cdot H \right).
\end{split}
\end{equation}

Again we have reduced our problem to the evaluation of another `homogeneous' determinant, namely 
\begin{equation}
\label{matrix}
\det_{2 \le i,j \le a} \left((b+i-j+2)_{j-2} (c-i+j+2)_{a-j} \cdot H \right).
\end{equation} 
But this 
one has the pleasing property to be a polynomial in $b$ and $c$ if we fix $a$ (since 
the entries of the underlying matrix are). In the following, $D_a(b,c)$ denotes 
the matrix underlying the determinant in \eqref{matrix} evaluated at $x=(a+b)/2$ and $y=(a+c-1)/2$.

When I computed this determinant for small values of $a$, I was led to the following conjecture 
\footnote{In fact, only the linear factors in front of the sums in \eqref{erg:polodd} and 
\eqref{erg:poleven} can immediately be guessed, by computing $\det( D_a(b,c))$ for 
$a=2, 3, \dots, 8$. In order to work out an conjecture concerning the two sums themselves 
I first computed them for $a=2,3, \dots, 23$.
On the basis of these data, I worked out (guesses for) `nice' evaluations 
of these polynomials at special values of $b$, for {\em all\/ } values of $a$. (See the Outline of 
the proof of Lemma~\ref{detlem} and Step~3 of the 
proof of Lemma~\ref{detlem} for more details, and, in particular, for an explicit listing of these 
special values of $b$).   Thereby I was helped by  Krattenthaler's {\it Mathematica} `guessing machine' RATE (which is available via 
Internet at {\tt http://radon.mat.univie.ac.at/People/kratt/rate/rate.html}; see also \cite[Appendix A]{kratdet}). Once I had got thus far, I computed 
the polynomials  
in \eqref{erg:polodd} and \eqref{erg:poleven} by Lagrange interpolation.}:  
There  holds
\begin{multline}
\label{erg:polodd} 
\det_{2 \le i,j \le a} \left(D_a(b,c)\right) 
=
\left(\prod_{i=2}^{a-1} \left(1+b+c\right)_{i-1} \right)  
  \left(\prod_{i=2}^{a} (i-1)!  \right) 2^{a-1} \\
  \left(\frac{b+1}{2}\right)_{(a-1)/2}^2 
  \left(\frac{c+2}{2}\right)_{(a-1)/2} 
  \left(\frac{1+b+c}{2}\right)_{(a-1)/2} \\
  \quad \times \sum_{k=0}^{(a-1)/2} \left[ \left(\frac{c+1}{2}\right)_k  
    \left(\frac{1+b+c}{2}\right)_k 
    \left(\frac{c+2k+2}{2}\right)_{(a-2k-1)/2} \right. \\ 
   \left.\hspace{2cm} \times \left(\frac{b+c+2k+3}{2}\right)_{(a-2k-1)/2}    
    \frac{(\frac{1}{2})_{(a-2k-1)/2}}{(1)_{(a-2k-1)/2}}  \right]   
\end{multline}
if $a$ is odd, and
\begin{multline}
\label{erg:poleven}
\det_{2 \le i,j \le a} \left(D_a(b,c)\right) 
=                                      
\left( \prod_{i=2}^{a-1} \left(1+b+c\right)_{i-1}  \right)
  \left(\prod_{i=2}^{a} (i-1)!  \right) 2^{a-2} \\
   b \left(\frac{b+2}{2}\right)^2_{(a-2)/2} \left(\frac{c+1}{2}\right)_{a/2} 
     \left(\frac{1+b+c}{2}\right)_{a/2}   \\
    \quad \times \sum_{k=0}^{(a-2)/2} \left[ \left(\frac{c+2}{2}\right)_k 
     \left(\frac{1+b+c}{2}\right)_{k} 
     \left(\frac{c+2k+3}{2}\right)_{(a-2k-2)/2} \right. \\
    \left. \hspace{2cm} \times \left(\frac{b+c+2k+3}{2}\right)_{(a-2k-2)/2}
     \frac{(\frac{1}{2})_{(a-2k-2)/2}}{(1)_{(a-2k-2)/2}}   \right]        
\end{multline}     
if $a$ is even.
If we remember the assertion in Lemma~\ref{detlem} and the factors we have taken out 
of the determinant in \eqref{factor1} and 
\eqref{factor2}, and if we then specialise $x=(a+b)/2$ and $y=(a+c-1)/2$, a tedious but 
straightforward calculation  yields that in order to complete a proof of 
Lemma~\ref{detlem} it remains to show \eqref{erg:polodd} and 
\eqref{erg:poleven}. 

As described in the outline of the proof, we first prove that the claimed linear factors, which 
are 
\begin{equation}
\label{lf:odd}
\left( \prod_{i=2}^{a-1} \left(1+b+c\right)_{i-1}  \right)
\left(\frac{b+1}{2}\right)_{(a-1)/2}^2 
\left(\frac{c+2}{2}\right)_{(a-1)/2} 
\left(\frac{1+b+c}{2}\right)_{(a-1)/2} 2^{2a-2}
\end{equation}
if $a$ is odd, and, respectively, 
\begin{equation} 
\label{lf:even}
\left( \prod_{i=2}^{a-1} \left(1+b+c\right)_{i-1} \right)
b \left(\frac{b+2}{2}\right)^2_{(a-2)/2} \left(\frac{c+1}{2}\right)_{a/2} 
     \left(\frac{1+b+c}{2}\right)_{a/2} 2^{2a-2}
      \end{equation}     
if $a$ is even, divide $\det(D_a(b,c))$ (see Step 2). In Step 3 we finally
calculate the remaining irreducible  polynomial in $b$ and $c$ (over $\mathbb{Z}$), which reads 
\begin{multline}  
\label{ir:1}
\left(\prod_{i=2}^{a} (i-1)! \right) \left(\frac{1}{2}\right)^{a-1} \sum_{k=0}^{(a-1)/2} \left[ \left(\frac{c+1}{2}\right)_k  
    \left(\frac{1+b+c}{2}\right)_k 
    \left(\frac{c+2k+2}{2}\right)_{(a-2k-1)/2} \right. \\ 
   \left.\hspace{2cm} \times \left(\frac{b+c+2k+3}{2}\right)_{(a-2k-1)/2}    
    \frac{(\frac{1}{2})_{(a-2k-1)/2}}{(1)_{(a-2k-1)/2}}  \right] 
\end{multline}
if $a$ is odd, and respectively,
\begin{multline}
\label{ir:2}
\left(\prod_{i=2}^{a} (i-1)! \right) \left( \frac{1}{2} \right)^{a}
\sum_{k=0}^{(a-2)/2} \left[ \left(\frac{c+2}{2}\right)_k 
     \left(\frac{1+b+c}{2}\right)_{k} 
     \left(\frac{c+2k+3}{2}\right)_{(a-2k-2)/2} \right. \\
    \left. \hspace{2cm} \times \left(\frac{b+c+2k+3}{2}\right)_{(a-2k-2)/2}
     \frac{(\frac{1}{2})_{(a-2k-2)/2}}{(1)_{(a-2k-2)/2}}   \right] 
\end{multline}
if $a$ is even.

\medskip

{\bf Step 2: The linear factors in \eqref{lf:odd}, respectively in \eqref{lf:even},  divide 
$\det(D_a(b,c))$ as a polynomial in $b$ and $c$.} 

Essentially we have four different types of linear factors:

\begin{enumerate}
\item \label{c+k} Factors of the form $(c+k)$: If $a$ is odd, $k=2,4,\dots, a-1$. Otherwise  
$k=1,3,\dots, a-1$. 

\item \label{b+k} Factors of the form $(b+k)^2$: If $a$ is odd, $k=1,3,\dots, a-2$. Otherwise 
$k=2,4,\dots, a-2$. The factor $b$ occurs once if $a$ even. 

\item \label{b+c+i,i-ung}Factors of the form $(b+c+k)^{a-k}$ for 
$k=1,3,\dots, 2\lceil(a-1)/2\rceil-1$.

\item \label{b+c+i,i-ger}Factors of the form $(b+c+k)^{a-k-1}$ for 
$k=2,4,\dots, 2\lfloor(a-1)/2\rfloor$.
\end{enumerate}

{\it re 1. --- The factors of the form $(c+k)$ divide $\det(D_a(b,c))$:} 
In order  to show that $(c+k)$ is a linear factor of 
$\det(D_a(b,c))$, where $1 \le k \le a$ and 
$k  \not\equiv a \hspace{2mm} (\text{mod} \hspace{1mm} 2)$, it suffices to show that 
$\det(D_a(b,-k))=0$ for those special $k$'s.  Therefore we search for 
linearly dependent rows or columns in the matrices $D_a(b,-k)$.

By computer experiments
\footnote{Again, I computed these linear combinations for small values of $a$, i.e., I solved the following 
system of linear equations in $c_i(a,b,k)$: 
$ \sum_{i=2}^{a} c_i(a,b,k) D_a(b,-k)_{(i,j)} = 0$, $j=2,3, \dots,a$, for 
$a=2,3, \dots, 15$. I was led to an conjecture for the coefficient $c_i(a,b,k)$ for all values of 
$a$ by using Krattenthaler's {\it Mathematica} `guessing machine' RATE 
(which is available via 
Internet at {\tt http://radon.mat.univie.ac.at/People/kratt/rate/rate.html}).}
, I found  the                                                
following non-trivial linear combination of rows:
\begin{equation}
\label{c+i:Lk}
\sum_{i=(a-k+3)/2}^{a-k+1} 
\frac{(-1)^{i} (b+i)_{a-k+1-i} \left(\frac{-a+k-1+2i}{2}\right)_{a-k+1-i}}
{(1)_{a-k+1-i} \left(\frac{b-a+2i-2}{2}\right)_{a-k+1-i}} D_a(b,-k)_{(i,j)}=0,
\end{equation}
for $j=2,3, \dots, a$.

In order to complete the proof that $(c+k)$ divides the determinant $\det(D_a(b,c))$, 
we just have to prove this identity. {\it Gosper's algorithm} \cite{gosper} 
for hypergeometric sums
(which are sums, where the quotient of two successive summands is a rational 
function in the summation index) helps us to recognize that this 
is actually a telescoping sum: First we examine the factor
$((-a+k-1+2i)/2)_{a-k+1-i}$ in the linear combination. I claim 
that this factor is equal to zero if $i<(a-k+3)/2$ and because of that we can  
omit the lower bound on the summation index  on the left-hand side of \eqref{c+i:Lk}: 
In this case $(-a+k-1+2i)/2 \le 0$ and 
$(a-k-1)/2 \ge 0$. Because of that, and since  $(-a+k-1+2i)/2$ is an integer, one of the factors 
of the Pochhammer symbol
\begin{equation*}
\left(\frac{-a+k-1+2i}{2} \right)_{a-k+1-i}= 
\left(\frac{-a+k-1+2i}{2} \right)
\left(\frac{-a+k-1+2i}{2}+1 \right) \dots
\left(\frac{a-k-1}{2} \right)
\end{equation*} 
is equal to zero.
If we then reverse the order of summation, we 
obtain that we must show the following identity:
\begin{equation}
\label{c+i:id}
\sum_{i=0}^{\infty}
\frac{(-1)^j \left(\frac{1}{2} - \frac{a}{2}  +  \frac{k}{2} \right)_i
              (-a-b+k)_{i+j-2} }
     {(1)_{-a+i+j} \left( 1 - \frac{a}{2} - \frac{b}{2} + k  \right)_{i}} 
     H \left|_{x=(a+b)/2,y=(a+c-1)/2,c=-k,i \rightarrow a-k+1-i} \right. =0.
\end{equation}     
Let $f(i)$ denote the summand of the sum in \eqref{c+i:id}. 
The magnificent algorithm due to Gosper
decides whether there exists a hypergeometric $g(i)$ with
$$
g(i+1)-g(i)=f(i).
$$
And in case of its existence the algorithm also computes $g(i)$. 
If such an $g(i)$ was found we would have 
$$
\sum_{i=0}^{\infty} f(i) = \sum_{i=0}^{\infty} \left( g(i+1)-g(i) \right) = 
\lim_{i \to \infty} \left( g(i) - g(0) \right).
$$
($\lim_{i \to \infty} g(i)$ exists in our cases, since $f(i)=0$ for all but a finite number 
of $i$ and therefore $g(i)$ is finally constant.)   

A computer implementation of Gosper's  algorithm \cite{paule} prints out 
\begin{equation*}
\frac{(-1)^{j} (-a-b+k)_{i+j-2} \left( \frac{1}{2} - \frac{a}{2} + \frac{k}{2} \right)_i (-1-b+k)
(i-j+k)}{2 (1)_{-1-a+i+j} \left( 1 - \frac{a}{2} - \frac{b}{2} + k \right)_{i-1}}
\end{equation*}
to be a suitable $g(i)$ for our $f(i)$. One may check the identity 
$g(i+1)-g(i)=f(i)$ by dividing the left-hand side by the right-hand side
and simplifying the resulting rational function to  1.

This implies $\sum\limits_{i=0}^{\infty} f(i)=0$, i.e., the truth of \eqref{c+i:id}, since 
$g(0)=0$, caused by the factor $1/(1)_{-1-a+j}$, and $g(i)=0$ for 
$i \ge (1+a-k)/2$, caused by $((1+a-k-2i)/2)_i$ (again we use the 
fact that $1 \le k \le a$ and $k \not\equiv a \hspace{2mm} (\text{mod} \hspace{1mm} 2)$).
Thus we have proved that $c+k$ is a factor of $\det (D_a(b,c))$.

\smallskip

{\it re 2. --- The factors of the form $(b+k)$ divide $\det(D_a(b,c))$:}
Next we have to prove that $(b+k)^2$ is a factor of 
$\det(D_a(b,c))$ if $0<k<a$ and $k \equiv a \hspace{2mm} (\text{mod} \hspace{1mm} 2)$, and, 
in addition, that $b$  is a factor of $\det(D_a(b,c))$  if $a$ is even.

In case that $k=a-2$, this is easy: If we examine  the matrix 
$D_a(-a+2,c)$ carefully, we observe that the
$(a-1)$th and $a$th row vanish. (For $j>2$ this is caused by 
$(b+i-j+2)_{j-2}$, and for $j=2$ by the nasty polynomial 
$H$.) Therefore $(b+a-2)$ is a factor of the $(a-1)$th row and a factor of the 
$a$th row of $D_a(b,c)$. Thus, we have proved that $(b+a-2)^2$ is a factor of the 
determinant $\det(D_a(b,c))$. 

In case that $k < a-2$ we use the same method as for the factors of form $(c+k)$ above.
But since we are now dealing with factors of higher multiplicity, it is not 
enough to find just one linear combination of 
rows or columns of $D_a(-k,c)$. In fact we have to find two linearly independent 
combinations if we want to prove that $(b+k)^2$ is a factor. 

We claim that  
\begin{multline} 
\label{b+k:Lk1}
\sum_{i=k+2}^{(a+k+2)/2}
 \frac{(-1)^{i-k} (c+a-i+2)_{i-k-2} \left(\frac{a+k-2i+4}{2}\right)_{i-k-2}}
{(1)_{i-k-1} \left(\frac{c+a-2i+3}{2}\right)_{i-k-2}^2} \\ 
\times p_{i-k-1}(c+a-k-i+1) D_a(-k,c)_{(i,j)} =0
\end{multline}
if $0 \le k < a-2$ and $k \equiv a \hspace{2mm} (\text{mod} \hspace{1mm} 2)$, and  that
\begin{multline}
\label{b+k:Lk2}
\sum_{i=k+3}^{(a+k+2)/2} 
\frac{(-1)^{i-k} (c+a-i+2)_{i-k-1} \left(\frac{a+k-2i+4}{2}\right)_{i-k-1}}
{(1)_{i-k-1} \left(\frac{c+a-2i+3}{2}\right)_{i-k-2}^2}  \\ 
\times p_{i-k-2}(c+a-k-i) D_a(-k,c)_{(i,j)}  =- D_a(-k,c)_{(k+1,j)}
\end{multline}
if $0 < k < a-2$ and $k \equiv a \hspace{2mm} (\text{mod} \hspace{1mm} 2)$ 
are such linearly independent linear combinations of rows, 
$p_n(c)$ being the sequence of polynomials given by
\begin{equation}
\label{seq}
p_n(c)=\sum_{h=0}^{n-1} \left(\frac{1+c-n}{2} \right)_{n-h-1} 
\left(\frac{1+c-2h+n}{2}\right)_{h}.
\end{equation}
Notice that the exceptional factor $b$ of $\det(D_a(b,c))$ for $a$ is even, is included in the 
first linear combination \eqref{b+k:Lk1}. 

In order to see that the two linear 
combinations in \eqref{b+k:Lk1} and \eqref{b+k:Lk2} are linearly independent we remark that 
the first linear combination \eqref{b+k:Lk1} involves the $(k+1)$th row of 
$D_a(-k,c)$ whereas the second linear combination \eqref{b+k:Lk2} does not (here we use 
$k<a-2$).

\smallskip

It remains to show the two hypergeometric identities \eqref{b+k:Lk1} and  \eqref{b+k:Lk2}. 
The situation is a bit more complicated this 
time compared to the hypergeometric identity \eqref{c+i:Lk} that proved  
that $(c+k)$ divides $\det(D_a(b,c))$: Actually these 
two identities are double sum identities, since they involve the 
polynomials $p_n(c)$. Luckily the proofs of the two identities are 
quite similar. We start with \eqref{b+k:Lk1}.

To begin with, we  split the sum into 
four smaller sums --- each one corresponding to one of the four 
summands of $H$ (according to  the representation of $H$ in \eqref{H}): 
\begin{multline} 
\label{b+k:1:1}
\sum_{i=k+2}^{(a+k+2)/2} \frac{(-1)^{i-k-2} (c+a-i+2)_{i-k-2} \left(\frac{a+k-2i+4}{2}\right)_{i-k-2}}
{(1)_{i-k-1} \left(\frac{c+a-2i+3}{2}\right)_{i-k-2}^2} \\ 
\times p_{i-k-1}(c+a-k-i+1) D_a(-k,c)_{(i,j)}=\\
\sum_{i=k+2}^{(a+k+2)/2} \frac{(-1)^{i-k-2} (c+a-i+2)_{i-k-2} \left(\frac{a+k-2i+4}{2}\right)_{i-k-2}}
{(1)_{i-k-1} \left(\frac{c+a-2i+3}{2}\right)_{i-k-2}^2} 
p_{i-k-1}(c+a-k-i+1)  \\ \times (-k+i-j+1)_{j-1} (c-i+j+1)_{a-j+1}
((-1-a+c+2j)/2)
((-2-a+2i-k)/2)\\
 +\sum_{i=k+2}^{(a+k+2)/2} \frac{(-1)^{i-k-2} (c+a-i+2)_{i-k-2} \left(\frac{a+k-2i+4}{2}\right)_{i-k-2}}
{(1)_{i-k-1} \left(\frac{c+a-2i+3}{2}\right)_{i-k-2}^2} 
p_{i-k-1}(c+a-k-i+1) \\ \times (-k+i-j+2)_{j-2} (c-i+j)_{a-j+2}
((-2-a+2i-k)/2)
((2+a-2j-k)/2)\\
+\sum_{i=k+2}^{(a+k+2)/2} \frac{(-1)^{i-k-2} (c+a-i+2)_{i-k-2} \left(\frac{a+k-2i+4}{2}\right)_{i-k-2}}
{(1)_{i-k-1} \left(\frac{c+a-2i+3}{2}\right)_{i-k-2}^2} 
p_{i-k-1}(c+a-k-i+1) \\ \times (-k+i-j)_{j} (c-i+j+2)_{a-j} 
((3+a+c-2i)/2)
((-1-a+c+2j)/2)\\
+\sum_{i=k+2}^{(a+k+2)/2} \frac{(-1)^{i-k-2} (c+a-i+2)_{i-k-2} \left(\frac{a+k-2i+4}{2}\right)_{i-k-2}}
{(1)_{i-k-1} \left(\frac{c+a-2i+3}{2}\right)_{i-k-2}^2} 
p_{i-k-1}(c+a-k-i+1) \\ \times (-k+i-j+1)_{j-1} (c-i+j+1)_{a-j+1} 
((3+a+c-2i)/2)((2+a-2j-k)/2).
\end{multline}

Next we examine the sum of the first and the third sum of the right-hand side in  
\eqref{b+k:1:1}. After some simplifications and interchange of summations one gets 
\begin{multline*}
\sum_{h=0}^{(a-k-2)/2} \left[ \sum_{i=h+k+2}^{(a+k+2)/2} 
(-1)^{a-j+1} \frac{\left(-\frac{1}{2}-\frac{a}{2}+\frac{c}{2}+j \right)
                   \left(\frac{1}{2}-\frac{a}{2}-\frac{c}{2}+k \right)_{h}}
                  {\left(\frac{3}{2}-\frac{a}{2}-\frac{c}{2}+k\right)_{h}}   \right. \\
\left. \times \frac{\left(1-\frac{a}{2}+\frac{k}{2} \right)_{-1+i-k}
             (1-a-c+k)_{-1+a+i-j-k}}
             {(1)_{i-j-k} \left(\frac{3}{2}-\frac{a}{2}-\frac{c}{2}+k\right)_{-2+i-k}} \right] \\
- \sum_{h=0}^{(a-k-2)/2} \left[ \sum_{i=h+k+2}^{(a+k+2)/2} 
(-1)^{a-j+1} \frac{\left(-\frac{1}{2}-\frac{a}{2}+\frac{c}{2}+j \right)
                   \left(\frac{1}{2}-\frac{a}{2}-\frac{c}{2}+k \right)_{h}}
                  {\left(\frac{3}{2}-\frac{a}{2}-\frac{c}{2}+k\right)_{h}}   \right.\\
\left. \times \frac{\left(1-\frac{a}{2}+\frac{k}{2} \right)_{-2+i-k}
             (1-a-c+k)_{-2+a+i-j-k}}
             {(1)_{-1+i-j-k} \left(\frac{3}{2}-\frac{a}{2}-\frac{c}{2}+k\right)_{-3+i-k}} \right].
\end{multline*}

Notice that these two double sums  
are quite similar. If we combine the exterior sums, take some common factors out of 
the two inner sum and shift the index of the first inner sum by one, we obtain 
\begin{multline*}               
\sum_{h=0}^{(a-k-2)/2}  
(-1)^{a-j+1} \frac{\left(-\frac{1}{2}-\frac{a}{2}+\frac{c}{2}+j \right)
                   \left(\frac{1}{2}-\frac{a}{2}-\frac{c}{2}+k \right)_{h}}
                  {\left(\frac{3}{2}-\frac{a}{2}-\frac{c}{2}+k\right)_{h}} \\
\left[ \sum_{i=h+k+3}^{(a+k+4)/2}   
\frac{\left(1-\frac{a}{2}+\frac{k}{2} \right)_{-2+i-k}
             (1-a-c+k)_{-2+a+i-j-k}}
             {(1)_{-1+i-j-k} \left(\frac{3}{2}-\frac{a}{2}-\frac{c}{2}+k\right)_{-3+i-k}} \right. \\                                                                         
\left.- \sum_{i=h+k+2}^{(a+k+2)/2}
 \frac{\left(1-\frac{a}{2}+\frac{k}{2} \right)_{-2+i-k}
             (1-a-c+k)_{-2+a+i-j-k}}
             {(1)_{-1+i-j-k} \left(\frac{3}{2}-\frac{a}{2}-\frac{c}{2}+k\right)_{-3+i-k}} \right].
\end{multline*}
This makes us discover, that the two inner sums nearly cancel out each other, because the two summands are exactly the same. 
Therefore the latter expression simplifies further to
\begin{multline*}
\sum_{h=0}^{(a-k-2)/2}  
(-1)^{a-j+1} \frac{\left(-\frac{1}{2}-\frac{a}{2}+\frac{c}{2}+j \right)
                   \left(\frac{1}{2}-\frac{a}{2}-\frac{c}{2}+k \right)_{h}}
                  {\left(\frac{3}{2}-\frac{a}{2}-\frac{c}{2}+k\right)_{h}} \\
\times  
\frac{\left( 1 - \frac{a}{2}+\frac{k}{2} \right)_{h} (1-a-c+k)_{a+h-j}}
{(1)_{1+h-j} \left(\frac{3}{2}-\frac{a}{2}-\frac{c}{2}+k\right)_{-1+h}}.   
\end{multline*}
Thus we have already reduced the sum of the two  double sums to the following 
single sum
\begin{multline}  
\label{simple1}
\sum_{h=0}^{(a-k-2)/2}  
(-1)^{a-j} \frac{\left(-\frac{1}{2}-\frac{a}{2}+\frac{c}{2}+j \right)
                   \left(\frac{1}{2}-\frac{a}{2}-\frac{c}{2}+k \right)
                   \left( 1 - \frac{a}{2}+\frac{k}{2} \right)_{h} (1-a-c+k)_{a+h-j}}
                  {\left(\frac{3}{2}-\frac{a}{2}-\frac{c}{2}+k\right)_{h}(1)_{1+h-j}}.     
\end{multline} 

Analogously one can show that the sum of the second und the fourth sum in \eqref{b+k:1:1}, 
which reads  
\begin{multline*}
\sum_{h=0}^{(a-k-2)/2} \sum_{i=2+h+k}^{(a+k+2)/2} 
(-1)^{-1+a-j} 
\frac{\left(1+\frac{a}{2}-j-\frac{k}{2}\right)
\left(\frac{1}{2}-\frac{a}{2}-\frac{c}{2}+k\right)_h}
{\left(\frac{3}{2}-\frac{a}{2}-\frac{c}{2}+k\right)_{h}} \\
\times 
\frac{\left(1-\frac{a}{2}+\frac{k}{2}\right)_{-1+i-k} 
      \left(1-a-c+k\right)_{a+i-j-k}}
{(1)_{1+i-j-k} \left(\frac{3}{2}-\frac{a}{2}-\frac{c}{2}+k \right)_{-2+i-k}} \\
-  
\sum_{h=0}^{(a-k-2)/2} \sum_{i=2+h+k}^{(a+k+2)/2} 
(-1)^{-1+a-j} 
\frac{\left(1+\frac{a}{2}-j-\frac{k}{2}\right)
\left(\frac{1}{2}-\frac{a}{2}-\frac{c}{2}+k\right)_h}
{\left(\frac{3}{2}-\frac{a}{2}-\frac{c}{2}+k\right)_{h}} \\
\times 
\frac{\left(1-\frac{a}{2}+\frac{k}{2}\right)_{-2+i-k} 
      \left(1-a-c+k\right)_{-1+a+i-j-k}}
{(1)_{i-j-k} \left(\frac{3}{2}-\frac{a}{2}-\frac{c}{2}+k \right)_{-3+i-k}} ,   
\end{multline*}
after interchange of summation and some cancellations, simplifies to 
\begin{equation}
\label{simple2}
\sum_{h=0}^{(a-k-2)/2}   (-1)^{a-j}
\frac{\left(1+\frac{a}{2}-j-\frac{k}{2}\right)
      \left(\frac{1}{2}-\frac{a}{2}-\frac{c}{2} + k \right)
      \left(1-\frac{a}{2}+\frac{k}{2}\right)_{h}
      \left(1-a-c+k\right)_{1+a+h-j}}
      {(1)_{2+h-j} \left(\frac{3}{2}-\frac{a}{2}-\frac{c}{2}+k \right)_{h}}.
\end{equation}

Now it remains to show that the two (single) sums \eqref{simple1} and 
\eqref{simple2} sum up to zero. Because of the factor 
$1/(1)_{1+h-j}$, we are able to change the lower bound 
on the summation index of the sum in \eqref{simple1} to $j-1$, since  $j \ge 2$ and
$1/(1)_n=0$ if $n$ is a negative integer.  Analogously 
we change the lower bound of the summation index of the sum in \eqref{simple2} 
to $j-2$ according to the factor $1/(1)_{2+h-j}$. Therefore the sum  
of \eqref{simple1} and  \eqref{simple2} is
\begin{multline*}   \left(\frac{1}{2}-\frac{a}{2}-\frac{c}{2}+k \right)   \\
 \times \left[ \sum_{h=j-1}^{(a-k-2)/2}  
(-1)^{a-j} \frac{\left(-\frac{1}{2}-\frac{a}{2}+\frac{c}{2}+j \right)
                 \left( 1 - \frac{a}{2}+\frac{k}{2} \right)_{h} (1-a-c+k)_{a+h-j}}
                  {\left(\frac{3}{2}-\frac{a}{2}-\frac{c}{2}+k\right)_{h}(1)_{1+h-j}}
                  \right. \\
+ \left.  \sum_{h=j-2}^{(a-k-2)/2} (-1)^{a-j}
\frac{\left(1+\frac{a}{2}-j-\frac{k}{2}\right)
      \left(1-\frac{a}{2}+\frac{k}{2}\right)_{h}
      \left(1-a-c+k\right)_{1+a+h-j}}
      {(1)_{2+h-j} \left(\frac{3}{2}-\frac{a}{2}-\frac{c}{2}+k \right)_{h}} \right].
\end{multline*}   

Using the standard hypergeometric notation  
\begin{equation*}
{} _{r} F _{s} \!\left[
                       \begin{matrix} 
                          a_1,\dots, a_r  \\
                          b_1,\dots, b_s
                       \end{matrix}; {\displaystyle z} \right]=
\sum_{k=0}^{\infty} \frac{(a_1)_{k} \dots (a_r)_{k}}{k! (b_1)_k \dots (b_s)_{k}} z^k,                          
\end{equation*} 
and, after canceling the factor
$\left(\frac{1}{2}-\frac{a}{2}-\frac{c}{2}+k \right)$, it remains to 
show that
\begin{multline}
\label{b+k:1:2}
{{({ \textstyle 1 - {a\over 2} + {k\over 2}}) _{-1 + j} \,
     ({ \textstyle 1 - a - c + k}) _{-1 + a} }\over 
   {({ \textstyle {3\over 2} - {a\over 2} - {c\over 2} + k}) _{-1 + j}
     }}  \\
\times \left(
{} _{2} F _{1} \!\left [ \begin{matrix} { -{{a}\over 2} + j + {k\over
   2}, -c + k}\\ { {1\over 2} - {a\over 2} - {c\over 2} + j +
   k}\end{matrix} ; {\displaystyle 1}\right ] \,
  ({ \textstyle -{1\over 2} - {a\over 2} + {c\over 2} + j})  \right. \\
\left. - {} _{2} F _{1} \!\left [ \begin{matrix} { -1 - {a\over 2} + j
     + {k\over 2}, -c + k}\\ { -{1\over 2} - {a\over 2} - {c\over 2} +
     j + k}\end{matrix} ; {\displaystyle 1}\right ] \,
    ({ \textstyle -{1\over 2} - {a\over 2} - {c\over 2} + j + k})  
        \right) = 0.
\end{multline} 
The factor $((2-a+k)/2)_{-1+j}$ tells us that this is true at least for 
$j \ge (-2+a-k)/2$.
Otherwise, we use Vandermonde's summation formula (see \cite[(1.7.7), Appendix (III.4)]{slater})  
\begin{equation} 
\label{Vandermonde}
{} _{2} F _{1} \!\left[
                       \begin{matrix} 
                          a,-n  \\
                          c 
                      \end{matrix}; {\displaystyle 1} \right]=   
\frac{(c-a)_{n}}{(c)_{n}},
\end{equation}
which is valid if $n$ is a nonnegative integer.
Namely, if we apply this formula to the left-hand side of  \eqref{b+k:1:2}
with $n=(a-2j-k)/2$, respectively $n=(2+a-2j-k)/2$  (here we use 
$j < (-2+a-k)/2$), we see that the two $_2 F _1$--series in 
\eqref{b+k:1:2}  sum up to zero.
Thus, we have proved the first linear combination,  \eqref{b+k:Lk1}, for the factors of 
form $(b+k)$. 

\smallskip

The proof of the second linear combination \eqref{b+k:Lk2} is quite similar:
Again we split the double sum into four smaller sums according to the summands of the 
polynomial $H$. After some simplifications and interchange of summation we have reduced our 
problem to the following hypergeometric identity:
\begin{multline}
\label{b+k:2}
\sum_{h = 0}^{(-4+a-k)/2}
   \sum_{i = 3 + h + k}
     ^{(2+a+k)/2}
     {{{{\left( -1 \right) }^{1 + a + h - j}}\,
         ({ \textstyle -{1\over 2} - {a\over 2} +
          {c\over 2} + j}) _{1} \,
         ({ \textstyle {{-a}\over 2} + {k\over
          2}}) _{i - k} \,
         ({ \textstyle -a - c + k}) _{a + i - j -
          k} }\over 
       {({ \textstyle 1}) _{i - j - k} \,
         ({ \textstyle -{3\over 2} + {a\over 2} +
          {c\over 2} - h - k}) _{1 + h} \,
         ({ \textstyle {3\over 2} - {a\over 2} -
          {c\over 2} + h + k}) _{-2 - h + i - k} }}    \\
+
\sum_{h = 0}^{(-4+a-k)/2}
   \sum_{i = 3 + h + k}
     ^{(2+a+k)/2}
     {{{{\left( -1 \right) }^{1 + a + h - j}}\,
         ({ \textstyle 1 + {a\over 2} - j -
          {k\over 2}}) _{1} \,
         ({ \textstyle {{-a}\over 2} + {k\over
          2}}) _{i - k} \,
         ({ \textstyle -a - c + k}) _{1 + a + i -
          j - k} }\over 
       {({ \textstyle 1}) _{1 + i - j - k} \,
         ({ \textstyle -{3\over 2} + {a\over 2} +
          {c\over 2} - h - k}) _{1 + h} \,
         ({ \textstyle {3\over 2} - {a\over 2} -
          {c\over 2} + h + k}) _{-2 - h + i - k} 
         }} \\
+
\sum_{h = 0}^{(-4+a-k)/2}
   \sum_{i = 3 + h + k}^{(2+a+k)/2}
     {{{{\left( -1 \right) }^{a + h - j}}\,
         ({ \textstyle -{1\over 2} - {a\over 2} + {c\over 2} + j})
          _{1} \,({ \textstyle {{-a}\over 2} + {k\over 2}}) _{-1 + i -
          k} \,({ \textstyle -a - c + k}) _{-1 + a + i - j - k} }\over
         {({ \textstyle 1}) _{-1 + i - j - k} \,
         ({ \textstyle -{3\over 2} + {a\over 2} + {c\over 2} - h - k})
          _{1 + h} \,({ \textstyle {3\over 2} - {a\over 2} - {c\over
          2} + h + k}) _{-3 - h + i - k} }}  \\
+
\sum_{h = 0}^{(-4+a-k)/2}
   \sum_{i = 3 + h + k}^{(2+a+k)/2}
     {{{{\left( -1 \right) }^{a + h - j}}\,
         ({ \textstyle 1 + {a\over 2} - j - {k\over 2}}) _{1} \,
         ({ \textstyle {{-a}\over 2} + {k\over 2}}) _{-1 + i - k} \,
         ({ \textstyle -a - c + k}) _{a + i - j - k} }\over 
       {({ \textstyle 1}) _{i - j - k} \,
         ({ \textstyle -{3\over 2} + {a\over 2} + {c\over 2} - h - k})
          _{1 + h} \,({ \textstyle {3\over 2} - {a\over 2} - {c\over
          2} + h + k}) _{-3 - h + i - k} }} \\
=-D_a(-k,c)_{(k+1,j)}.                                  
\end{multline}  

As before, the sum of the first and the third summand on the left-hand side 
of \eqref{b+k:2} simplifies to a single sum:
\begin{equation}
\label{simple3}
\sum_{h=0}^{(-4+a-k)/2}  (-1)^{a+h-j} 
\frac{\left(-\frac{1}{2}-\frac{a}{2}+\frac{c}{2}+j\right)_1 
     \left(-\frac{a}{2}+\frac{k}{2}\right)_{2+h} (-a-c+k)_{2+a+h-j}}
     {\left(-\frac{3}{2}+\frac{a}{2}+\frac{c}{2}-h-k\right)_{1+h}(1)_{2+h-j}}.     
\end{equation}
Analogously, the second and fourth summand of \eqref{b+k:2} give 
\begin{equation}
\label{ausnahme}
\sum_{h=0}^{(-4+a-k)/2} (-1)^{a+h-j}
           \frac{\left(1+\frac{a}{2}-j - \frac{k}{2} \right)_1 
           \left(-\frac{a}{2}+\frac{k}{2}\right)_{2+h} 
           (-a-c+k)_{3+a+h-j}}
           {\left(-\frac{3}{2}+\frac{a}{2}+\frac{c}{2}-h-k\right)_{1+h}
            (1)_{3+h-j}}.
\end{equation}

Again we want to change the lower bound of the remaining single  
sums  in \eqref{simple3} and \eqref{ausnahme} to $j-2$, respectively to $j-3$, 
according to the factor 
$1/(1)_{2+h-j}$, respectively $1/(1)_{3+h-j}$. This is because then the 
two sums are $_{2} F _{1}$--series and this 
makes it possible to apply Vandermonde's summation formula \eqref{Vandermonde} again. In case of 
\eqref{ausnahme}  and $j=2$ this change is problematic since $j-3 < 0$ for $j=2$. 
But the right-hand side of 
\eqref{b+k:Lk2} compensates this missing term.

In terms of  hypergeometric notation, it remains to show that
\begin{multline*}
{{{{\left( -1 \right) }^a}\,
     ({ \textstyle {{-a}\over 2} + {k\over 2}}) _{j} \,
     ({ \textstyle -a - c + k}) _{a} }\over 
   {({ \textstyle {1\over 2} + {a\over 2} + {c\over 2} - j - k}) 
      \,({ \textstyle {3\over 2} + {a\over 2} + {c\over 2} - j - k})
      _{-2 + j} }}      \\
\times \left[
{} _{2} F _{1} \!\left [ \begin{matrix} { -{{a}\over 2} + j + {k\over
    2}, -c + k}\\ { {1\over 2} - {a\over 2} - {c\over 2} + j +
    k}\end{matrix} ; {\displaystyle 1}\right ] \,
   ({ \textstyle -{1\over 2} - {a\over 2} + {c\over 2} + j})  \right. \\ +  \left.
  {} _{2} F _{1} \!\left [ \begin{matrix} { -1 - {a\over 2} + j +
    {k\over 2}, -c + k}\\ { -{1\over 2} - {a\over 2} - {c\over 2} + j
    + k}\end{matrix} ; {\displaystyle 1}\right ] \,
   ({ \textstyle {1\over 2} + {a\over 2} + {c\over 2} - j - k})  
\right] =0.  
\end{multline*}  
For $j \ge (a-k+1)/2$ this is obvious because of  the factor 
$((-a+k)/2)_j$. 
Otherwise, we use again Vandermonde's summation formula \eqref{Vandermonde} with 
$n=(a-2j-k)/2$, respectively $n=(2+a-2j-k)/2$.

Finally we have completely proved that $(b+k)^2$ is a factor of 
$\det(D_a(b,c))$.

\smallskip

{\it re 3., 4. ---  $\prod\limits_{i=2}^{a} (b+c+1)_{i-1}/ ((b+c+2)/2)_{\lceil (a-2)/2 \rceil}$ 
is a factor of $\det(D_a(b,c))$:}
In order to prove that the third and the fourth type of factors (see 
the beginning  of Step 2) divide $\det(D_a(b,c))$, we could
use the same procedure as for the factors of form $(c+k)$ and 
the factors of form $(b+k)^2$. But, fortunately, we will see that the triple 
sum that is  the subject of Section~\ref{dreisumsec} provides an easier proof.
`Fortunately,' since in this case the procedure leads to 
quite complicated hypergeometric identities. Furthermore, the 
underlying linear combinations are more difficult to guess because 
an additional parameter originating in the higher multiplicity of the factors is involved. 

If we look at \eqref{factor1} and \eqref{factor2}, we 
see that \eqref{matrix} is equal
to \eqref{drei} evaluated at $x=(a+b)/2$ and $y=(a+c-1)/2$ modulo some factors. Therefore 
a combination of \eqref{factor1}, \eqref{factor2} and  
Lemma~\ref{dreisum} gives 
\begin{multline}
\label{dreisum:1}
\det ( D_a(b,c) )
=\left( \prod_{i=2}^{a} (b+c+1)_{i-1} (i-1)! \right) \left( \frac{c-a+3}{2} \right)_{a-1} 
\left(\frac{b-a+2}{2} \right)_{a-1} \\
 \quad \times
\sum_{n=1}^{a} \sum_{m=1}^{a} \sum_{s=1}^{m}  \left[ (-1)^{n+s} 
\frac{\left(\frac{1-a+c+2n}{2}\right)_{(a+b-2)/2}}
     {\left(\frac{3-a+c}{2}\right)_{(a+b-2)/2}}
\frac{\left(\frac{-a+b+2s}{2}\right)_{(a+c-1)/2}} 
     {\left(\frac{2-a+b}{2}\right)_{(a+c-1)/2}}
\binom{m-1}{s-1}  \right. \\
 \hspace{3cm} \left.
\frac{(b+1)_{s-1}}{(b+c+1)_{s-1}}\frac{(c+1)_{n-1}}{(n-1)!} \frac{(b+c+n)_{m-n}}{(m-n)!}  \right].
\end{multline}

We are aiming to show that 
$\prod_{i=2}^{a} (b+c+1)_{i-1} / ((b+c+2)/2)_{\lceil (a-2)/2 \rceil}$ 
is a factor of 
\eqref{matrix}. Using \eqref{dreisum:1}, 
it remains to show that the 
following sum is equal to a polynomial in $b$ and $c$ for all integers 
$b$ and $c$ if we fix $a$: 
\begin{multline}
\label{pol}
\left( \prod_{i=2}^{a} (i-1)! \right) 
\left(\frac{b+c+2}{2}\right)_{\lceil (a-2)/2 \rceil}
\left( \frac{c-a+3}{2} \right)_{a-1} 
\left(\frac{b-a+2}{2} \right)_{a-1} \\
\times
\sum_{n=1}^{a} \sum_{m=1}^{a} \sum_{s=1}^{m}  \left[ (-1)^{n+s} 
\frac{\left(\frac{1-a+c+2n}{2}\right)_{(a+b-2)/2}}
     {\left(\frac{3-a+c}{2}\right)_{(a+b-2)/2}}
\frac{\left(\frac{-a+b+2s}{2}\right)_{(a+c-1)/2}} 
     {\left(\frac{2-a+b}{2}\right)_{(a+c-1)/2}}
\binom{m-1}{s-1}  \right. \\
\left.
\frac{(b+1)_{s-1}}{(b+c+1)_{s-1}}\frac{(c+1)_{n-1}}{(n-1)!} \frac{(b+c+n)_{m-n}}{(m-n)!}  \right]. 
\end{multline}
As it stands, this does not appear as a polynomial in $b$ and $c$, since $b$ and $c$ occur 
in the second arguments of some Pochhammer symbols. 
But by using 
\begin{equation}
\label{pochid}
\frac{(r)_{n}}{(s)_{n}}=\frac{(s+n)_{r-s}}{(s)_{r-s}},
\end{equation}
which is valid for all integers $n$, 
we are able to ban $b$ and $c$  from the second arguments of the 
Pochhammer symbols in \eqref{pol}. Then, after some cancellations, we 
obtain 
\begin{multline*}
\prod_{i=2}^{a} (i-1)!
       \sum_{n=1}^{a} \sum_{m=1}^{a} \sum_{s=1}^{m} (-1)^{n+s}
       \left( \frac{1-a+c+2n}{2} \right)_{a-n} 
       \left( \frac{b+c+1}{2}    \right)_{n-1}
       \left( \frac{-a+b+2s}{2}  \right)_{a-s}     \\
\times \binom{m-1}{s-1} 
              \frac{\left(\frac{b+c+2}{2}\right)_{\lceil (a-2)/2 \rceil} 
              \left( \frac{b+c+1}{2} \right)_{s-1}}
             {(b+c+1)_{s-1}}
         \frac{ (b+1)_{s-1} (c+1)_{n-1} (b+c+n)_{m-n} }
             { (n-1)! (m-n)! } 
\end{multline*}                 
for \eqref{pol}, and this is manifestly a polynomial in $b$ and $c$. 
Indeed, because of $s \le m \le a$, thanks to the binomial 
coefficient, the term $(b+c+1)_{s-1}$ cancels with the two 
terms on top of the fraction.

\smallskip
 
Therefore we have finally proved that the linear factors in \eqref{lf:odd}, 
respectively in \eqref{lf:even}, divide the 
determinant $\det(D_a(b,c))$ as a polynomial in $b$ and $c$, and we may now turn to 
the computation of the remaining irreducible polynomial.

\medskip

{\bf Step 3: Computation of the irreducible polynomial.} 
We emphasize once more that Step 1 and Step 2  show that 
\begin{equation}
\label{all}
det \left( D_a(b,c) \right) = ( \mbox{Linear factors in 
\eqref{lf:odd}, respectively \eqref{lf:even}} ) \times P_a(b,c), 
\end{equation}
where $P_a(b,c)$ is a certain polynomial in $b$ and $c$ if we fix $a$. 
In this final step we prove that $P_a(b,c)$ equals 
\eqref{ir:1}, respectively \eqref{ir:2}.

The computation of the irreducible polynomial, which will be denoted by 
$P_a(b,c)$, is done in the following way: We consider the polynomial $P_a(b,c)$ 
as a polynomial in $b$ over $\mathbb{Z}[c]$ and find that it has 
`nice' evaluations  at $b=-c-k$, $k$ odd and $1 \le k \le a$, see 
\eqref{auswert2} and \eqref{auswert1}. We are able to 
prove these `nice' evaluations with the help of Lemma~\ref{dreisum}. Furthermore, 
we will see that $P_a(b,c)$ is a polynomial in $b$ with a degree smaller or 
equal to $\lfloor (a-1)/2 \rfloor$. Therefore these `nice' evaluations are 
just enough so that we can compute $P_a(b,c)$ by using Lagrange's interpolation 
formula, see \eqref{poodd}, respectively \eqref{poeven}.

\medskip

First we will convince ourselves that the assertion about the degree of 
$P_a(b,c)$ in $b$ is true: 
The entry   in the $i$th row and the $j$th column of the matrix $D_a(b,c)$,
see \eqref{matrix}, is a polynomial in 
$b$ of degree $j+1$, since $H$ is a polynomial in $b$ of degree $2$. 
Hence, in the defining expansion of the determinant,  
each summand  has  degree $2+3+\dots+a=(a+2)(a-1)/2$ in $b$ and 
therefore the degree of the determinant itself is at most 
$(a+2)(a-1)/2$.   One easily checks that 
the product of the linear factors is a polynomial in $b$ with degree $a^2/2$ in 
case of $a$ is even and $(a+1)(a-1)/2$ otherwise. Because of that 
$(a+2)(a-1)/2-a^2/2=\lfloor (a-1)/2 \rfloor$  respectively 
$(a+2)(a-1)/2-(a+1)(a-1)/2=\lfloor (a-1)/2 \rfloor$ is an upper estimation
for the degree of $P_a(b,c)$ in $b$. Later we will see that this is in fact the exact degree. 

\smallskip

We claim that
\begin{multline}
\label{auswert2}
P_a(-c-k,c)= \left( \frac{1}{2}\right)^{a-1}
\left( \prod_{i=2}^{a} (i-1)! \right) \left( \frac{1}{2} \right)_{(a-1)/2}
\frac{\left(\frac{k-1}{2}\right)! (-1)^{(k-1)/2}}{\left(\frac{a-k+1}{2}\right)_{(k-1)/2}} \\
\times \left(\frac{c+k+1}{2}\right)_{(a-k)/2} \left(\frac{c+1}{2}\right)_{(k-1)/2}
\end{multline}
if $a$ is odd and $k=1,3,\dots,a$, and
\begin{multline}
\label{auswert1}
P_a(-c-k,c)= \left( \frac{1}{2}\right)^{a}
\left( \prod_{i=2}^{a} (i-1)! \right) \left( \frac{1}{2} \right)_{(a-2)/2}
\frac{\left(\frac{k-1}{2}\right)! (-1)^{(k-1)/2}}{\left(\frac{a-k}{2}\right)_{(k-1)/2}} \\
 \times \left(\frac{c+k+2}{2}\right)_{(a-k-1)/2} \left(\frac{c+2}{2}\right)_{(k-1)/2}
\end{multline}
if $a$ is even and $k=1,3,\dots,a-1$.      

Assuming the truth of the claim, we would have
\begin{multline}
\label{poodd}
P_a(b,c)=\left(\frac{1}{2}\right)^{a-1} \prod_{i=2}^{a} (i-1)!
\sum_{k=0}^{(a-1)/2} \left[ \left(\frac{c+1}{2}\right)_k 
    \left(\frac{1+b+c}{2}\right)_k 
      \left(\frac{c+2k+2}{2}\right)_{(a-2k-1)/2}  \right. \\ \left.
      \times \left(\frac{b+c+2k+3}{2}\right)_{(a-2k-1)/2}    
      \frac{(\frac{1}{2})_{(a-2k-1)/2}}{(1)_{(a-2k-1)/2}}  \right]
\end{multline}  
if $a$ is odd, and, respectively,    
\begin{multline}
\label{poeven}
P_a(b,c)=\left(\frac{1}{2}\right)^{a}  \prod_{i=2}^{a} (i-1)!
\sum_{k=0}^{(a-2)/2} \left[ \left(\frac{c+2}{2}\right)_k 
    \left(\frac{1+b+c}{2}\right)_k 
      \left(\frac{c+2k+3}{2}\right)_{(a-2k-2)/2}  \right. \\ \left.
      \times \left(\frac{b+c+2k+3}{2}\right)_{(a-2k-2)/2}    
      \frac{(\frac{1}{2})_{(a-2k-2)/2}}{(1)_{(a-2k-2)/2}}  \right]  
\end{multline} 
if $a$ is even, by Lagrange interpolation.

Thus, it remains to show \eqref{auswert2} and \eqref{auswert1}. By definition, $P_a(b,c)$ is 
the quotient of $\det (D_a(b,c))$ and the linear factors in
\eqref{lf:odd}, respectively in \eqref{lf:even} (see \eqref{all}). 
This, together with a tedious but straightforward  calculation 
shows that the claims in \eqref{auswert2} and \eqref{auswert1} 
conflate to the following 
\begin{multline} 
\label{zz}
\left. \frac{\det ( D_a(b,c) ) \displaystyle 
 \binom{\frac{b+c-1}{2}}{\frac{a+b-2}{2}}
 \binom{\frac{b+c-1}{2}}{\frac{a+c-1}{2}}}
{\displaystyle \left( \prod_{i=2}^{a}  (1+b+c)_{i-1} (i-1)!  \right)
 \left( \frac{c-a+3}{2} \right)_{a-1} 
 \left( \frac{b-a+2}{2} \right)_{a-1}}  \right|_{b=-c-k}
 =
(-1)^a\frac{(c+1)_{k-1}}{(k-1)!}.
\end{multline}
The reader should notice that we are not able to directly set $b=-c-k$ on the left-hand side of \eqref{zz}, 
because the denominator vanishes for $b=-c-k$.
As mentioned before we want to replace the determinant  on the left-hand sides 
of \eqref{zz} by the triple sum in Lemma~\ref{dreisum}. 
Still even after cancellations we are  not able to directly set $b=-c-k$ in this 
triple sum, since then the denominator of the summand becomes 
zero if $k > n$. 

In order to avoid an indefinite expression we act as follows: If we reexamine the proof of 
Lemma~\ref{dreisum}, we obtain the following generalization of Lemma~\ref{dreisum}:

\it Let $a$,$b$,$c$ be positive integers and $b'$,$x$,$x$ be integers. Then the determinant

\unitlength1cm
\begin{picture}(0,3.5)
\put(-0.5,1.5){\parbox{4cm}{
\begin{equation} 
\label{m:2}
-\det_{1 \le i, j \le a+1} \left(
\begin{array}{ccc|c}
      &                    & &   \\
      &\binom{b+c}{c-i+j}  & & \binom{b'-x+y}{y-i+1} \\
      &                    & &   \\ \hline
      &\binom{c+x-y-1}{x-j}& & 0 
\end{array}
\right)
\end{equation}}}
\put(6.2,1.8){\makebox(2,1)[t]{ \scriptsize $1 \le j \le a$}}
\put(8.6,1.8){\makebox(2,1)[t]{ \scriptsize $j=a+1$}}
\put(10.3,0.9){\makebox(2,1)[t]{ \scriptsize $1 \le i \le a$}}
\put(10.3,-0.1){\makebox(2,1)[t]{ \scriptsize $i=a+1$}}
\put(11.5,0.6){\makebox(2,1)[t]{.}}
\end{picture} 

equals 
\begin{multline}
\label{dreisum1}
\left( \prod_{i=1}^{a} \prod_{j=1}^b \prod_{k=1}^c \frac{i+j+k-1}{i+j+k-2}  \right)
\frac{(1)_c}{(b+1)_c}\\
\times
\sum_{n=1}^{a} \sum_{m=1}^{a} \sum_{s=1}^{m}  \left[ (-1)^{n+s} 
\binom{c+x-y+n-2}{x-1} \binom{b'-x+y+s-1}{b'-x+s-1} \frac{(b+1)_{s-1}}{(b+c+1)_{s-1}}  \right. \\
\left.
\binom{m-1}{s-1} \frac{(c+1)_{n-1}}{(n-1)!} \frac{(b+c+n)_{m-n}}{(m-n)!}  \right].
\end{multline}  \rm

We set $b=-c-k+\varepsilon$ and $b'=-c-k$ in \eqref{dreisum1}. In view of the fact 
that \eqref{m:2} is continuous in $b$, and if we compare the determinant in  
\eqref{matrix}  with the determinant in \eqref{matrix:allg:1} evaluated at $x=(a+b)/2$ and 
$y=(a+c-1)/2$, we obtain that the left-hand side of \eqref{zz}  
is equal to
\begin{multline}
\label{left}
\lim_{\varepsilon \to 0} \sum_{n=1}^{a} \sum_{m=1}^{a} \sum_{s=1}^{m}  \left[ (-1)^{n+s} 
\binom{\frac{-k+2n-3}{2}}{\frac{a-c-k-2}{2}} \binom{\frac{-k+2s-3}{2}}{\frac{a+c-1}{2}} 
\binom{m-1}{s-1}  \right. \\
\times \left.
\frac{(-c-k+\varepsilon+1)_{s-1}}{(-k+\varepsilon+1)_{s-1}}
\frac{(c+1)_{n-1}}{(n-1)!} \frac{(-k+\varepsilon+n)_{m-n}}{(m-n)!}  \right].
\end{multline}
Now it remains to show that the right-hand side of \eqref{zz} equals \eqref{left} for 
$k=1,3,\dots,2 \lfloor a-1/2 \rfloor.$

Therefore let us consider the triple sum \eqref{left}. First of all, we are allowed 
to extend the sum over $s$ to all positive integers, since $\binom{m-1}{s-1}=0$ if $s>m$. Using 
hypergeometric notation for this innermost sum, we get
\begin{multline*}
\lim_{\varepsilon \to 0} \sum_{n = 1}^{a}\sum_{m = 1}^{a}
     \left[ {{\left( -1 \right) }^{n+1}}\,
          {} _{3} F _{2} \!\left [ \begin{matrix} { 1 - m, 1 - c - k+\varepsilon,
           {1\over 2} - {k\over 2}}\\ { 1 - k+ \varepsilon, 1 - {a\over 2} -
           {c\over 2} - {k\over 2}}\end{matrix} ; {\displaystyle
           1}\right ] \,\right. \\ \left. {({ \textstyle 1 + c}) _{-1 + n} \,
          ({ \textstyle 1 - {a\over 2} - {c\over 2} - {k\over 2}})
           _{-{1\over 2} + {a\over 2} + {c\over 2}} \,
          ({ \textstyle {1\over 2} - {a\over 2} + {c\over 2} + n})
           _{-1 + {a\over 2} - {c\over 2} - {k\over 2}} \,
          ({ \textstyle -k+ \varepsilon + n}) _{m - n} }\over 
        {({ \textstyle 1}) _{-{1\over 2} + {a\over 2} + {c\over 2}} \,
          ({ \textstyle 1}) _{-1 + {a\over 2} - {c\over 2} - {k\over
           2}} \,({ \textstyle 1}) _{m - n} \,
          ({ \textstyle 1}) _{-1 + n} } \right] .
\end{multline*}  
We apply the following transformation formula due to Thomae \cite[(3.1.1)]{gasper} to this 
hypergeometric sum,  
\begin{equation}
\label{thom} 
{} _{3} F _{2} \!\left[
                       \begin{matrix} 
                          a,b,-n  \\
                          d,e
                      \end{matrix}; {\displaystyle 1} \right] =   
\frac{(-b+e)_{n}}{(e)_{n}}
{} _{3} F _{2} \!\left[
                       \begin{matrix} 
                          -n,b,-a+d  \\
                          d,1+b-e-n
                      \end{matrix}; {\displaystyle 1} \right]
\end{equation}
with $a=1-m$, $b=1-c-k+\varepsilon$, $n=(k-1)/2$, $d=1-k+\varepsilon$ and $e=1-a/2-c/2-k/2$.

Writing the resulting $_3 F _2$--series as a sum over $s$, after some simplifications and 
cancellations this gives
\begin{multline*}
\lim_{\varepsilon \to 0} \sum_{n = 1}^{a}\sum_{m = 1}^{a}
     \sum_{s = 0}^{\infty} \left[ {{\left( -1 \right) }^
              {{1\over 2} - {{3\,a}\over 2} + {c\over 2} + k + n + s}}
             \,\left( -{{a}\over 2} + {c\over 2} + {k\over 2} \right)
             \, \right. \\   \left.
             \times {{({ \textstyle 1 + c}) _{-1 + n} \,
            ({ \textstyle 1 - c - k+\varepsilon}) _{s} \,
            ({ \textstyle {1\over 2} - {a\over 2} + {c\over 2} + n})
             _{-1 + {a\over 2} - {c\over 2} - {k\over 2}} \,
            ({ \textstyle -k +\varepsilon + n}) _{m - n + s} }\over 
          {({ \textstyle 1}) _{m - n} \,({ \textstyle 1}) _{-1 + n} \,
            ({ \textstyle 1}) _{-{1\over 2} + {k\over 2} - s} \,
            ({ \textstyle 1}) _{s} \,
            ({ \textstyle 1}) _{{1\over 2} + {a\over 2} - {c\over 2} -
             k + s} \,({ \textstyle 1 - k+ \varepsilon}) _{s} }} \right].
\end{multline*} 

Next we interchange the two inner sums and reverse the order of 
summation in the innermost sum. We obtain
\begin{multline*}
\lim_{\varepsilon \to 0} \sum_{n = 1}^{a}\sum_{s = 0}^{\infty}
     \sum_{m = 0}^{-1 + a} \left[ {{\left( -1 \right) }^
             {{1\over 2} - {3 a\over 2} + {c\over 2} + k + n + s}}\,
           \left( -{{a}\over 2} + {c\over 2} + {k\over 2} \right) \, \right. \\
           \left. \times {{({ \textstyle 1 + c}) _{-1 + n} \,
           ({ \textstyle 1 - c - k+\varepsilon}) _{s} \,
           ({ \textstyle {1\over 2} - {a\over 2} + {c\over 2} + n})
            _{-1 + {a\over 2} - {c\over 2} - {k\over 2}} 
            \left( -k + \varepsilon +n \right)_{a-m-n+s} }\over 
         {({ \textstyle 1}) _{a - m - n} \,
           ({ \textstyle 1}) _{-1 + n} \,
           ({ \textstyle 1}) _{-{1\over 2} + {k\over 2} - s} \,
           ({ \textstyle 1}) _{s} \,
           ({ \textstyle 1}) _{{1\over 2} + {a\over 2} - {c\over 2} -
            k + s}  \,
           ({ \textstyle 1 - k + \varepsilon }) _{s} }} \right].
\end{multline*} 

The reason why reversing was so helpful, is that we now can extend the summation 
with respect to $m$ to all nonnegative integers (since $1/(1)_{a-m-n}=0$ if $m>a-n$).
If we use hypergeometric notation this yields
\begin{multline*}
\lim_{\varepsilon \to 0} \sum_{n = 1}^{a}\sum_{s = 0}^{\infty}
    \left[  {{\left( -1 \right) }^
           {{1\over 2} - {3 a\over 2} + {c\over 2} + k + n + s}}\,
         {} _{2} F _{1} \!\left [ \begin{matrix} { 1,-a+n}\\ { 1 -
          a + k + \varepsilon - s}\end{matrix} ; {\displaystyle 1}\right ] \,  \right. \\ \left.
        \times {{({ \textstyle 1 + c}) _{-1 + n} \,
         ({ \textstyle 1 - c - k+\varepsilon}) _{s} \,
         ({ \textstyle {{-a}\over 2} + {c\over 2} + {k\over 2}})
          \,({ \textstyle {1\over 2} - {a\over 2} + {c\over 2} + n})
          _{-1 + {a\over 2} - {c\over 2} - {k\over 2}} 
           \left( -k + \varepsilon +n \right)_{a-n+s}}\over 
       {({ \textstyle 1}) _{a - n} \,({ \textstyle 1}) _{-1 + n} \,
         ({ \textstyle 1}) _{-{1\over 2} + {k\over 2} - s} \,
         ({ \textstyle 1}) _{s} \,
         ({ \textstyle 1}) _{{1\over 2} + {a\over 2} - {c\over 2} - k
          + s} \,
         ({ \textstyle 1 - k + \varepsilon }) _{s} }} \right].                   
\end{multline*} 

We apply Vandermonde's summation formula \eqref{Vandermonde} to this 
hypergeometric sum. After some manipulations this gives
\begin{multline*}
\lim_{\varepsilon \to 0} \sum_{n = 1}^{a}\sum_{s = 0}^{\infty} \left[
     {{\left( -1 \right) }^{-{1\over 2} + {k\over 2} }}\,
         \left( -{{a}\over 2} + {c\over 2} + {k\over 2} \right) \, 
         \frac{({ \textstyle -a + k - \varepsilon  - s}) _{a - n} \,
         ({ \textstyle 1 + k - \varepsilon - n - s}) _{s}}
         {({ \textstyle 1}) _{{1\over 2} + {k\over 2} - n}
          ({ \textstyle 1}) _{-{1\over 2} + {k\over 2} - s}}  \, 
           \right. \\ 
         \left. \times {{({ \textstyle 1 + c}) _{-1 + n} \,
         ({ \textstyle 1 - c - k+\varepsilon}) _{s} \,
         ({ \textstyle {3\over 2} + {a\over 2} - {c\over 2} - k + s})
          _{-1 + k - n - s} }\over 
       {({ \textstyle 1}) _{a - n} \,
         ({ \textstyle 1}) _{-1 + n} \,
          \,({ \textstyle 1}) _{s} \,({ \textstyle 1-k + \varepsilon }) _{s} }} \right].         
\end{multline*} 
Now we are able to perform the limit $\varepsilon \to 0$, since, because of the factor 
$1/(1)_{-1/2+k/2-s}$, we have $k > s$. `Performing the limit' means that we simply 
set $\varepsilon=0$.

With pleasure we notice that our triple sum has simplified to a double sum.
Big surprise arises when I finally claim that the summand of the sum is only 
different from zero if and only if $n=k/2+1/2$ and $s=k/2-1/2$, and we therefore get 
rid of all sums: 
The factor $1/(1)_{1/2+k/2-n}$ implies $n \le 1/2+k/2$, and the factor 
$1/(1)_{-1/2+k/2-s}$ implies $s \le -1/2+k/2$. Next we notice that
\begin{equation*}
(-a+k-s)_{a-n} (1+k-n-s)_{s} = \frac{(-a+k-s)_{1+a-n-s}}{(k-n-s)}.
\end{equation*}
The numerator on the right-hand side is zero since $-a+k-s \le -a+k \le 0$ 
and $k-n \ge 0$. Therefore the summand is only different from zero if 
$k-n-s=0$. This together with our first observation gives my assertion.

It remains to compute the summand in the last double sum for 
$n=k/2+1/2$  and $s=k/2-1/2$ (and, of course, $\varepsilon=0$). After some simplifications 
one does indeed get the right-hand side of \eqref{zz} and Lemma~\ref{detlem} is finally proved.  

\bigskip

Now we turn to the case that the side lengths $a$,$b$,$c$ of the hexagon have the same parity and 
therefore to the evaluation of the determinant in Lemma~\ref{det} with $x=(a+b)/2$ and 
$y=(a+c)/2$ (see the beginning of this section).

\begin{lem} 
\label{detlem2}
Let $a$,$b$,$c$ be integers and $a \equiv c \hspace{2mm} (\text{\rm mod} \hspace{1mm} 2)$. Then the determinant 

\unitlength1cm
\begin{picture}(0,5)
\put(-1,2.5){\parbox{5cm}{
\begin{equation*} 
-\det_{1 \le i, j \le a+1} \left(
\begin{array}{ccc|c}
      &                    & &   \\
      &\displaystyle \binom{b+c}{c-i+j}  & & \displaystyle \binom{\frac{b+c}{2}}{\frac{c}{2}-i+\frac{a+2}{2}} \\
      &                    & &   \\ \hline
      &\displaystyle \binom{\frac{b+c-2}{2}}{\frac{c}{2}-\frac{a+2}{2}+j}& & 0 
\end{array}
\right)
\end{equation*}}}
\put(5.2,3.3){\makebox(2,1)[t]{ \scriptsize $1 \le j \le a$}}
\put(8.2,3.3){\makebox(2,1)[t]{ \scriptsize $j=a+1$}}
\put(11,2.1){\makebox(2,1)[t]{ \scriptsize $1 \le i \le a$}}
\put(11,0.8){\makebox(2,1)[t]{ \scriptsize $i=a+1$}}
\end{picture}

{\parindent0cm is equal to}
 
\begin{multline*}
\left( \prod_{i=1}^a \prod_{j=1}^b \prod_{k=1}^c \frac{i+j+k-1}{i+j+k-2}\right)
   \frac{(1)_c}{(b+1)_{c+a-1}}  
\binom{\frac{b+c-2}{2}}{\frac{b-1}{2}} \binom{\frac{a+b+c-1}{2}}{\frac{b-1}{2}} 2^{a-1}\\ 
\times \left( \left(\frac{c+1}{2} \right)_{(a-1)/2} \left( \frac{b+c+2}{2} \right)_{(a-1)/2} 
       \frac{\left( \frac{1}{2} \right)_{(a-1)/2}}
            {(1)_{(a-1)/2}}  \right. \\
+ \sum_{k=1}^{(a-1)/2}  \left(\frac{c+2}{2}\right)_{k-1} 
      \left(\frac{b+c}{2}\right)_{k}  
      \left(\frac{c+2k+1}{2}\right)_{(a-2k+1)/2} \\ 
       \left. \times  
      \left(\frac{b+c+2k+2}{2}\right)_{(a-2k-1)/2}
      \frac{(\frac{1}{2})_{(a-2k-1)/2}}{(1)_{(a-2k-1)/2}}  \right)
\end{multline*}
in case that $a$ is odd, and 
\begin{multline*}
\left( \prod_{i=1}^a \prod_{j=1}^b \prod_{k=1}^c \frac{i+j+k-1}{i+j+k-2}\right)
   \frac{ (1)_c}{(b)_{c+a}}       
\binom{\frac{b+c-2}{2}}{\frac{b-2}{2}} \binom{\frac{a+b+c-2}{2}}{\frac{b-2}{2}} 2^{a} \\
\times \left( \left(\frac{c+2}{2} \right)_{(a-2)/2} \left( \frac{b+c+2}{2} \right)_{a/2} 
       \frac{\left( \frac{1}{2} \right)_{a/2}}
            {(1)_{(a-2)/2}}  \right. \\
+ \sum_{k=1}^{a/2}  \left(\frac{c+1}{2}\right)_k           
      \left(\frac{b+c}{2}\right)_{k} 
      \left(\frac{c+2k+2}{2}\right)_{(a-2k)/2}  \\  
       \left. \times 
      \left(\frac{b+c+2k+2}{2}\right)_{(a-2k)/2}
      \frac{(\frac{1}{2})_{(a-2k)/2}}{(1)_{(a-2k)/2}} \right)
\end{multline*}
in case that $a$ is even.
\end{lem} 

As already mentioned at the beginning of this section the proofs of Lemma~\ref{detlem} and Lemma~\ref{detlem2} are 
quite similar. Therefore we restrict ourselves to just pointing out the differences to 
the proof of Lemma~\ref{detlem}, but we omit the details. The only major difference to 
Lemma~\ref{detlem} is that in the present case the irreducible polynomial $P_a(b,c)$ (see, e.g., \eqref{all}) 
has one `exceptional evaluation'  which has to be treated in a different way than the rest. 
By the way, this exceptional evaluation
also causes the formulas for the number of rhombus tilings which contain the `almost central' rhombus above the 
centre to be  not 
as simple compared to the 
formulas  for the number of rhombus tilings which contain  the central rhombus.

{\bf Proof of Lemma~\ref{detlem2}. }

{\bf Step 1: From our determinant to a determinant with polynomial entries.}

We can straightforwardly adopt Step 1 from the proof of Lemma~\ref{detlem}, where we have reduced the problem to the evaluation a polynomial determinant, since 
we have not even specialised the coordinates of the fixed rhombus $(x,y)$ there. 
In the following $\hat{D}_a(b,c)$ denotes the matrix underlying the 
determinant in \eqref{matrix} evaluated for $x=(a+b)/2$ and $y=(a+c)/2$.

Again I was led to a conjecture concerning the determinant $\det ( \hat{D}_a(b,c) )$
by computing it for small values of $a$:

\begin{multline}
\label{erg:polodd2} 
\det_{2 \le i ,j \le a} \left( \hat{D}_a(b,c) \right) =
\left( \prod_{i=2}^{a-1} \left(1+b+c\right)_{i-1}  \right)
\left(  \prod_{i=2}^{a} (i-1)! \right) \\
  \left(\frac{b+1}{2}\right)_{(a-1)/2}^2 
  \left(\frac{c+1}{2}\right)_{(a-1)/2} 
  \left(\frac{2+b+c}{2}\right)_{(a-1)/2} 2^{a-1} \\
\times \left( \left(\frac{c+1}{2} \right)_{(a-1)/2} \left( \frac{b+c+2}{2} \right)_{(a-1)/2} 
       \frac{\left( \frac{1}{2} \right)_{(a-1)/2}}
            {(1)_{(a-1)/2}}  \right. \\
+ \sum_{k=1}^{(a-1)/2}  \left(\frac{c+2}{2}\right)_{k-1} 
      \left(\frac{b+c}{2}\right)_{k}  
      \left(\frac{c+2k+1}{2}\right)_{(a-2k+1)/2} \\ 
     \displaystyle  \left. \times  
      \left(\frac{b+c+2k+2}{2}\right)_{(a-2k-1)/2}
      \frac{(\frac{1}{2})_{(a-2k-1)/2}}{(1)_{(a-2k-1)/2}}  \right)
\end{multline}                                  
if $a$ is odd, and
\begin{multline}
\label{erg:poleven2}
\det_{2 \le i ,j \le a} \left( \hat{D}_a(b,c) \right) 
=     
\left( \prod_{i=2}^{a-1} \left(1+b+c\right)_{i-1} \right)
 \left(  \prod_{i=2}^{a} (i-1)! \right) \\
    \left(\frac{b}{2}\right)_{a/2} \left(\frac{b+2}{2}\right)_{(a-2)/2} 
    \left(\frac{c+2}{2}\right)_{(a-2)/2} 
     \left(\frac{2+b+c}{2}\right)_{(a-2)/2}  2^{a-1} \\
   \times \left( \left(\frac{c+2}{2} \right)_{(a-2)/2} \left( \frac{b+c+2}{2} \right)_{a/2} 
       \frac{\left( \frac{1}{2} \right)_{a/2}}
            {(1)_{(a-2)/2}}  \right. \\
 + \sum_{k=1}^{a/2}  \left(\frac{c+1}{2}\right)_k           
      \left(\frac{b+c}{2}\right)_{k} 
      \left(\frac{c+2k+2}{2}\right)_{(a-2k)/2}  \\  
      \left. \times 
      \left(\frac{b+c+2k+2}{2}\right)_{(a-2k)/2}
      \frac{(\frac{1}{2})_{(a-2k)/2}}{(1)_{(a-2k)/2}} \right)  
\end{multline}
if a is even. Just as in the proof of Lemma~\ref{detlem} the two assertions 
\eqref{erg:polodd2} and \eqref{erg:poleven2} are equivalent to the assertion in 
Lemma~\ref{detlem2}.

The linear factors of 
$\det ( \hat{D}_a(b,c) )$ are 
\begin{equation}
\label{lf:odd2}
\prod_{i=2}^{a-1} \left(1+b+c\right)_{i-1}  
  \left(\frac{b+1}{2}\right)_{(a-1)/2}^2 
  \left(\frac{c+1}{2}\right)_{(a-1)/2} 
  \left(\frac{2+b+c}{2}\right)_{(a-1)/2} 2^{2a-2}
\end{equation}
if $a$ is odd, and, respectively,
\begin{equation} 
\label{lf:even2}
\prod_{i=2}^{a-1} \left(1+b+c\right)_{i-1}
    \left(\frac{b}{2}\right)_{a/2} \left(\frac{b+2}{2}\right)_{(a-2)/2} 
    \left(\frac{c+2}{2}\right)_{(a-2)/2} 
     \left(\frac{2+b+c}{2}\right)_{(a-2)/2}  2^{2a-3}
\end{equation} 
if $a$ is even. The irreducible polynomial, denoted by $\hat{P}_a(b,c)$, reads
\begin{multline}  
\label{ir:12}
\left( \prod_{i=2}^{a} (i-1)! \right) \left(\frac{1}{2} \right)^{a-1} \left( \left(\frac{c+1}{2} \right)_{(a-1)/2} \left( \frac{b+c+2}{2} \right)_{(a-1)/2} 
       \frac{\left( \frac{1}{2} \right)_{(a-1)/2}}
            {(1)_{(a-1)/2}}  \right. \\
+ \sum_{k=1}^{(a-1)/2}  \left(\frac{c+2}{2}\right)_{k-1} 
      \left(\frac{b+c}{2}\right)_{k}  
      \left(\frac{c+2k+1}{2}\right)_{(a-2k+1)/2} \\ 
     \displaystyle  \left. \times  
      \left(\frac{b+c+2k+2}{2}\right)_{(a-2k-1)/2}
      \frac{(\frac{1}{2})_{(a-2k-1)/2}}{(1)_{(a-2k-1)/2}}  \right)
\end{multline}
if $a$ is odd, and respectively,
\begin{multline}
\label{ir:22}
\left( \prod_{i=2}^{a} (i-1)! \right) \left( \frac{1}{2} \right)^{a-2}
\left( \left(\frac{c+2}{2} \right)_{(a-2)/2} \left( \frac{b+c+2}{2} \right)_{a/2} 
       \frac{\left( \frac{1}{2} \right)_{a/2}}
            {(1)_{(a-2)/2}}  \right. \\
\displaystyle + \sum_{k=1}^{a/2}  \left(\frac{c+1}{2}\right)_k           
      \left(\frac{b+c}{2}\right)_{k} 
      \left(\frac{c+2k+2}{2}\right)_{(a-2k)/2}  \\  
\displaystyle       \left. \times 
      \left(\frac{b+c+2k+2}{2}\right)_{(a-2k)/2}
      \frac{(\frac{1}{2})_{(a-2k)/2}}{(1)_{(a-2k)/2}} \right)
\end{multline}
if $a$ is even.

\medskip

{\bf Step 2: The linear factors in \eqref{lf:odd2}, respectively in \eqref{lf:even2},  divide 
$\det(\hat{D}_a(b,c))$ as a polynomial in $b$ and $c$.} 

Again we have four different types of linear factors:

\begin{enumerate}
\item Factors of the form $(c+k)$: If $a$ is odd, $k=1,3,\dots, a-2$. Otherwise  
$k=2,4,\dots, a-2$. 

\item Factors of the form $(b+k)^2$: If $a$ is odd, $k=1,3,\dots, a-2$. Otherwise 
$k=2,4,\dots, a-2$. The factor $b$ occurs once if $a$ even. 

\item Factors of the form $(b+c+k)^{a-k-1}$ for 
$k=1,3,\dots, 2\lceil(a-2)/2\rceil-1$.

\item Factors of the form $(b+c+k)^{a-k}$ for 
$k=2,4,\dots, 2\lfloor(a-2)/2\rfloor$.
\end{enumerate}

{\it re 1. --- The factors of the form $(c+k)$ divide $\det(\hat{D}_a(b,c))$:} 
As described in the analogous passage of the proof of 
Lemma~\ref{detlem}, each linear factor of  $\det(\hat{D}_a(b,c))$ corresponds to 
a linear combination of rows of a certain matrix. The linear combinations 
for the factors of the form $(c+k)$ are 
\begin{equation}
\label{c+i:Lk2}
\sum_{i=(a-k+2)/2}^{a-k+1} 
\frac{(-1)^{i-1} (b+i)_{a-k+1-i} \left(\frac{-a+k-2+2i}{2}\right)_{a-k+1-i}}
{(1)_{a-k+1-i} \left(\frac{b-a+2i-2}{2}\right)_{a-k+1-i}} D_a(b,-k)_{(i,j)}=0
\end{equation}
for $j=2,3, \dots,a$, $1 \le k \le a-1$ and  
$k  \equiv a \hspace{2mm} (\text{mod} \hspace{1mm} 2)$.
In order to prove that $(c+k)$ divides the determinant 
$\det ( \hat{D}_a(b,c) )$ we just have to prove 
the identity \eqref{c+i:Lk2}.

The comparison of the identity in \eqref{c+i:Lk} and the identity in \eqref{c+i:Lk2} shows 
that these  two
identities are quite similar. And in fact the proofs do not differ 
essentially from each other, either:  Again we are able to apply {\it Gosper's 
algorithm } \cite{gosper} for hypergeometric sums to the left-hand side of \eqref{c+i:Lk2}
and recognize that this is actually a telescoping sum.
 
The expression  analogous to \eqref{c+i:id} is 
\begin{equation}
\label{c+i:22}
\sum_{i=0}^{\infty} \frac{(-1)^j \left( 1 - \frac{a}{2} + \frac{k}{2} \right)_i 
\left( -a - b + k \right)_{i+j-2}}
                         { (1)_{-a+i+j} \left( 1 - \frac{a}{2} - \frac{b}{2} +k \right)_i}  
\left. H \right|_{x=(a+b)/2, y=(a+c)/2, c=-k, i \rightarrow a - k + 1 -i}=0,
\end{equation}
which is the identity in \eqref{c+i:Lk2} after reversing the summation order.
If $\hat{f}(i)$ denotes the summand in the previous  sum and $\hat{g}(i)$ denotes the following expression  
\begin{multline*}
\frac{(-1)^j \left( 1 -\frac{a}{2} + \frac{k}{2} \right)_i \left( -a-b + k \right)_{i+j-2}}
     { (1)_{-1-a+i+j} \left( 1 - \frac{a}{2} - \frac{b}{2} + k \right)_{i-1}} \\
\times \left( \frac{-2 - a - b  - b i + 2 j + b j - b k + i k - j k + k^2}{2} \right),
\end{multline*}
then we have
\begin{equation*}
\hat{f}(i)=\hat{g}(i+1)-\hat{g}(i).
\end{equation*}
By using this identity it is easy to compute the left-hand side of \eqref{c+i:22} 
and to show that it is equal to zero.

{\it re 2. --- The factors of the form $(b+k)$ divide $\det(\hat{D}_a(b,c))$:}

The two linearly independent linear combination for the factors of the form are
$(b+k)$ are 
\begin{multline} 
\label{b+k:Lk12}
\sum_{i=k+2}^{(a+k+2)/2}
 \frac{(-1)^{i-k-2} (c+a-i+2)_{i-k} \left(\frac{a+k-2i+4}{2}\right)_{i-k-2}}
{(1)_{i-k-1} \left(\frac{c+a-2i+4}{2}\right)_{i-k-2}^2} \\ 
\times p_{i-k-1}(c+a-k-i+2) \hat{D}_a(-k,c)_{(i,j)} =0
\end{multline}
if $0 \le k < a-2$ and $k \equiv a \hspace{2mm} (\text{mod} \hspace{1mm} 2)$, and
\begin{multline}
\label{b+k:Lk22}
\sum_{i=k+3}^{(a+k+2)/2} 
\frac{(-1)^{i-k} (c+a-i+2)_{i-k-1} \left(\frac{a+k-2i+4}{2}\right)_{i-k-1}}
{(1)_{i-k-1} \left(\frac{c+a-2i+4}{2}\right)_{i-k-2}^2}  \\ 
\times p_{i-k-2}(c+a-k-i+1) \hat{D}_a(-k,c)_{(i,j)}  =- \hat{D}_a(-k,c)_{(k+1,j)}
\end{multline}
if $0 < k < a-2$ and $k \equiv a \hspace{2mm} (\text{mod} \hspace{1mm} 2)$, 
where $p_n(c)$ is the same sequence of polynomials as in Lemma~\ref{detlem}, see 
\eqref{seq}. The first identity, \eqref{b+k:Lk12}, corresponds to \eqref{b+k:Lk1}, and
the second identity, \eqref{b+k:Lk22}, corresponds to \eqref{b+k:Lk2}. Again these two 
linear combinations do not cover the case that $k=a-2$.  But similar to the 
situation in Lemma~\ref{detlem} it can be verified directly that 
$(b+a-2)^2$ is a factor of the determinant $\det(\hat{D}_a(b,c))$.

The proofs of the identities  \eqref{b+k:Lk12} and  \eqref{b+k:Lk22} 
are analogous to the proofs of their corresponding identity: Again 
we first have to interchange the summation and then split the double sum 
into four smaller sums according to the polynomial $H$. Just as in Lemma~\ref{detlem} the  inner 
sums of the first and the third double sum  nearly cancel out each other, and so do the inner 
sums
of the second and the fourth summand. This reduces the problem to identities only involving $_2 F _1$--series and 
therefore Vandermonde's summation formula \eqref{Vandermonde} finishes the proof of 
these identities.

{\it re 3., 4. --- $\prod\limits_{i=2}^{a} (b+c+1)_{i-1}/ ((b+c+1)/2)_{\lceil (a-1)/2 \rceil}$ 
is a factor of $\det(\hat{D}_a(b,c))$:}
Just as in Lemma~\ref{detlem} the triple sum from Lemma~\ref{dreisum} provides an easy proof that
the factors of type $3$ and type $4$ divide the determinant $\det(\hat{D}_a(b,c))$ as a polynomial in $b$ and 
$c$. Again a combination of \eqref{factor1}, \eqref{factor2} and  
Lemma~\ref{dreisum} gives the equation  analogous to \eqref{dreisum:1}, namely 
\begin{multline}
\label{dreisum:12}
\det  (\hat{D}_a(b,c) ) 
=\left( \prod_{i=2}^{a} (b+c+1)_{i-1} (i-1)! \right) \left( \frac{c-a+2}{2} \right)_{a-1} 
\left(\frac{b-a+2}{2} \right)_{a-1} \\
 \quad \times
\sum_{n=1}^{a} \sum_{m=1}^{a} \sum_{s=1}^{m}  \left[ (-1)^{n+s} 
\frac{\left(\frac{-a+c+2n}{2}\right)_{(a+b-2)/2}}
     {\left(\frac{2-a+c}{2}\right)_{(a+b-2)/2}}
\frac{\left(\frac{-a+b+2s}{2}\right)_{(a+c)/2}} 
     {\left(\frac{2-a+b}{2}\right)_{(a+c)/2}}
\binom{m-1}{s-1}  \right. \\
 \hspace{3cm} \left.
\frac{(b+1)_{s-1}}{(b+c+1)_{s-1}}\frac{(c+1)_{n-1}}{(n-1)!} \frac{(b+c+n)_{m-n}}{(m-n)!}  \right].
\end{multline}    
By using \eqref{pochid} we obtain that the quotient of \eqref{dreisum:12} and the product 
of the linear factors $\prod\limits_{i=2}^{a} (b+c+1)_{i-1}/ ((b+c+1)/2)_{\lceil (a-1)/2 \rceil}$ is 
the following  polynomial in $b$ and $c$ if we fix $a$:
\begin{multline*}
\prod_{i=2}^{a} (i-1)!
       \sum_{n=1}^{a} \sum_{m=1}^{a} \sum_{s=1}^{m} (-1)^{n+s}
       \left( \frac{-a+c+2n}{2} \right)_{a-n} 
       \left( \frac{b+c}{2}    \right)_{n-1}
       \left( \frac{-a+b+2s}{2}  \right)_{a-s}     \\
\times \binom{m-1}{s-1} 
              \frac{\left(\frac{b+c+1}{2}\right)_{\lceil (a-1)/2 \rceil} 
              \left( \frac{b+c+2}{2} \right)_{s-1}}
             {(b+c+1)_{s-1}}
         \frac{ (b+1)_{s-1} (c+1)_{n-1} (b+c+n)_{m-n} }
             { (n-1)! (m-n)! }.
\end{multline*}

{\bf Step 3: Computation of the irreducible polynomial.} 
Analogously to the situation in Lemma~\ref{dreisum} we now evaluate the quotient 
of the linear factors in \eqref{lf:odd2}, respectively \eqref{lf:even2}, and the determinant
$\det ( \hat{D}_a(b,c) )$. This quotient, which will be denoted by $\hat{P}_a(b,c)$, 
is a polynomial in $b$ and $c$  if we fix $a$.
Again we will find that the polynomial  $\hat{P}_a(b,c)$ has enough `nice' evaluations so 
that we can compute the polynomial by using Lagrange's interpolation formula. 
Namely, these `nice' evaluations arise for $b=-c-k$, $0 \le k \le a$ and $k$ is even, 
and the degree of $\hat{P}_a(b,c)$ as a polynomial in $b$ happens to be 
$\lfloor a/2 \rfloor$. 

The assertion about the degree of $\hat{P}_a(b,c)$ as a polynomial in $b$ can 
be checked routinely, just as in Lemma~\ref{detlem}. Concerning 
the computation of `nice' evaluations of the polynomial  $\hat{P}_a(b,c)$ 
there arises the one and only major difference to Lemma~\ref{detlem}:
For $k=2,4, \dots, \lfloor a/2 \rfloor$ the evaluations 
$\hat{P}_a(-c-k,c)$ can be proved in the same way as in Lemma~\ref{detlem}, but 
for $k=0$ the situation is different, and therefore this exceptional evaluation needs a separate proof.  
This fact already shows up  if we 
look at the following conjecture for the `nice' evaluations:

We claim that
\begin{equation}
\label{aus1}
\hat{P}_a(-c,c)=\left( \frac{1}{2} \right)^{a-1} \left( \prod_{i=2}^{a} (i-1)!  \right) \left( \frac{1}{2} \right)_{(a-1)/2} 
          \left( \frac{c+1}{2} \right)_{(a-1)/2} 
\end{equation}
if $a$ is odd, and 
\begin{equation}
\label{aus2}
\hat{P}_a(-c,c)=\left( \frac{1}{2} \right)^{a-2} \left( \prod_{i=2}^{a} (i-1)!  \right)
\left( \frac{a}{2} \right)  \left( \frac{1}{2} \right)_{a/2} 
          \left( \frac{c+2}{2} \right)_{(a-2)/2}
\end{equation} 
if $a$ is even, and, furthermore,   
\begin{multline}
\label{auswert22}
\hat{P}_a(-c-k,c)= \left( \frac{1}{2}\right)^{a-1}
\left( \prod_{i=2}^{a} (i-1)! \right) \left( \frac{1}{2} \right)_{(a-1)/2}
\frac{\left(\frac{k}{2}\right)! (-1)^{k/2}}{\left(\frac{a-k}{2}\right)_{k/2}} \\
\times \left(\frac{c+k+1}{2}\right)_{(a-k+1)/2} \left(\frac{c+2}{2}\right)_{(k-2)/2}
\end{multline}
if $a$ is odd and $k=2,4,\dots,a-1$, and
\begin{multline}
\label{auswert12}
\hat{P}_a(-c-k,c)= \left( \frac{1}{2}\right)^{a-2}
\prod_{i=1}^{a-1} i! \left( \frac{1}{2} \right)_{(a-2)/2}
\frac{\left(\frac{k}{2}\right)! (-1)^{k/2}}{\left(\frac{a-k+1}{2}\right)_{(k-2)/2}} \\
 \times \left(\frac{c+k+2}{2}\right)_{(a-k)/2} \left(\frac{c+1}{2}\right)_{k/2}
\end{multline}
if $a$ is even and $k=2,4,\dots,a$.

Assuming the truth of the claim, we would have
\begin{multline*}
\hat{P}_a(b,c)=\left(\frac{1}{2}\right)^{a-1} \left( \prod_{i=2}^{a} (i-1)!   \right)
\left( \left(\frac{c+1}{2} \right)_{(a-1)/2} \left( \frac{b+c+2}{2} \right)_{(a-1)/2} 
       \frac{\left( \frac{1}{2} \right)_{(a-1)/2}}
            {(1)_{(a-1)/2}}  \right. \\
+ \sum_{k=1}^{(a-1)/2}  \left(\frac{c+2}{2}\right)_{k-1} 
      \left(\frac{b+c}{2}\right)_{k}  
      \left(\frac{c+2k+1}{2}\right)_{(a-2k+1)/2} \\ 
       \left. \times  
      \left(\frac{b+c+2k+2}{2}\right)_{(a-2k-1)/2}
      \frac{(\frac{1}{2})_{(a-2k-1)/2}}{(1)_{(a-2k-1)/2}}  \right)
\end{multline*}  
if $a$ is odd, and, respectively    
\begin{multline*}
\hat{P}_a(b,c)=\left(\frac{1}{2}\right)^{a-2}  \prod_{i=2}^{a} (i-1)!
\left( \left(\frac{c+2}{2} \right)_{(a-2)/2} \left( \frac{b+c+2}{2} \right)_{a/2} 
       \frac{\left( \frac{1}{2} \right)_{a/2}}
            {(1)_{(a-2)/2}}  \right. \\
\displaystyle + \sum_{k=1}^{a/2}  \left(\frac{c+1}{2}\right)_k           
      \left(\frac{b+c}{2}\right)_{k} 
      \left(\frac{c+2k+2}{2}\right)_{(a-2k)/2}  \\  
\displaystyle       \left. \times 
      \left(\frac{b+c+2k+2}{2}\right)_{(a-2k)/2}
      \frac{(\frac{1}{2})_{(a-2k)/2}}{(1)_{(a-2k)/2}} \right)
\end{multline*} 
if $a$ is even, by Lagrange interpolation.      

Thus, it remains to show \eqref{aus1}, \eqref{aus2}, \eqref{auswert12} and \eqref{auswert22}. In these claims 
we replace the polynomial $\hat{P}_a(b,c)$ by the quotient of $\det ( \hat{D}_a(b,c) )$
and the linear factors in  
\eqref{lf:odd2}, respectively \eqref{lf:even2}, in order to see that the four claims conflate to the following two:
\begin{multline} 
\label{aus}
\left.  \frac{\displaystyle 2 \det ( \hat{D}_a(b,c) )
 \binom{\frac{b+c-2}{2}}{\frac{a+b-2}{2}}
 \binom{\frac{b+c}{2}}{\frac{a+c}{2}}}
{\displaystyle \left( \prod_{i=2}^{a}  (1+b+c)_{i-1} (i-1)!  \right)
 \left( \frac{c-a+2}{2} \right)_{a-1} 
 \left( \frac{b-a+2}{2} \right)_{a-1} (b+c) }  \right|_{b=-c} \\
={{{-a + {{\left( -1 \right) }^{a+1}}\,a + c + 
     {{\left( -1 \right) }^{a+1}}\,c}}
     \over {c\,\left( a + c \right) }},
\end{multline}
and 
\begin{multline} 
\label{zz2}
\left. \frac{\displaystyle \det ( \hat{D}_a(b,c) )
 \binom{\frac{b+c-2}{2}}{\frac{a+b-2}{2}}
 \binom{\frac{b+c}{2}}{\frac{a+c}{2}}}
{\displaystyle \left( \prod_{i=2}^{a}  (1+b+c)_{i-1} (i-1)!  \right)
 \left( \frac{c-a+2}{2} \right)_{a-1} 
 \left( \frac{b-a+2}{2} \right)_{a-1}}  \right|_{b=-c-k}
=(-1)^{a-1} \frac{(c+1)_{k-1}}{(k-1)!}
\end{multline}
for $k=2,4,\dots , \lfloor a/2 \rfloor$.
Again the reader should know that we are not able to directly set $b=-c-k$ on the left hand 
sides of \eqref{zz2} and \eqref{aus}, because the denominator vanishes for $b=-c-k$.

First we consider the case that $k \not= 0$. As already mentioned, the situation in this case is 
quite the same as in Lemma~\ref{detlem}. Therefore I only explain the essential steps and omit details.

For the left-hand side of \eqref{zz2}  we use again the 
modification of Lemma~\ref{dreisum}, see \eqref{m:2} and \eqref{dreisum1}. 
Thus, this left-hand side is equal to 
\begin{multline}
\label{left2}
\lim_{\varepsilon \to 0} \sum_{n=1}^{a} \sum_{m=1}^{a} \sum_{s=1}^{m}  \left[ (-1)^{n+s} 
\binom{\frac{-k+2n-4}{2}}{\frac{a-c-k-2}{2}} \binom{\frac{-k+2s-2}{2}}{\frac{a+c}{2}} 
\binom{m-1}{s-1}  \right. \\
\times \left.
\frac{(-c-k+\varepsilon+1)_{s-1}}{(-k+\varepsilon+1)_{s-1}}
\frac{(c+1)_{n-1}}{(n-1)!} \frac{(-k+\varepsilon+n)_{m-n}}{(m-n)!}  \right].
\end{multline}
Hence, it remains to show that the right hand side of \eqref{zz2} is equal to \eqref{left2} for 
$k=2,4, \dots, \lfloor a/2 \rfloor$. 

Just as in Lemma~\ref{dreisum} we apply the transformation formula due 
to Thomae, see \eqref{thom}, to the innermost sum of the triple sum, with slightly changed 
parameters, i.e., 
$a=1-m$, $b=1-c-k+\varepsilon$, $n=(k-2)/2$, $d=1-k+\varepsilon$ and $e=1-a/2-c/2-k/2$, 
compared to the proof of Lemma~\ref{dreisum}.
This yields
\begin{multline*}
\lim_{\varepsilon \to 0} \sum_{n = 1}^{a}\sum_{m = 1}^{a}
     \sum_{s = 0}^{\infty}{{\left( -1 \right) }^
             {{a\over 2} + {c\over 2} + n + s}}\,
           \left( {{-a}\over 2} + {c\over 2} + {k\over 2} \right) \,    \\
          \times {{({ \textstyle 1 + c}) _{-1 + n} \,
           ({ \textstyle 1 - c - k}) _{s} \,
           ({ \textstyle {{-a}\over 2} + {c\over 2} + n}) _{-1 +
            {a\over 2} - {c\over 2} - {k\over 2}} \,
           ({ \textstyle \varepsilon - k + n}) _{m - n + s} }\over 
         {({ \textstyle 1}) _{m - n} \,({ \textstyle 1}) _{-1 + n} \,
           ({ \textstyle 1}) _{-1 + {k\over 2} - s} \,
           ({ \textstyle 1}) _{s} \,
           ({ \textstyle 1}) _{1 + {a\over 2} - {c\over 2} - k + s} \,
           ({ \textstyle 1 + \varepsilon - k}) _{s} }}.
\end{multline*} 
Next we interchange the two inner sums, reverse the order of the summation in the 
new innermost sum and use hypergeometric notation for this innermost sum. 
We obtain
\begin{multline*}
\lim_{\varepsilon \to 0} \sum_{n=1}^{m} \sum_{s=0}^{\infty} {{\left( -1 \right) }^{{a\over 2} + {c\over 2} - n + s}}\,
     {} _{2} F _{1} \!\left [ \begin{matrix} { 1, -a + n}\\ { 1 - a -
      \varepsilon + k - s}\end{matrix} ; {\displaystyle 1}\right ] \, \\
     \times {{({ \textstyle 1 + c}) _{-1 + n} \,
     ({ \textstyle 1 - c - k}) _{s} \,
     ({ \textstyle {{-a}\over 2} + {c\over 2} + {k\over 2}}) \,
     ({ \textstyle {{-a}\over 2} + {c\over 2} + n}) _{-1 + {a\over 2}
      - {c\over 2} - {k\over 2}} \,
     ({ \textstyle \varepsilon - k + n}) _{a - n + s} }\over 
   {({ \textstyle 1}) _{a - n} \,({ \textstyle 1}) _{-1 + n} \,
     ({ \textstyle 1}) _{-1 + {k\over 2} - s} \,
     ({ \textstyle 1}) _{s} \,
     ({ \textstyle 1}) _{1 + {a\over 2} - {c\over 2} - k + s} \,
     ({ \textstyle 1 + \varepsilon - k}) _{s} }}.                 
\end{multline*} 
Then we apply Vandermonde's summation formula \eqref{Vandermonde} to this hypergeometric sum. 
After some  manipulations this gives
\begin{multline*}
\lim_{\varepsilon \to 0} \sum_{n = 1}^{a}\sum_{s = 0}^{\infty}
     {{\left( -1 \right) }^{-1 - {k\over 2}}}\,
         \left( {{-a}\over 2} + {c\over 2} + {k\over 2} \right) \,
         ({ \textstyle 1 + c}) _{-1 + n} \,
         ({ \textstyle 1 - c - k}) _{s} \,\\
         \times {{({ \textstyle -a - \varepsilon + k - s}) _{a - n} \,
         ({ \textstyle 1 - \varepsilon + k - n - s}) _{s} \,
         ({ \textstyle 2 + {a\over 2} - {c\over 2} - k + s}) _{-1 + k
          - n - s} }\over 
       {({ \textstyle 1}) _{a - n} \,
         ({ \textstyle 1}) _{1 + {k\over 2} - n} \,
         ({ \textstyle 1}) _{-1 + n} \,
         ({ \textstyle 1}) _{-1 + {k\over 2} - s} \,
         ({ \textstyle 1}) _{s} \,({ \textstyle 1 + \varepsilon - k}) _{s} }}.         
\end{multline*} 
Now we are able to perform the limit $\varepsilon \to 0$ and 
then observe that the summand of the double sum is only different from zero if 
$n=k/2+1$ and $s=k/2-1$. The evaluation of the summand of 
the double sum for these special values of $n$ and $s$ finishes the proof of 
\eqref{zz2} for $k=2,4, \dots \lfloor a / 2 \rfloor$.

Finally we consider the case that $k=0$, see \eqref{aus}. 
By Lemma~\ref{dreisum}, \eqref{factor1} and \eqref{factor2} 
the left-hand side of \eqref{aus} is equal to 
\begin{multline*}
\frac{2}{b+c} \sum_{n = 1}^{a}\sum_{m = 1}^{a}
     \sum_{s = 1}^{m}{{{{\left( -1 \right) }^{n + s}}\,
           ({ \textstyle 1 + b}) _{-1 + s} \,  ({ \textstyle 1 + c}) _{-1 + n} \,
           \over  {({ \textstyle 1}) _{-1 + s} \,
                  ({ \textstyle 1}) _{-1 + n} \,  }}} \\
         \times  {{({ \textstyle 1 + {a\over 2} + {c\over 2}}) _{-1 - {a\over
            2} + {b\over 2} + s} \,
           ({ \textstyle {{-a}\over 2} + {c\over 2} + n}) _{-1 +
            {a\over 2} + {b\over 2}} \,
           ({ \textstyle b + c + n}) _{m - n} \,
           ({ \textstyle 1 + m - s}) _{-1 + s} }\over 
         {({ \textstyle 1}) _{-1 + {a\over 2} + {b\over 2}} \,
           ({ \textstyle 1}) _{m - n} \,
           ({ \textstyle 1}) _{-1 - {a\over 2} + {b\over 2} + s} \,
           ({ \textstyle 1 + b + c}) _{-1 + s} }}.
\end{multline*}
before evaluating it at $b=-c$.
We take the factor $(b+c)/2$ out of the Pochhammer symbol 
\begin{equation*}
\left( \frac{2+a+c}{2} \right) _{(-2-a+b+2s)/2} = 
\left( \frac{2+a+c}{2} \right)_{(-2+a+b)/2} 
\left( \frac{b}{2} + \frac{c}{2} \right) 
\left( \frac{b+c+2}{2} \right)_{s-1}
\end{equation*}
and set $b=-c$. After some cancellations and simplifications 
concerning the summand of this triple sum, this yields 
\begin{multline}
\sum_{n = 1}^{a}\sum_{m = 1}^{a}
     \sum_{s = 1}^{m}{{{{\left( -1 \right) }^
             {-1 + {a\over 2} - {c\over 2} + n + s}}\,
           ({ \textstyle 1 - c}) _{-1 + s} \,
           ({ \textstyle 1 + c}) _{-1 + n} \, \over
            {({ \textstyle 1}) _{-1 + n} \,
           ({ \textstyle 1}) _{-1 + s} \,  }}} \\
           \times {{({ \textstyle {a\over 2} - {c\over 2}}) _{1 - n} \,
           ({ \textstyle 1 + {a\over 2} + {c\over 2}}) _{-1 - {a\over
            2} - {c\over 2}} \,
           ({ \textstyle n}) _{m - n} \,
           ({ \textstyle 1 + m - s}) _{-1 + s} }\over 
         {({ \textstyle 1}) _{1 - n} \,({ \textstyle 1}) _{m - n} \,
           ({ \textstyle 1}) _{-1 - {a\over 2} - {c\over 2} + s} }}.
\end{multline} 
The factor $1/(1)_{1-n}$ shows that the summand of the triple sum is 
only different from  zero if $n=1$. Therefore we get rid of the sum 
over $n$. We use hypergeometric notation for the inner
sum and obtain
\begin{equation*}
\sum_{m = 1}^{a}{{{{\left( -1 \right) }^{1 + c}}\,
       {} _{2} F _{1} \!\left [ \begin{matrix} { 1 - c, 1 - m}\\ { 1 -
        {a\over 2} - {c\over 2}}\end{matrix} ; {\displaystyle 1}\right
        ] }\over {{a\over 2} + {c\over 2}}}.
\end{equation*} 
Next we apply Vandermonde's summation formula \eqref{Vandermonde} to the $_2 F _1$--series.
This yields 
\begin{equation}
\label{ss}
\sum_{m = 1}^{a}{{{{\left( -1 \right) }^{c}}\,
       ({ \textstyle {-{a}\over 2} + {c\over 2}}) _{-1 + m} }\over 
     {({ \textstyle {-{a}\over 2} - {c\over 2}}) _{m} }}.
\end{equation}  
Once more,  Gosper's algorithm \cite{gosper} can be used to see that this single sum is actually 
telescoping. Namely, there holds
\begin{multline*}
(-1)^{c} \frac{\left( -\frac{a}{2}+\frac{c}{2} \right)_{-1+m}}{\left(-\frac{a}{2}-\frac{c}{2}\right)_m}\\
=\frac{(-1)^{c-1}}{\left(1 + \frac{a}{2} + \frac{c}{2} \right) c } 
\left( \frac{ \left( -\frac{a}{2} - \frac{c}{2} + m \right) \left( -1 - \frac{a}{2} + \frac{c}{2} \right)_{m+1}}
            { \left( -\frac{a}{2} - \frac{c}{2} \right)_{m+1}} -
       \frac{ \left( -\frac{a}{2} - \frac{c}{2} + m-1 \right) \left( -1 - \frac{a}{2} + \frac{c}{2} \right)_{m}}
            { \left( -\frac{a}{2} - \frac{c}{2} \right)_{m}} \right).
\end{multline*}
Because of that,  and because $a$ and $c$ have the same parity, the single sum \eqref{ss} is equal to 
\begin{equation*}
{{{-a + {{\left( -1 \right) }^{a+1}}\,a + c + 
     {{\left( -1 \right) }^{a+1}}\,c}}
     \over {c\,\left( a + c \right) }}.
\end{equation*}
If we compare this to the right-hand side of \eqref{aus} 
we see that this is just what we claimed. This finally completes  the proof of Lemma~\ref{detlem2}.

\section{The Proofs of Theorem~\ref{asym:th} and Theorem~\ref{asym:th2} }
\label{asymsec}

We end this article with the proofs of Theorem~\ref{asym:th} and Theorem~\ref{asym:th2}.
We start with  Theorem~\ref{asym:th}.
In Theorem~\ref{th:hp} we showed that the probability to choose a rhombus tilings of a 
hexagon with side lengths $a$,$b$,$c$,$a$,$b$,$c$ which contains 
a rhombus in the centre is, using hypergeometric notation,  
\begin{multline}
\label{hyp:odd}   
\,{} _{4} F _{3} \!\left [ \matrix { {1\over 2} -
      {a\over 2}, {1\over 2} + {c\over 2}, {1\over 2} + {b\over 2} +
      {c\over 2}, 1}\\ { 1 - {a\over 2}, 1 + {c\over 2}, {3\over 2} +
      {b\over 2} + {c\over 2}}\endmatrix ; {\displaystyle 1}\right ]  \\
     \times  {{ \,({ \textstyle {1\over 2}}) _{-{1\over 2} + {a\over 2}} \,
     ({ \textstyle 1}) _{c} \,
     ({ \textstyle 1 + {c\over 2}}) _{-{1\over 2} + {b\over 2}} \,
     ({ \textstyle 1 + {c\over 2}}) _{-1 + {a\over 2} + {b\over 2}} \,
     ({ \textstyle {3\over 2} + {b\over 2} + {c\over 2}}) _{-{1\over
      2} + {a\over 2}} {2^{a-1}}}\over 
   {({ \textstyle 1}) _{-{1\over 2} + {a\over 2}} \,
     {{({ \textstyle 1}) _{-{1\over 2} + {b\over 2}} }^2}\,
     ({ \textstyle 1 + b}) _{-1 + a + c} }}
\end{multline} 
in case that $a$ is odd, and 
\begin{multline}
\label{hyp:even}
\,{} _{4} F _{3} \!\left [ \matrix { 1 - {a\over 2}, 1 +
      {c\over 2}, {1\over 2} + {b\over 2} + {c\over 2}, 1}\\ { {3\over
      2} - {a\over 2}, {3\over 2} + {c\over 2}, {3\over 2} + {b\over
      2} + {c\over 2}}\endmatrix ; {\displaystyle 1}\right ] \,  \\
    \times  {{({ \textstyle {1\over 2}}) _{-1 + {a\over 2}} \,
     ({ \textstyle 1}) _{c} \,({ \textstyle b}) \,
     ({ \textstyle {1\over 2} + {c\over 2}}) _{{b\over 2}} \,
     ({ \textstyle {3\over 2} + {c\over 2}}) _{-1 + {a\over 2} +
      {b\over 2}} \,({ \textstyle {3\over 2} + {b\over 2} + {c\over
      2}}) _{-1 + {a\over 2}} {2^{a-2}}} \over 
    {({ \textstyle 1}) _{-1 + {a\over 2}} \,
     {{({ \textstyle 1}) _{{b\over 2}} }^2}\,
     ({ \textstyle 1 + b}) _{-1 + a + c} }}
\end{multline} 
in case that $a$ is even.
In order to see this, divide \eqref{number:odd}, respectively \eqref{number:even}, by 
MacMahon's formula for the total number of rhombus tilings of a hexagon with side 
lengths $a$,$b$,$c$,$a$,$b$,$c$ given in \eqref{mac}.

Next we apply Bailey's transformation formula (see \cite[(4.3.5.1)]{slater}) between 
two balanced $_4 F _3$--series,
\begin{multline}
\label{bai}
_{4} F _{3} \!\left [ \matrix { a, b, c, -n}\\ 
                              { e, f, 1+a+b+c-e-f-n}\endmatrix ; 
                              {\displaystyle 1}\right ] =
\frac{(e-a)_n (f-a)_n}{(e)_n (f)_n} \times \\
_{4} F _{3} \!\left [ \matrix { -n, a, 1+a+c-e-f-n, 1+a+b-e-f-n}\\ 
                              { 1+a+b+c-e-f-n,1+a-e-n,1+a-f-n}\endmatrix ; 
                              {\displaystyle 1}\right ], 
\end{multline}
which is valid if $n$ is a positive integer,  
with $a=1,b=(1+c)/2,c=(1+b+c)/2,n=(a-1)/2,e=(2-a)/2,f=(c+2)/2$ in case of 
$a$ is odd, and with $a=1,b=(2+c)/2,c=(1+b+c)/2,n=(a-2)/2,e=(3-a)/2,f=(c+3)/2$ in 
case of $a$ is even. This gives
\begin{multline*}
_{4} F _{3} \!\left [ \matrix { 1, 1, \frac{1}{2}-\frac{a}{2}, 1+\frac{b}{2}}\\ 
                              { \frac{3}{2}, 
                                \frac{3}{2}+\frac{b}{2}+\frac{c}{2}, 
                                \frac{3}{2}-\frac{a}{2}-\frac{c}{2}}\endmatrix ; 
                              {\displaystyle 1}\right ] \\
\times
\frac{a!b!c! ((a+b+c-2)/2)! ((a+b+c)/2)!}
     {((a-1)/2)!^2((b-1)/2)!^2 (c/2)! ((c-2)/2)! (a+b+c-1)! ((-1+a+c)/2)((1+b+c)/2)}
\end{multline*}
if $a$ is odd, and
\begin{multline*}
_{4} F _{3} \!\left [ \matrix { 1, 1, 1-\frac{a}{2}, \frac{1}{2}+\frac{b}{2}}\\ 
                              { \frac{3}{2}, \frac{3}{2}+\frac{b}{2}+\frac{c}{2}, 
                                \frac{3}{2}-\frac{a}{2}-\frac{c}{2}}\endmatrix ; 
                              {\displaystyle 1}\right ] \\
\times
\frac{(a-1)!b!c! ((a+b+c-1)/2)!^2}
     {((a-2)/2)!^2(b/2)!^2 ((c-1)/2)!^2 (a+b+c-1)! ((-1+a+c)/2)((1+b+c)/2)}
\end{multline*}
if $a$ is even, after transforming the Pochhammer symbol $(a)_{n}=(a+n-1)!/(a-1)!$ to factorials. 

Now we substitute $a \sim \alpha N$, $b \sim \beta N$ and $c \sim \gamma N$ and perform 
the limit $N \rightarrow \infty$. We use Stirling's formula to determine the limit for 
the quotient of the factorials in the second line as 
$\sqrt{\alpha}\sqrt{\beta}\sqrt{\gamma}
\sqrt{\alpha+\beta+\gamma}/(2 \pi (\alpha + \gamma) (\beta + \gamma))$. For the 
$_{4} F _{3}$--series, we may exchange limit and summation by uniform convergence:
\begin{equation*}
\lim_{N \rightarrow \infty} \!
_{4} F _{3} \!\left [ \matrix { 1, 1, \frac{1}{2}-\frac{a}{2}, 1+\frac{b}{2}}\\ 
                              { \frac{3}{2}, \frac{3}{2}+\frac{b}{2}+\frac{c}{2}, 
                                \frac{3}{2}-\frac{a}{2}-\frac{c}{2}}\endmatrix ; 
                              {\displaystyle 1}\right ] =
{_{2} F _{1}} \!\left[ \matrix {1,1}\\
                             {\frac{3}{2}} \endmatrix;
                             {\displaystyle 
                             \frac{\alpha \beta}{(\beta + \gamma)(\alpha + \gamma)}} \right].
\end{equation*}    
A combination of  these results and the use  of the identity 
\begin{equation*}
_{2} F _{1} \!\left[ \matrix {1,1}\\
                             {\frac{3}{2}} \endmatrix;
                             {\displaystyle z} \right]=
\frac{\arcsin \sqrt{z}}{\sqrt{z (1-z)}}
\end{equation*}
establishes Theorem~\ref{asym:th}.

\medskip

The proof of  Theorem~\ref{asym:th2} is analogous: Again we divide the formulas in Theorem~\ref{th:hp2},
see  \eqref{number:odd2} and \eqref{number:even},
by MacMahon's formula for the total number of rhombus tilings of a hexagon with side lengths $a$, $b$, $c$, $a$, $b$, $c$
given in \eqref{mac}, and obtain 
the probability to choose a rhombus tiling of a hexagon with side 
lengths $a$, $b$, $c$, $a$, $b$, $c$ which contain the `almost central' rhombus above the centre. 
 Next we apply Bailey's transformation formula \eqref{bai} on the hypergeometric sum and transform the 
Pochhammer symbols to factorials. We then substitute  
$a \sim \alpha N$, $b \sim \beta N$ and $c \sim \gamma N$
and perform the limit $N \to \infty$ in the same way as before.   
The additional summands (caused by the `exceptional' evaluation of the irreducible polynomial 
$\hat{P}_a(b,c)$) that appear in the formulas in Theorem~\ref{th:hp2} (when compared  
to the formulas in Theorem~\ref{th:hp}) vanish when this limit is performed.

\end{document}